\documentclass[11pt,leqno]{amsart}
\usepackage{epsfig}
\usepackage{amssymb}
\usepackage{amscd}
\usepackage[matrix,arrow]{xy}

\newtheorem{theo}{Theorem}[section]
\newtheorem{lemma}[theo]{Lemma}
\newtheorem{defi}[theo]{Definition}
\newtheorem{prop}[theo]{Proposition}

\newtheorem{cor}[theo]{Corollary}
\newtheorem{remark}[theo]{Remark}

\newtheorem{example}[theo]{Example}
\numberwithin{equation}{section}

\def\P{{\mathcal{P}}}
\def\D{{\mathcal{D}}}

\def\bR{{\mathbf R}}
\def\bL{{\mathbf L}}

\def\pre-tr{\operatorname{pre-tr}}
\def\h{\operatorname{h}}
\def\Hom{\operatorname{Hom}}
\def\End{\operatorname{End}}
\def\gr{\operatorname{gr}}

\newcommand{\bbZ}{{\mathbb Z}}

\newcommand{\cJ}{{\mathcal J}}
\newcommand{\cQ}{{\mathcal Q}}
\newcommand{\cF}{{\mathcal F}}
\newcommand{\cG}{{\mathcal G}}
\newcommand{\cO}{{\mathcal O}}
\newcommand{\cP}{{\mathcal P}}
\newcommand{\cL}{{\mathcal L}}
\newcommand{\cM}{{\mathcal M}}

\newcommand{\cD}{{\mathcal D}}

\newcommand{\cA}{{\mathcal A}}
\newcommand{\cB}{{\mathcal B}}
\newcommand{\cI}{{\mathcal I}}
\newcommand{\cC}{{\mathcal C}}

\newcommand{\cR}{{\mathcal R}}

\newcommand{\cH}{{\mathcal H}}

\newcommand{\cl}{\operatorname{cl}}
\newcommand{\qu}{\operatorname{qu}}

\newcommand{\DG}{\operatorname{DG}}
\newcommand{\Fun}{\operatorname{Fun}}

\newcommand{\Def}{\operatorname{Def}}
\newcommand{\Perf}{\operatorname{Perf}}

\newcommand{\Ker}{\operatorname{Ker}}
\newcommand{\im}{\operatorname{Im}}

\newcommand{\Ext}{\operatorname{Ext}}

\newcommand{\Id}{\operatorname{Id}}

\newcommand{\Ind}{\operatorname{Ind}}
\newcommand{\Res}{\operatorname{Res}}

\newcommand{\dgart}{\operatorname{dgart}}
\newcommand{\art}{\operatorname{art}}

\newcommand{\coDef}{\operatorname{coDef}}
\newcommand{\cart}{\operatorname{cart}}
\newcommand{\Alg}{\operatorname{Alg}}

\newcommand{\Ho}{\operatorname{Ho}}

\newcommand{\id}{\operatorname{id}}

\newcommand{\dgalg}{\operatorname{dgalg}}
\newcommand{\DEF}{\operatorname{DEF}}
\newcommand{\coDEF}{\operatorname{coDEF}}
\newcommand{\ev}{\operatorname{ev}}
\newcommand{\adgalg}{\operatorname{adgalg}}

\title{DG deformation theory of objects  in homotopy and derived categories I}

\author{Valery A.~Lunts}
\address{Department of Mathematics, Indiana University,
Bloomington, IN 47405, USA} \email{vlunts@indiana.edu}
\author{Dmitri Orlov}
\address{Steklov Mathematical Institute, 8 Gubkina St. Moscow, Russia}
\email{orlov@mi.ras.ru}

\begin{document}

\begin{abstract} We develop a general deformation theory of objects
in homotopy and derived categories of DG categories. The main result
is a general pro-representability theorem for the corresponding
deformation functor.
\end{abstract}

\maketitle

\tableofcontents

\section{Introduction}

\subsection{}
It is well known ([De],[Dr2]) that for many mathematical objects $X$
(defined over a field of characteristic zero) the formal deformation
theory of $X$ is controlled by a DG Lie algebra
$\frak{g}=\frak{g}(X)$ of (derived) infinitesimal automorphisms of
$X$. This is so in case $X$ is an algebra, a compact complex
manifold, a principal $G$-bundle, etc..

Let $\cM(X)$ denote the base of the universal deformation of $X$ and
$o\in \cM(X)$ be the point corresponding to $X$. Then (under some
conditions on $\frak{g}$) the completion of the local ring
$\hat{\cO}_{\cM(X),o}$ is naturally isomorphic to the linear dual of
the homology space $H_0(\frak{g})$. The space $H_0(\frak{g})$ is a
co-commutative coalgebra, hence its dual is a commutative algebra.

The homology $H_0(\frak{g})$ is the zero cohomology group of
$B\frak{g}$ -- the bar construction of $\frak{g}$, which is a
co-commutative DG coalgebra. It is therefore natural to consider the
DG "formal moduli space" $\cM ^{DG}(X)$, so that the "local ring"
$\hat{\cO}_{\cM^{DG}(X),o}$ is the linear dual $(B\frak{g})^*$,
which is a commutative DG algebra. The space $\cM ^{DG}(X)$ is thus
the "true" universal deformation space of $X$; it coincides with
$\cM(X)$ in case $H^i(B\frak{g})=0$ for $i\neq 0$. See [Ka],
[Ci-FoKa1],[Ci-FoKa2] for some examples. In particular, it appears
that the primary object is not the DG algebra $(B\frak{g})^*$, but
rather the DG coalgebra $B\frak{g}$ (this is the point of view in
[Hi]).

Note that the passage from a DG Lie algebra $\frak{g}$ to the
commutative DG algebra $(B\frak{g})^*$ is an example of the Koszul
duality for operads [GiKa]. Indeed, the operad of DG Lie algebras is
Koszul dual to that of commutative DG algebras.

Given all that, no general deformation theory has been developed
yet. For example, it has not been proved that the complete DG
algebra $(B\frak{g})^*$ pro-represents the functor of infinitesimal
DG deformations of $X$.

\subsection{} This paper is concerned with a general deformation
theory in a slightly different context. Namely, we consider
deformations of "linear" objects $E$, such as objects in a homotopy
or a derived category. More precisely, $E$ is a DG module over a DG
category $\cA$. In this case the deformation theory of $E$ is
controlled by $\cB =\End (E)$ which is a DG {\it algebra} (and not a
DG Lie algebra). (This works equally well in positive
characteristic.) Then the DG formal deformation space of $E$ is the
"Spec" of the (noncommutative!) DG algebra $(B\cB )^*$ -- the linear
dual of the bar construction $B\cB$ which is a DG coalgebra. Again
this is in agreement with the Koszul duality for operads, since the
operad of DG algebras is self-dual. (All this was already
anticipated in [Dr2].)

More precisely, let $\dgart$ be the category of local artinian (not
necessarily commutative) DG algebras and $\bf{Gpd}$ be the
2-category of groupoids. We define a 2-functor
$$\Def (E):\dgart \to \bf{Gpd},$$
which assigns to an artinian DG algebra $\cR$ the groupoid $\Def
_{\cR}(E)$  of $\cR$-deformations of $E$ in the derived category
$D(\cA)$. Actually we prefer to work with the 2-functor of {\it
co-deformations}
$$\coDef (E):\dgart \to \bf{Gpd},$$
which in many cases is equivalent to $\Def(E)$. The main result of
this paper is the pro-representability theorem for the functor
$\coDef(E)$ (under some finiteness assumptions on the DG algebra
$\End(E)$) by the complete local DG algebra $(B\cB)^*$.

Classically one defines representability for functors with values in
the category of sets. However, a deformation functor naturally takes
values in the 2-category of groupoids. Therefore in order to achieve
the correct pro-representability of such a functor its source
category must be also a 2-category. The category of DG categories is
naturally a 2-category with 1-morphisms being quasi-functors [Ke].
In particular this makes the category of DG algebras a 2-category.
We introduce the 2-category $2\text{-}\dgalg$ of {\it augmented} DG
algebras and prove the existence of the 2-functor
$\coDEF(E):2\text{-}\dgalg \to \bf{Gpd}$ which is a "lift" of the
functor $\coDef (E)$. Then the main theorem asserts that under some
finiteness conditions the 2-functor $\coDEF(E)$ is equivalent to the
2-functor
$$h_{(B\cB)^*}(\cdot):=\Hom
_{2\text{-}\dgalg}((B\cB)^*,\cdot):2\text{-}\dgart \to \bf{Gpd}.$$

\subsection{}
Let us briefly describe the contents of the paper.

Part 1 is a rather lengthy review of basics of DG categories and DG
modules over them with some minor additions that we did not find in
the literature. The reader who is familiar with basic DG categories
is suggested to go directly to Part 2, except for looking up the
definition of the DG functors $i^*$ and $i^!$.

Part 2 contains the definition and study of various deformation
functors. First we introduce the {\it homotopy} deformation and
co-deformation functors $\Def ^{\h}(E)$ and $\coDef ^{\h}(E)$ and
prove that they are always equivalent. We present a version of the
well known invariance theorem (Deligne's theorem) for these
functors. Then we introduce the {\it derived} deformation and
co-deformation functors $\Def (E)$ and $\coDef (E)$. We prove that
under some boundedness conditions the functors $\Def ^{\h}(E)$ and
$\Def (E)$ are equivalent (resp. $\coDef ^{\h}(E)$ and $\coDef (E)$
are equivalent). Actually, we only prove the equivalence of the
restrictions of these functors to the category $\dgart _-$ of {\it
non-positive} artinian DG algebras.

Part 3 contains the main pro-representability theorem. We start with
reviewing the bar construction and its relevance to the deformation
theory (known at least since [Q]). Then for an augmented DG algebra
$\cC$ (satisfying some finiteness conditions) and its bar
construction $B\cC$ we study some natural functors between the
derived categories $D((B\cC)^*)$ and $D(\cC)$ defined by the bar
complex $B\cC \otimes \cC$. Finally, we prove the
pro-representability theorem.

\subsection{}
In the sequel [LO] of this paper we plan to include the following:
pro-representability of the 2-functor $\DEF$, pro-representability
of the deformation functor controlled by a DG Lie algebra (as
explained above), applications of the theory to deformations of
objects in derived categories of {\it abelian} categories.

In is our pleasure to thank A.Bondal, P.Deligne, M.Mandell, M.Larsen
and P.Bressler for useful discussions. We especially appreciate the
generous help of B.Keller. We also thank W.Goldman and V.Schechtman
for sending us copies of letters [De] and [Dr2] respectively and
W.Lowen for sending us the preprint [Lo].

\part{Preliminaries on DG categories}

\section{Artinian DG algebras}

 We fix a field $k$. All algebras are assumed to be
$\bbZ$ graded $k$-algebras with unit and all categories are
$k$-linear. Unless mentioned otherwise $\otimes $ means $\otimes
_k$.

For a homogeneous element $a$ we denote its degree by $\bar{a}$.

A {\it module} always means a (left) graded module.

A DG algebra $\cB=(\cB ,d_{\cB})$ is a (graded) algebra with a map
$d=d_{\cB}:\cB \to \cB$ of degree 1 such that $d^2=0$, $d(1)=0$ and
$$d(ab)=d(a)b+(-1)^{\bar{a}}ad(b).$$

Given a DG algebra $\cB$ its opposite is the DG algebra $\cB ^0$
which has the same differential as $\cB$ and multiplication
$$a\circ b=(-1)^{\bar{a}\bar{b}}ba,$$
where $ba$ is the product in $\cB$.

We denote by $\dgalg$ the category of DG algebras.

A (left) DG module over a DG algebra $\cB$ is called a DG
$\cB$-module or, simply a $\cB$-module. A {\it right} $\cB$-module
is a DG module over $\cB ^0$.  We denote by $\cB \text{-mod}$ the
abelian category of $\cB$-modules.

If $\cB$ is a DG algebra and $M$ is a usual (not DG) module over the
algebra $\cB$, then we say that $M^{\gr}$ is a $\cB ^{\gr}$-module.

An {\it augmentation} of a DG algebra $\cB$ is a (surjective)
homomorphism of DG algebras $\cB \to k$. Its kernel is denoted
$\overline{\cB}$, this is  DG ideal (i.e. an ideal closed under the
differential) of $\cB$. Denote by $\adgalg$ the category of
augmented DG algebras (morphisms commute with the augmentation).

\begin{defi} Let $R$ be an algebra. We call $R$
{\it artinian}, if it is finite dimensional  and has a (graded)
nilpotent two-sided (maximal) ideal $m\subset R$, such that $R/m=k$.
\end{defi}

\begin{defi} Let $\cR $ be an augmented DG algebra.
We call $\cR$ {\it artinian} if $\cR$ is  artinian as an algebra and
the maximal ideal $m\subset R$ is a DG ideal, i.e. the quotient map
$\epsilon _{\cR} :R\to R/m$ is an augmentation of the DG algebra
$\cR$. (So $m=\overline{\cR}$). Note that a homomorphism of artinian
DG algebras automatically commutes with the augmentations. Denote by
$\dgart$ the category of artinian DG algebras.
\end{defi}

\begin{defi} An artinian DG algebra $\cR$ is called positive (resp. negative) if
negative (resp. positive) degree components of $\cR$ are zero.
Denote by $\dgart _+$ and $\dgart _-$ the corresponding full
subcategories of $\dgart$. Let $\art :=\dgart _-\cap \dgart _+$ be
the full subcategory of $\dgart$ consisting of (not necessarily
commutative) artinian algebras concentrated in degree zero. Denote
by $\cart \subset \art$ the full subcategory of commutative artinian
algebras.
\end{defi}

Given a DG algebra $\cB$ one studies the category $\cB\text{-mod}$
and the corresponding homotopy and derived categories. A
homomorphism of DG algebras induces various functors between these
categories. We will recall these categories and functors in the more
general context of DG categories in the next section.

\section{DG categories}

In this section we recall some basic facts about DG categories which
will be needed in this paper. Our main references here are
[BoKa],[Dr],[Ke].

 A DG category is a $k$-linear category $\cA$ in which the sets $\Hom (A,B)$, $A,B\in Ob\cA$,
 are
 provided
with a structure of a $\bbZ$-graded $k$-module and a differential
$d:\Hom(A,B)\to \Hom (A,B)$ of degree 1, so that for every $A,B,C\in
\cA$ the composition $\Hom (A,B)\times \Hom (B,C) \to \Hom (A,C)$
comes from a morphism of complexes $\Hom (A,B)\otimes \Hom (B,C) \to
\Hom (A,C)$. The identity morphism $1_A\in \Hom (A,A)$ is closed of
degree zero.

The simplest example of a DG category is the category $DG(k)$ of
complexes of $k$-vector spaces, or DG $k$-modules.

Note also that a DG algebra is simply a DG category with one object.

Using the supercommutativity isomorphism $S\otimes T\simeq
T\otimes S$ in the category of DG $k$-modules one defines for
every DG category $\cA$ the opposite DG category $\cA ^0$ with
$Ob\cA ^0=Ob\cA$, $\Hom_{\cA ^0}(A,B)=\Hom _{\cA}(B,A)$.
We denote by $\cA ^{\gr}$ the {\it graded }
 category which is obtained from $\cA$ by forgetting the differentials on $\Hom $'s.

The tensor product of DG-categories $\cA$ and $\cB$ is defined as
follows:

(i) $Ob(\cA \otimes \cB):=Ob\cA \times Ob\cB$; for $A\in Ob\cA$
and $B\in Ob\cB$ the corresponding object is denoted by $A\otimes
B$;

(ii) $\Hom(A\otimes B,A^\prime \otimes B^\prime):=\Hom
(A,A^\prime)\otimes \Hom (B,B^\prime)$ and the composition map is
defined by $(f_1\otimes g_1)(f_2\otimes g_2):=
(-1)^{\bar{g_1}\bar{f_2}}f_1f_2\otimes g_1g_2.$

Note that the DG categories $\cA \otimes \cB$ and $\cB \otimes \cA$
are canonically isomorphic. In the above notation the isomorphism DG
functor $\phi$ is
$$\phi (A\otimes B)=(B\otimes A), \quad \phi(f\otimes g)=(-1)^{\bar{f}\bar{g}}(g\otimes f).$$

Given a DG category $\cA$ one defines the graded category $\Ho^\cdot
(\cA)$ with $Ob\Ho^\cdot (\cA)=Ob\cA$ by replacing each $\Hom$
complex by the direct sum of its cohomology groups. We call
$\Ho^\cdot (\cA)$ the {\it graded homotopy category} of $\cA$.
Restricting ourselves to the 0-th cohomology of the $\Hom $
complexes we get the {\it homotopy category} $\Ho(\cA)$.

Two objects $A,B\in Ob\cA$ are called DG {\it isomorphic} (or,
simply, isomorphic) if there exists an invertible degree zero
morphism $f\in \Hom(A,B)$. We say that $A,B$ are {\it homotopy
equivalent} if they are isomorphic in  $\Ho(\cA)$.

A DG-functor between DG-categories $F:\cA \to \cB$ is said to be a
{\it quasi-equivalence} if $\Ho(F):\Ho(\cA)\to \Ho(\cB)$ is an
equivalence. We say that $F$ is a DG {\it equivalence} if it is
fully faithful and every object of $\cB$ is DG isomorphic to an
object of $F(\cA)$. Certainly, a DG equivalence is a
quasi-equivalence. DG categories $\cC$ and $\cD$ are called {\it
quasi-equivalent} if there exist DG categories $\cA _1,...,\cA _n$
and a chain of quasi-equivalences
$$\cC \leftarrow \cA _1 \rightarrow ...\leftarrow \cA _n \rightarrow \cD.$$

Given DG categories $\cA$ and $\cB$ the collection of covariant DG
functors $\cA \to \cB$ is itself the collection of objects of a DG
category, which we denote by $\Fun _{\DG}(\cA ,\cB)$. Namely, let
$\Phi $ and $\Psi$ be two DG functors. Put $\Hom ^k(\Phi ,\Psi)$
equal to the set of natural transformations $t:\Phi ^{\gr} \to
\Psi ^{\gr}[k]$ of graded functors from $\cA ^{\gr}$ to $\cB
^{\gr}$. This means that for any morphism $f \in
\Hom_{\cA}^s(A,B)$ one has
$$\Psi (f )\cdot t(A)=(-1)^{ks}t(B)\cdot \Phi (f).$$
On each $A\in \cA$ the differential of the transformation $t$ is
equal to $(dt)(A)$ (one easily checks that this is well defined).
Thus, the closed transformations of degree 0 are the DG
transformations of DG functors. A similar definition gives us the
DG-category
 consisting of the contravariant DG functors
 $\Fun _{\DG}(\cA ^0 ,\cB)=\Fun _{\DG}(\cA  ,\cB ^0)$
 from $\cA$ to $\cB$.

\subsection{DG modules over DG categories}
We denote the DG category $\Fun _{\DG}(\cA ,DG(k))$ by $\cA
\text{-mod}$ and call it the
 category
of DG $\cA$-modules. There is a natural covariant DG functor $h:\cA
\to \cA ^{0}\text{-mod}$ (the Yoneda embedding) defined by
$h^A(B):=\Hom _{\cA}(B,A)$. As in the "classical" case one verifies
that the functor $h$ is fully faithful, i.e. there is a natural
isomorphism of complexes
$$\Hom _{\cA}(A,A^\prime)=\Hom_{\cA ^{0}\text{-mod}}(h^A,h^{A^\prime}).$$
Moreover, for any $M\in \cA ^{0}\text{-mod}$, $A\in \cA$
$$\Hom _{\cA ^0\text{-mod}}(h^A,M)=M(A).$$

The DG $\cA^0$-modules $h^A$, $A\in \cA$ are called {\it
free}.

For $A\in \cA$ one may consider also the covariant DG functor
$h_A(B):=\Hom _{\cA}(A,B)$ and the contravariant DG functor
$h^*_A(B):=\Hom _k(h_A(B),k)$. For any $M\in \cA ^0\text{-mod}$ we
have
$$\Hom _{\cA ^0\text{-mod}}(M,h^*_A)=\Hom _k(M(A),k).$$

 A DG $\cA^{0}$-module $M$ is called acyclic, if the complex $M(A)$
is acyclic for all $A\in \cA$. Let $D(\cA)$ denote the {\it derived
category} of DG $\cA ^{0}$-modules, i.e. $D(\cA)$ is the Verdier
quotient of the homotopy category $\Ho(\cA^{0}\text{-mod})$ by the
subcategory of acyclic DG-modules. This is a triangulated category.

A DG $\cA ^{0}$-module $P$ is called h-{\it projective} if for any
acyclic DG $\cA ^{0}$-module $N$ the complex $\Hom (P,N)$ is
acyclic. A free DG module is h-projective. Denote by $\P(\cA)$ the
full DG subcategory of $\cA^{0}\text{-mod}$ consisting of
h-projective DG modules.

Similarly, a  DG $\cA ^{0}$-module $I$ is called h-{\it injective}
if for any acyclic DG $\cA ^{0}$-module $N$ the complex $\Hom (N,I)$
is acyclic. For any $A\in \cA$ the DG $\cA ^0$-module $h^*_A$ is
h-injective.  Denote by $\cI(\cA)$ the full DG subcategory of
$\cA^{0}\text{-mod}$ consisting of homotopically injective DG
modules.

For any DG category $\cA$ the DG categories $\cA^{0}\text{-mod}$,
$\P(\cA)$, $\cI(\cA)$ are (strongly) pre-triangulated ([Dr],[BoKa],
also see subsection 3.5 below). Hence the homotopy categories
$\Ho(\cA^{0}\text{-mod})$, $\Ho(\P(\cA))$, $\Ho(\cI(\cA))$ are
triangulated.

The following theorem was proved in [Ke].

\begin{theo}  The inclusion functors $\P(\cA )\hookrightarrow \cA
^{0}\text{-mod}$, $\cI(\cA )\hookrightarrow \cA ^{0}\text{-mod}$
 induce  equivalences of triangulated categories $\Ho(\P(\cA))\simeq
 D(\cA)$ and $\Ho(\cI(\cA))\simeq
 D(\cA)$.
 \end{theo}

Actually, it will be convenient for us to use some more precise
results from [Ke]. Let us recall the relevant definitions.

\begin{defi} A DG $\cA ^0$-module $M$ is called relatively
projective if $M$ is a direct summand of a direct sum of DG
$\cA^0$-modules of the form $h^A[n]$, $A\in \cA$, $n\in \bbZ$. A DG
$\cA ^0$-module $P$ is said to have property (P) if it admits a
filtration
$$0=F_{-1}\subset F_0\subset F_1\subset ... P$$
such that

\noindent(F1) $\cup_iF_i=P$;

\noindent(F2) the inclusion $F_i\hookrightarrow F_{i+1}$ splits as a
morphism of graded modules;

\noindent(F3) each quotient  $F_{i+1}/F_i$ is a relatively
projective DG $\cA ^0$-module.
\end{defi}

\begin{defi}  A DG $\cA ^0$-module $M$ is called relatively
injective if $M$ is a direct summand of a direct product of DG
$\cA^0$-modules of the form $h_A^*[n]$, $A\in \cA$, $n\in \bbZ$. A
DG $\cA ^0$-module $I$ is said to have property (I) if it admits a
filtration
$$I=F_{0}\supset F_1\supset ...$$
such that

\noindent(F1') the canonical morphism
$$I\to \lim_{\leftarrow}I/F_i$$
is an isomorphism;

\noindent(F2') the inclusion $F_{i+1}\hookrightarrow F_i$ splits as
a morphism of graded modules;

\noindent(F3') each quotient  $F_{i}/F_{i+1}$ is a relatively
injective DG $\cA ^0$-module.
\end{defi}

\begin{theo} ([Ke]) a) A DG $\cA ^0$-module with property (P) is
$h$-projective.

b) For any $M\in \cA ^0\text{-mod}$ there exists a quasi-isomorphism
$P\to M$, such that the DG $\cA ^0$-module $P$ has property (P).

c)  A DG $\cA ^0$-module with property (I) is $h$-injective.

d) For any $M\in \cA ^0\text{-mod}$ there exists a quasi-isomorphism
$M\to I$, such that the DG $\cA ^0$-module $I$ has property (I).
\end{theo}

\begin{remark} a) Assume that a DG $\cA ^0$-module $M$ has an
increasing filtration $M_1\subset M_2\subset ...$ such that $\cup
M_i=M$, each inclusion $M_i\hookrightarrow M_{i+1}$ splits as a
morphism of graded modules,  and each subquotient $M_{i+1}/M_i$ is
$h$-projective. Then $M$ is h-projective. b) Assume that a DG $\cA
^0$-module $N$ has a decreasing filtration $N=N_1\supset N_2\supset
...$ such that $\cap N_i=0$, each inclusion $N_{i+1}\hookrightarrow
N_i$ splits as a morphism of graded modules,  each subquotient
$N_i/N_{i+1}$ is h-injective (hence $N/N_i$ is h-injective for each
$i$) and the natural map
$$N\to \lim_{\leftarrow}N/N_i$$
is an isomorphism. Then $N$ is h-injective.
\end{remark}

\subsection{Some DG functors}
Let $\cB$ be a small DG category. The complex
$$\Alg _{\cB}:=\bigoplus _{A,B\in Ob \cB}\Hom(A,B)$$
has a natural structure of a DG algebra possibly without a unit. It
has the following property: every finite subset of $\Alg _{\cB}$ is
contained in $e\Alg _{\cB} e$ for some idempotent $e$ such that
$de=0$ and $\deg e=0$. We say that a DG module $M$ over $\Alg
_{\cB}$ is {\it quasi-unital} if every element of $M$ belongs to
$eM$ for some idempotent $e\in \Alg _{\cB}$ (which may be assumed
closed of degree $0$ without loss of generality). If $\Phi $ is a DG
$\cB$-module then
$$M_{\Phi}:=\oplus _{A\in Ob \cB}\Phi (A)$$
 is a quasi-unital
DG module over $\Alg _{\cB}$. Thus we get a DG equivalence between
DG category of DG $\cB$-modules and that of quasi-unital DG
modules over $\Alg _{\cB}$.

Recall that a homomorphism of (unitary) DG algebras $\phi :\cA \to
\cB$ induces functors
$$\phi _*:\cB^{0}\text{-mod}\to \cA^{0}\text{-mod},$$
$$\phi ^*:\cA^{0}\text{-mod}\to \cB^{0} \text{-mod}$$
$$\phi ^!:\cA^{0}\text{-mod}\to \cB^{0} \text{-mod}$$
where $\phi _*$ is the restriction of scalars, $\phi ^*(M)=M \otimes
_{\cA}\cB$ and $\phi ^!(M)=\Hom _{\cA ^0}(\cB,M)$. The DG functors
$(\phi ^*,\phi _*)$  and $(\phi _*,\phi ^!)$ are adjoint: for $M\in
\cA^{0}\text{-mod}$ and $N\in \cB^{0}\text{-mod}$ there exist
functorial isomorphisms of complexes
$$\Hom (\phi ^*M,N)=\Hom (M,\phi _*N),\quad
\Hom (\phi _*N,M)=\Hom (N,\phi ^!M). $$

This generalizes to a DG functor $F:\cA \to \cB$ between DG
categories. We obtain DG functors
$$F _*:\cB^0\text{-mod}\to \cA^0\text{-mod},$$
$$F ^*:\cA^0\text{-mod}\to \cB^0\text{-mod}.$$
$$F ^!:\cA^0\text{-mod}\to \cB^0\text{-mod}.$$

Namely, the DG functor $F$ induces a homomorphism of DG algebras
$F:\Alg _{\cA}\to \Alg _{\cB}$ and hence defines functors $F_*$,
$F^*$ between quasi-unital DG modules as above. (These functors
$F_*$ and $F^*$ are denoted in [Dr] by $\Res _F$ and $\Ind _F$
respectively.) The functor $F^!$ is defined as follows: for a
quasi-unital $\Alg _{\cA}^0$-module $M$ put
$$F^!(M)=\Hom _{\Alg _{\cA}^0}(\Alg _{\cB},M)^{\qu},$$
where $N^{\qu}\subset N$ is the {\it quasi-unital} part of a $\Alg
_{\cB}^0$-module $N$ defined by
$$N^{\qu}:=\im (N\times  \Alg _{\cB}\to N).$$

The DG functors $(F ^*,F _*)$ and $(F_*,F^!)$ are adjoint.

\begin{lemma} Let $F:\cA \to \cB$ be a DG functor. Then

a) $F_*$ preserves acyclic DG modules;

b) $F^*$ preserves h-projective DG modules;

c) $F^!$ preserves h-injective DG modules.

\end{lemma}

\begin{proof} The first assertion is obvious and the other two
follow by adjunction.
\end{proof}

By Theorem 3.1 above the DG subcategories $\P(\cA)$ and $\cI(\cA)$
of $\cA ^{0}\text{-mod}$ allow us to define (left and right) derived
functors of DG functors $G:\cA ^{0}\text{-mod}\to \cB
^{0}\text{-mod}$ in the usual way. Namely for a DG $\cA ^{0}$-module
$M$ choose quasi-isomorphisms $P\to M$ and $M\to I$ with $P\in
\P(\cA)$ and $I\in \cI(\cA)$. Put
$$\bL G(M):=G(P),\quad \quad \bR G(M):=G(I).$$
In particular for a DG functor $F:\cA \to \cB$ we will consider
 derived functors $\bL F^*:D(\cA)\to D(\cB)$, $\bR F^!:D(\cA)\to D(\cB)$. We also
have the obvious functor $F_*:D(\cB)\to D(\cA)$. The functors $(\bL
F^*,F_*)$  and $(F_*, \bR F^!)$ are adjoint.

\begin{prop} Assume that the DG functor $F:\cA \to \cB$ is a
quasi-equivalence. Then

a) $F^*:\P(\cA)\to \P(\cB)$ is a quasi-equivalence;

b) $\bL F ^*:D(\cA)\to D(\cB)$ is an equivalence;

c) $F_*:D(\cB)\to D(\cA)$ is an
equivalence.

d) $\bR F^!:D(\cA)\to D(\cB)$ is an equivalence.

e) $F^!:\cI (\cA)\to \cI (\cB)$ is a quasi-equivalence.

\end{prop}

\begin{proof} a) is proved in [Ke] and it implies b) by
Theorem 2.1. c) (resp. d)) follows from b) (resp. c) by adjunction.
Finally, e) follows from d) by Theorem 2.1.
\end{proof}

Given DG $\cA ^{0}$-modules $M,N$ we denote by $\Ext ^n(M,N)$
the group of morphisms $\Hom ^n _{D(\cA)}(M,N)$.

\subsection{DG category $\cA _{\cR}$} Let $\cR$ be a DG
algebra. We may and will consider $\cR$ as a DG category with one
object whose endomorphism DG algebra is $\cR$. We denote this DG
category again by $\cR$. Note that the DG category
$\cR^0\text{-mod}$ is just the category of right DG modules over the
DG algebra $\cR$.

For a DG category $\cA$ we denote the DG category $\cA \otimes \cR$
by $\cA _{\cR}$. Note that the collections of objects of $\cA$ and
$\cA _{\cR}$ are naturally identified. A homomorphism of DG algebras
$\phi :\cR\to \cQ$ induces the obvious DG functor $\phi=\id \otimes
\phi :\cA _{\cR}\to \cA _{\cQ}$ (which is the identity on objects),
whence the DG functors $\phi _*$, $ \phi ^*$, $\phi ^!$ between the
DG categories $\cA^0 _{\cR}\text{-mod}$ and $\cA ^0
_{\cQ}\text{-mod}$. For  $M \in \cA_{\cR}^0 \text{-mod}$
 we have
 $$\phi ^*(M)=M\otimes _{\cR}{\cQ}.$$
Also in case $\cQ^{\gr}$ is a finitely generated $\cR ^{\gr}$-module
we have
$$ \phi ^!(M)=\Hom _{\cR^0}(\cQ ,M).$$

In particular, if $\cR$ is  augmented then the canonical morphisms
of DG algebras $p:k\to \cR$ and $i:\cR \to k$ induce functors
$$p:\cA \to \cA _{\cR},\quad i:\cA _{\cR}\to \cA,$$
such that $i\cdot p=\Id_{\cA}$.  So for $S\in \cA ^0\text{-mod}$ and
$T\in \cA ^0_{\cR}\text{-mod}$ we have

$$p^*(S)=S\otimes _k\cR, \quad i^*(T)=T\otimes _{\cR}k, \quad i^!(T)=\Hom _{\cR^0}(k,T).$$

For an artinian DG algebra $\cR$ we denote by $R^*$ the DG $\cR
^0$-module $\Hom _k(\cR,k)$. This is a left $\cR$-module by the
formula
$$rf(q):=(-1)^{(\bar{f}+\bar{q})\bar{r}}f(qr)$$
and a right $\cR$-module by the formula
$$fr(p):=f(rp)$$
for $r,p\in \cR$ and $f\in \cR ^*$. The augmentation map $\cR \to k$
defines the canonical (left and right) $\cR$-submodule $k\subset R
^*$. Moreover, the embedding $k\hookrightarrow \cR ^*$ induces an
isomorphism $k\to \Hom _{\cR}(k,\cR ^*)$.

\begin{defi} Let $\cR$ be an artinian DG algebra. A DG $\cA^0
_{\cR}$-module $M$ is called graded $\cR$-free (resp. graded
$\cR$-cofree) if there exists a DG $\cA ^0$-module $K$ such that
$M^{\gr}\simeq (K\otimes \cR)^{\gr}$ (resp. $M^{\gr}\simeq (K\otimes
\cR^*)^{\gr}$). Note that for such $M$ one may take $K=i^*M$ (resp.
$K=i^!M$).
\end{defi}

\begin{lemma} Let $\cR$ be an artinian DG algebra.

a) The full DG subcategories of DG $\cA _{\cR}^0$-modules consisting
of graded $\cR$-free (resp. graded $\cR$-cofree)  modules are DG
isomorphic. Namely, if $M\in \cA_{\cR}^0\text{-mod}$ is graded
$\cR$-free (resp. graded $\cR$-cofree) then $M\otimes _{\cR}\cR ^*$
(resp. $\Hom _{\cR^0}(\cR ^*,M)$) is graded $\cR$-cofree (resp.
graded $\cR$-free).

b) Let $M$ be a graded $\cR$-free module. There is a natural
isomorphism of DG $\cA ^0$-modules
$$i^*M\stackrel{\sim}{\to}i^!(M\otimes _{\cR}\cR ^*).$$
\end{lemma}

\begin{proof} a) If $M$ is graded $\cR$-free, then obviously $M\otimes
_{\cR}\cR ^*$ is graded $\cR$-cofree. Assume that $N$ is graded
$\cR$-cofree, i.e. $N^{\gr}=(K\otimes \cR ^*)^{\gr}$. Then
$$(\Hom
_{\cR^0}(\cR ^*,N))^{\gr}=(K\otimes \Hom _{\cR^0}(\cR
^*,\cR^*))^{\gr},$$ since $\dim _k\cR <\infty$. On the other hand
$$\Hom _{\cR^0}(\cR ^*,\cR^*)=\Hom _{\cR^0}(\cR ^*,\Hom _k(\cR ,k))=\Hom _k(\cR ^*\otimes
_{\cR}\cR,k)=\cR,$$ so $(\Hom _{\cR^0}(\cR ^*,N))^{\gr}=(K\otimes
\cR )^{\gr}$.

b) For an arbitrary DG $\cA _{\cR}^0$-module $M$ we have a natural
(closed degree zero) morphism of DG $\cA ^0$-modules
$$i^*M\to i^!(M\otimes _{\cR}\cR^*),\quad m\otimes 1\mapsto (1\mapsto m\otimes i),$$
where $i:\cR \to k$ is the augmentation map. If $M$ is graded
$\cR$-free this map is an isomorphism.
\end{proof}

\begin{prop} Let $\cR$ be an artinian DG algebra. Assume that a DG $\cA ^0_{\cR}$-module $M$
satisfies property (P) (resp. property (I)). Then $M$ is graded
$\cR$-free (resp. graded $\cR$-cofree).
\end{prop}

\begin{proof}  Notice that the collection of graded $\cR$-free
objects in $\cA ^0_{\cR}\text{-mod}$ is closed under taking direct
sums, direct summands (since $\cR$ is a local ring) and direct
products (since $\cR$ is finite dimensional). Similarly for graded
$\cR$-cofree objects since the DG functors in Lemma 3.9 a) preserve
direct sums and products. Also notice that for any $A \in \cA_{\cR}$
the DG $\cA ^0_{\cR}$-module $h^A$ (resp. $h_A^*$) is graded
$\cR$-free (resp. graded $\cR$-cofree). Now the proposition follows
since a DG $\cA ^0_{\cR}$-module $P$ (resp. $I$) with property (P)
(resp. property (I)) as a graded module is a direct sum of
relatively projective DG modules (resp. a direct product of
relatively injective DG modules).
\end{proof}

\begin{cor} Let $\cR$ be an artinian DG algebra. Then for any DG
$\cA^0_{\cR}$-module $M$ there exist quasi-isomorphisms $P\to M$ and
$M\to I$ such that $P\in \cP(\cA _{\cR})$, $I\in \cI (\cA _{\cR})$
and $P$ is graded $\cR$-free, $I$ is graded $\cR$-cofree.
\end{cor}

\begin{proof} Indeed, this follows from Theorem 3.4 and Proposition
3.10 above.
\end{proof}

\begin{prop} Let $\cR$ be an artinian
DG algebra and $S,T\in \cA ^0_{\cR}\text{-mod}$ be graded $\cR$-free
(resp. graded $\cR$-cofree).

a) There is an isomorphism of graded algebras
$\Hom(S,T)=\Hom(i^*S,i^*T) \otimes \cR$, (resp.
$\Hom(S,T)=\Hom(i^!S,i^!T) \otimes \cR$). In particular, the map $i
^*:\Hom (S,T)\to \Hom (i^*S,i^*T)$ (resp. $i ^!:\Hom(S,T)\to
\Hom(i^!S,i^!T)$) is surjective.

b) The DG module $S$ has a finite filtration with subquotients
isomorphic to $i^*S$ as DG $\cA ^0$-modules (resp. to $i_*i^*S$ as
DG $\cA _{\cR}^0$-modules).

c) The DG algebra $\End(S)$ has a finite filtration by DG ideals
with subquotients isomorphic to $\End (i^*S)$.

d) If $f\in \Hom (S,T)$ is a closed morphism of degree zero such
that $i^*f$ (resp. $i^!f$) is an isomorphism or a homotopy
equivalence or a quasi-isomorphism, then $f$ is also such.
\end{prop}

\begin{proof} Because of Lemma 3.9 above it suffices to
prove the proposition for graded $\cR$-free modules. So assume that
$S$, $T$ are graded $\cR$-free.

a) This holds because $\cR$ is finite dimensional.

b) We can refine the filtration of $\cR$ by powers of the maximal
ideal to get a filtration  $F_i\cR$ by ideals with 1-dimensional
subquotients (and zero differential). Then the filtration
$F_iS:=S\cdot F_i\cR$ satisfies the desired properties.

c) Again the filtration $F_i\End (S):=\End(S)\cdot F_i\cR$ has the
desired properties.

d) If $i^*f$ is an isomorphism, then $f$ is surjective by the
Nakayama lemma for $\cR$. Also $f$ is injective since $T$ if graded
$\cR$-free.

Assume that $i^*f$ is a homotopy equivalence. Let $C(f)\in \cA
_{\cR}^0\text{-mod}$ be the cone of $f$. (It is also graded
$\cR$-free.) Then $i^*C(f)\in \cA ^0\text{-mod}$ is the cone
$C(i^*f)$ of the morphism $i^*f$. By assumption the DG algebra $\End
(C(i^*f))$ is acyclic. But by part c) the complex $\End (C(f))$ has
a finite filtration with subquotients isomorphic to the complex
$\End (C(i^*f))$. Hence $\End (C(f))$ is also acyclic, i.e. the DG
module $C(f)$ is null-homotopic, i.e. $f$ is a homotopy equivalence.

Assume that $i^*f$ is a quasi-isomorphism. Then in the above
notation $C(i^*f)$ is acyclic. Since by part b) $C(f)$ has a finite
filtration with subquotients isomorphic to $C(i^*f)$, it is also
acyclic. Thus $f$ is a quasi-isomorphism.
\end{proof}

\subsection{More DG functors} So far we considered DG functors $F_*$, $F^*$,
$F^!$ between the DG categories $\cA ^0$-mod and $\cB ^0$-mod which
came from a DG functor $F:\cA \to \cB$. We will also need to
consider a different type of DG functors.

\begin{example}  For an artinian DG algebra $\cR$ and a small DG category $\cA$
we will consider two types of "restriction of scalars" DG functors
$\pi _*,\pi _!:\cA ^0_{\cR}\text{-mod}\to \cR ^0\text{-mod}$.
Namely, for $M\in \cA _{\cR}^0\text{-mod}$ put
$$\pi _*M:=\prod_{A\in Ob\cA _{\cR}}M(A),\quad \pi _!M:=\bigoplus_{A\in Ob\cA
_{\cR}}M(A).$$ We will also consider the two "extension of scalars"
functors $\pi ^*,\pi ^!:\cR ^0\text{-mod}\to \cA
^0_{\cR}\text{-mod}$ defined by
$$\pi ^*(N)(A):=N\otimes \bigoplus_{B\in Ob\cA }\Hom_{\cA }(A,B)
\quad \pi ^!(N)(A):=\Hom _k(\bigoplus_{B\in Ob\cA }\Hom_{\cA
}(B,A),N)$$ for  $A\in Ob\cA _{\cR}$. Notice that the DG functors
$(\pi ^*,\pi _*)$ and $(\pi _!,\pi ^!)$ are adjoint, that is for
$M\in \cA ^0_{\cR}\text{-mod}$ and $N\in \cR ^0\text{-mod}$ there is
a functorial isomorphism of complexes
$$\Hom (\pi ^*N,M)=\Hom (N,\pi_*M), \quad \Hom (\pi _!M,N)=\Hom (M,\pi^!N).$$

The DG functors $\pi^*,\pi ^!$  preserve acyclic DG modules, hence
$\pi _*$ preserves h-injectives and $\pi _!$ preserves
h-projectives.

We have the following commutative functorial diagrams
$$\begin{array}{ccc}
\cA _{\cR}^0\text{-mod} & \stackrel{i^*}{\longrightarrow} & \cA
^0\text{-mod}\\
\pi _!\downarrow & & \pi _!\downarrow \\
\cR ^0\text{-mod} & \stackrel{i^*}{\longrightarrow} & DG(k),
\end{array}$$
$$\begin{array}{ccc}
\cA _{\cR}^0\text{-mod} & \stackrel{i^!}{\longrightarrow} & \cA
^0\text{-mod}\\
\pi _*\downarrow & & \pi _*\downarrow \\
\cR ^0\text{-mod} & \stackrel{i^!}{\longrightarrow} & DG(k).
\end{array}$$
\end{example}

\begin{example} Fix $E\in \cA ^0\text{-mod}$ and put $\cB=\End(E)$.
Consider the DG functor
$$\Sigma =\Sigma ^E:\cB ^0\text{-mod}\to \cA ^0\text{-mod}$$
defined by $\Sigma (M)=M\otimes _{\cB}E$. Clearly, $\Sigma (\cB)=E$.
This DG functor gives rise to the functor
$$\bL \Sigma :D(\cB)\to D(\cA),\quad \quad \bL \Sigma
(M)=M\stackrel{\bL}{\otimes }_{\cB}E.$$
 \end{example}

\subsection{Pre-triangulated DG categories} For any DG category $\cA$ there
exists a DG category $\cA ^{pre-tr}$ and a canonical full and
faithful DG functor $F:\cA \to \cA^{\pre-tr}$ (see [BoKa],[Dr]). The
homotopy category $\Ho (\cA^{\pre-tr})$ is canonically triangulated.
The DG category $\cA$ is called {\it pre-triangulated} if the DG
functor $F$ is a quasi-equivalence. The DG category $\cA^{\pre-tr}$
is pre-triangulated.

Let $\cB$ be another DG category and $G:\cA \to \cB$ be a quasi-equivalence. Then
$G^{\pre-tr}:\cA ^{\pre-tr}\to \cB ^{\pre-tr}$ is also a quasi-equivalence.

The DG functor $F$ induces a DG isomorphism of DG categories
$F_*:(\cA^{\pre-tr})^0\text{-mod} \to \cA ^0\text{-mod}$. Hence the
functors $F_*:D(\cA^{\pre-tr})\to D(\cA)$ and $\bL F^*:D(\cA)\to
D(\cA ^{\pre-tr})$ are equivalences. We obtain the following
corollary.

\begin{cor} Assume that a DG functor $G_1:\cA \to \cB$ induces a quasi-equivalence
$G_1^{\pre-tr}:\cA ^{\pre-tr}\to \cB ^{\pre-tr}$. Let $\cC$ be another DG category and consider
the DG functor $G:=G_1\otimes \id :\cA \otimes \cC\to \cB \otimes \cC$.
 Then the functors $G_*,\bL G^*, \bR G^!$ between
the derived categories $D(\cA \otimes \cC)$ and $D(\cB \otimes \cC)$ are
equivalences.
\end{cor}

\begin{proof} The DG functor $G$ induces the quasi-equivalence
$G ^{\pre-tr}:(\cA \otimes \cC)^{\pre-tr}\to (\cB \otimes \cC
)^{\pre-tr}$. Hence the corollary follows from the above discussion
and Proposition 3.6.
\end{proof}

\begin{example} Suppose $\cB$ is a pre-triangulated DG category. Let
$G_1:\cA \hookrightarrow \cB$ be an embedding of a full DG
subcategory so that the triangulated category $\Ho(\cB)$ is
generated by the collection of objects $Ob\cA$. Then the assumptions
of the previous corollary hold.
\end{example}

\subsection{A few lemmas}

\begin{lemma} Let $\cR$, $\cQ$ be DG algebras and $M$ be a DG
$\cQ\otimes \cR ^0$-module. Then for any DG modules $N$, $S$ over
the DG algebras $\cQ ^0$ and $\cR ^0$ respectively there is a
natural isomorphism of complexes $$\Hom _{\cR^0}(N\otimes
_{\cQ}M,S)=\Hom _{\cQ}(N,\Hom _{\cR^0}(M,S)).$$
\end{lemma}

\begin{proof} Indeed, for $f\in \Hom _{\cR^0}(N\otimes _{\cQ}M,S)$
define $\alpha (f)\in \Hom _{\cQ}(N,\Hom _{\cR^0}(M,S))$ by the
formula $\alpha (f)(n)(m)=f(n\otimes m)$. Conversely, for $g\in \Hom
_{\cQ}(N,\Hom _{\cR^0}(M,S))$ define $\beta (g)\in \Hom
_{\cR^0}(N\otimes _{\cQ}M,S)$ by the formula $\beta(g)=g(n)(m)$.
Then $\alpha $ and $\beta$ are mutually inverse isomorphisms of
complexes.
\end{proof}

\begin{lemma} Let $\cR$ be an artinian DG algebra. Then in the DG
category $\cR ^0\text{-mod}$ a direct sum of copies of $\cR ^*$ is
h-injective.
\end{lemma}

\begin{proof} Let $V$ be a graded vector space, $M=V\otimes \cR
^*\in \cR ^0\text{-mod}$ and $C$ an acyclic DG $\cR^0$-module.
Notice that $M=\Hom _k(\cR ,V)$ since $\dim \cR <\infty$. Hence the
complex
$$\Hom _{\cR ^0}(C,M)=\Hom _{\cR ^0}(C,\Hom _k(\cR ,V))=\Hom
_k(C\otimes _{\cR}\cR,V)=\Hom _k(C,V)$$ is acyclic.
\end{proof}

\begin{lemma} Let $\cB$ be a DG algebra, such that $\cB ^i=0$ for
$i>0$. Then the category $D(\cB)$ has truncation functors: for any
DG $\cB$-module $M$ there exists a short exact sequence in
$\cB\text{-mod}$
$$\tau _{<0}M\to M\to \tau _{\geq 0}M,$$
where $H^i(\tau _{<0}M)=0$ if $i\geq 0$ and $H^i(\tau _{\geq 0}M)=0$
for $i<0$.
\end{lemma}

\begin{proof} Indeed, put $\tau _{<0}M:=\oplus _{i<0}M^i\oplus d(M^{-1})$.
\end{proof}

\begin{lemma} Let $\cB$ be a DG algebra, s.t. $\cB ^i=0$ for $i>0$
and $\dim \cB ^i<\infty$ for all $i$. Let $N$ be a DG $\cB$-module
with finite dimensional cohomology. Then there exists an
h-projective DG $\cB$-module $P$ and a quasi-isomorphism $P\to N$,
where $P$ in addition satisfies the following conditions

a) $P^i=0$ for $i>>0$,

b) $\dim P^i<\infty$ for all $i$.
\end{lemma}

\begin{proof} First assume that $N$ is concentrated in one degree,
say $N^i=0$ for $i\neq 0$. Consider $N$ as a $k$-module and put
$P_0:=\cB \otimes N$. We have a natural surjective map of DG
$\cB$-modules $\epsilon :P_0\to N$ which is also surjective on the
cohomology. Let $K:=\Ker \epsilon$. Then $K^i=0$ for $i>0$ and $\dim
K^i<\infty$ for all $i$. Consider $K$ as a DG $k$-module and put
$P_{-1}:=\cB \otimes K$. Again we have a surjective map of DG
$\cB$-modules $P_{-1}\to K$ which is surjective and surjective on
cohomology. And so on. This way we obtain an exact sequence of DG
$\cB$-modules
$$...\to P_{-1}\to P_0\stackrel{\epsilon}{\to}N\to 0,$$
where $P_{-j}^i=0$ for $i<0$ and $\dim P_{-j}^i<\infty$ for all $j$.
Let $P:=\oplus _jP_{-j}[j]$ be the "total" DG $\cB$-module of the
complex $...\to P_{-1} \to P_0\to 0$. Then $\epsilon :P\to N$ is a
quasi-isomorphism. Since each DG $\cB$-module $P_{-j}$ has the
property (P), the module $P$ is h-projective by Remark 3.5a). Also
$P^i=0$ for $i<0$ and $\dim P^i<infty $ for all $i$.

How consider the general case. Let $H^s(N)=0$ and $H^i(N)=0$ for all
$i<s$. Replacing $N$ by $\tau _{\geq s}N$ (Lemma 3.19) we may and
will assume that $N^i=0$ for $i<s$. Then $M:=(\Ker d_N)\cap N^s$ is
a DG $\cB$-submodule of $N$ which is not zero. If the embedding
$M\hookrightarrow N$ is a quasi-isomorphism, then we may replace $N$
by $M$ and so we are done by the previous argument. Otherwise we
have a short exact sequence of DG $\cB$-modules
$$o\to M\to N\to N/M\to 0$$
with $\dim H(M), \dim H(N/M)<\dim H(N)$. By the induction on $\dim
H(N)$ we may assume that the lemma holds for $M$ and $N/M$. But then
it also holds for $N$.
\end{proof}

\begin{cor} Let $\cB$ be a DG algebra, s.t. $\cB ^i=0$ for $i>0$
and $\dim \cB ^i<\infty$ for all $i$. Let $N$ be a DG $\cB$-module
with finite dimensional cohomology. Then $N$ is quasi-isomorphic to
a finite dimensional DG $\cB$-module.
\end{cor}

\begin{proof} By Lemma 3.20 there exists a bounded above and locally
finite DG $\cB$-module $P$ which is quasi-isomorphic to $N$. It
remains to apply the appropriate truncation functor to $P$ (Lemma
3.19).
\end{proof}

\begin{cor} Let $\cB$ be an augmented DG algebra, s.t. $\cB ^i=0$ for $i>0$
and $\dim \cB ^i<\infty$ for all $i$. Denote by $\langle k\rangle
\subset D(\cB^0)$ the triangulated envelope of the DG $\cB
^0$-module $k$. Let $N$ be a DG $\cB$-module with finite dimensional
cohomology. Then $N\in \langle k\rangle$.
\end{cor}

\begin{proof} By the previous corollary we may assume that $N$ is
finite dimensional. But then $N$ has a filtration by DG
$\cB^0$-modules with subquotients isomorphic to $k$.
\end{proof}

\begin{lemma} Let $\cB$ and $\cC$ be DG algebras. Consider the DG
algebra $\cB \otimes \cC$ and a homomorphism of DG algebras $F:\cB
\to \cB \otimes \cC$, $F(b)=b\otimes 1$. Let $N$ be an h-projective
(resp. h-injective) DG $\cB \otimes \cC$-module. Then the DG
$\cB$-module $F_*N$ is also h-projective (resp. h-injective).
\end{lemma}

\begin{proof} The assertions follow from the fact that the DG
functor $F_*:\cB \otimes \cC\text{-mod}\to \cB\text{-mod}$ has a
left adjoint DG functor $F^*$ (resp. right adjoint DG functor $F^!$)
which preserves acyclic DG modules. Indeed,
$$F^*(M)=\cC \otimes _{k}M,\quad \quad F^!(M)=\Hom _k(\cC ,M).$$
\end{proof}

\part{Deformation functors}

\section{The homotopy deformation and co-deformation functors}

Denote by ${\bf Gpd}$ the 2-category of groupoids.

Let $\cA$ be a DG category and $E$ be a DG $\cA^0$-module. Let us
define the homotopy deformation 2-functor $\Def ^{\h} (E):\dgart \to
{\bf Gpd}$. This functor describes "infinitesimal" (i.e. along
artinian DG algebras) deformations of $E$ in the homotopy category
of DG $\cA^0$-modules.

\begin{defi} Let $\cR$ be an artinian DG algebra. An object in the
groupoid $\Def _{\cR}^{\h} (E)$ is a pair $(S,\sigma)$, where $S\in
\cA _{\cR}^0\text{-mod}$ and $\sigma :i^*S\to E$ is an isomorphism
of DG $\cA^0$-modules such that the following holds: there exists an
isomorphism of graded $\cA^0 _{\cR}$-modules $\eta :(E\otimes
\cR)^{\gr} \to S^{\gr}$ so that the composition
$$E= i^*(E\otimes \cR)
\stackrel{i^*(\eta)}{\to} i^*S\stackrel{\sigma}{\to}E$$ is the
identity.

Given objects $(S,\sigma),(S^\prime ,\sigma ^\prime)\in \Def
_{\cR}^{\h}(E)$ a map $f:(S,\sigma)\to (S^\prime,\sigma ^\prime)$ is
an isomorphism $f:S\to S^\prime$ such that $\sigma ^\prime \cdot
i^*f=\sigma$. An allowable homotopy between maps $f,g$ is a homotopy
$h:f\to g$ such that $i^*(h)=0$. We define morphisms in $\Def
_{\cR}^{\h}(E)$ to be classes of maps modulo allowable homotopies.

Note that a homomorphism of artinian DG algebras $\phi :\cR \to \cQ$
induces the functor $\phi ^*:\Def _{\cR}^{\h}(E)\to \Def
_{\cQ}^{\h}(E)$. This defines the 2-functor
$$\Def {^h}(E):\dgart \to {\bf Gpd}.$$
\end{defi}

We refer to objects of $\Def _{\cR}^{\h} (E)$ as homotopy
$\cR$-deformations of $E$.

\begin{example} We call $(p^*E,\id)\in \Def _{\cR}^{\h}(E)$ the
trivial $\cR$-deformation of $E$.
\end{example}

\begin{defi} Denote by $\Def _+^{\h}(E)$, $\Def _-^{\h}(E)$, $\Def
_0^{\h}(E)$, $\Def ^{\h}_{\cl}(E)$ the restrictions of the 2-functor
$\Def ^{\h}(E)$ to subcategories $\dgart _+$, $\dgart _-$, $\art$,
$\cart$ respectively.
\end{defi}

Let us give an alternative description of the same deformation
problem. We will define the homotopy {\it co-deformation} 2-functor
$\coDef^{\h}(E)$ and show that it is equivalent to $\Def ^{\h}(E)$.
The point is that in practice one should use $\Def ^{\h}(E)$ for a
h-projective $E$ and $\coDef ^{\h}(E)$ for a h-injective $E$ (see
Section 11).

For an artinian DG algebra $\cR$ recall the $\cR ^0$-module $\cR ^*
=\Hom _k(\cR ,k)$.

\begin{defi} Let $\cR$ be an artinian DG algebra. An object in the groupoid
$\coDef^{\h}_{\cR}(E)$ is a pair $(T, \tau)$, where $T$ is a DG $\cA
^0_{\cR}$-module and $\tau :E\to i^!T$ is an isomorphism of DG
$\cA^0$-modules so that the following holds: there exists an
isomorphism of graded $\cA ^0_{\cR}$-modules $\delta :T^{\gr}\to
(E\otimes \cR ^*)^{\gr}$ such that the composition
$$E \stackrel{\tau}{\to}i^!T \stackrel{i^!(\delta)}{\to} i^!(E\otimes \cR ^*)
 =E$$ is the identity.

Given objects $(T,\tau)$ and $(T^\prime,\tau ^\prime)\in \coDef
^h_{\cR}(E)$ a map $g:(T,\tau)\to (T^\prime ,\tau ^\prime)$ is an
isomorphism $f:T\to T^\prime$ such that $i^!f \cdot \tau =\tau
^\prime$. An allowable homotopy between maps $f,g$ is a homotopy
$h:f\to g$ such that $i^!(h)=0$. We define morphisms in $\coDef
_{\cR}^{\h}(E)$ to be classes of maps modulo allowable homotopies.

Note that a homomorphism of DG algebras $\phi :\cR \to \cQ$ induces
the functor $\phi ^!:\Def _{\cR}^{\h}(E)\to \Def _{\cQ}^{\h}(E)$.
This defines the 2-functor
$$\Def ^{\h}(E):\dgart \to {\bf Gpd}.$$
\end{defi}

We refer to objects of $\coDef _{\cR}^{\h} (E)$ as homotopy
$\cR$-co-deformations of $E$.

\begin{example} For example we can take $T=E\otimes \cR ^*$ with the
differential $d_{E,R^*}:=d_E\otimes 1+1\otimes d_{R^*}$ (and $\tau =\id$). This we
consider as the {\it trivial} $\cR$-co-deformation of $E$.
\end{example}

\begin{defi} Denote by $\coDef _+^{\h}(E)$, $\coDef _-^{\h}(E)$, $\coDef
_0^{\h}(E)$, $\coDef _{\cl}^{\h}(E)$ the restrictions of the
2-functor $\coDef ^{\h}(E)$ to subcategories $\dgart _+$, $\dgart
_-$, $\art$, $\cart$ respectively.
\end{defi}

\begin{prop} There exists a natural equivalence of 2-functors
$$\delta =\delta ^E:\Def ^{\h} (E)\to \coDef ^{\h}(E).$$
\end{prop}

\begin{proof} We use Lemma 3.9 above. Namely, let $S$ be an
$\cR$-deformation of $E$. Then $S\otimes _{\cR}\cR ^*$ is an
$\cR$-co-deformation of $E$. Conversely, given an
$\cR$-co-deformation $T$ of $E$ the  DG $\cA ^0_{\cR}$-module $\Hom
_{\cR^0}(\cR^*,T)$ is an $\cR$-deformation of $E$. This defines
mutually inverse equivalences $\delta _{\cR}$ and $\delta
_{\cR}^{-1}$  between $\Def _{\cR}^{\h}(E)$ and $\coDef ^{\h}
_{\cR}(E)$, which extend to morphisms between 2-functors $\Def ^{\h}
(E)$ and $\coDef ^{\h}(E)$. Let us be a little more explicit.

Let $\phi :\cR \to \cQ$ be a homomorphism of artinian DG algebras
and $S\in \Def ^{\h}(E)$. Then
$$\delta _{\cQ} \cdot \phi ^*(S)=S\otimes _{\cR}Q\otimes _{\cQ}\cQ
^*=S\otimes _{\cR}\cQ^*,\quad \quad \phi ^!\cdot \delta
_{\cR}(S)=\Hom _{\cR ^0}(\cQ ,S\otimes _{\cR}\cR ^*).$$ The
isomorphism $\alpha _{\phi}$ of these DG $\cA _{\cQ}^0$-modules is
defined by $\alpha _{\phi}(s\otimes f)(q)(r):=sf(q\phi (r))$ for
$s\in S$, $f\in \cQ^*$, $q\in \cQ$, $r\in \cR$. Given another
homomorphism $\psi :\cQ \to \cQ ^\prime$ of DG algebras one checks
the cocycle condition $\alpha _{\psi \phi}=\psi ^!(\alpha
_{\phi})\cdot \alpha _{\psi}$ (under the natural isomorphisms $(\psi
\phi)^*=\psi ^* \phi ^*$, $(\psi \phi)^!=\psi ^! \phi ^!$).
\end{proof}

\section{ Maurer-Cartan functor}

\begin{defi} For a DG algebra $\cC$ with the differential $d$ consider
the (inhomogeneous) quadratic map
$$Q:\cC ^1 \to \cC ^2; \quad Q(\alpha )=d\alpha +\alpha ^2.$$
We denote by $MC(\cC)$ the (usual) Maurer-Cartan cone
$$MC(\cC)=\{ \alpha \in \cC ^1\vert Q(\alpha )=0\}.$$
\end{defi}

Note that $\alpha \in MC(\cC)$ is equivalent to the operator
$d+\alpha :\cC \to \cC$ having square zero. Thus the set $MC(\cC)$
describes the space of "internal" deformations of the differential
in the complex $\cC$.

\begin{defi} Let $\cB$ be a DG algebra with the differential $d$
and  a nilpotent DG ideal $\cI \subset \cB$. We define the
Maurer-Cartan groupoid $\cM \cC (\cB ,\cI)$ as follows. The set of
objects of $\cM \cC (\cB ,\cI)$ is the cone $MC(\cI)$. Maps
between objects are defined by means of the gauge group $G(\cB
,\cI):=1+\cI ^0$ ($\cI ^0$ is the degree zero component of $\cI$) acting on $\cM \cC (\cB ,\cI)$ by the formula
$$g:\alpha \mapsto g\alpha g^{-1}+gd(g^{-1}),$$
where $g\in G(\cB  ,\cI)$, $\alpha \in MC(\cI)$. (This comes from
the conjugation action on the space of differentials $g:d+\alpha
\mapsto g(d+\alpha )g^{-1}$.) So if $g(\alpha)=\beta$, we call $g$
a map from $\alpha $ to $\beta$. Denote by $G(\alpha ,\beta)$ the
collection of such maps. We define the set $\Hom (\alpha , \beta)$
in the category $\cM \cC (\cB ,\cI)$ to consist of homotopy
classes of maps, where the homotopy relation is defined as
follows. There is an action of the group $\cI ^{-1}$ on the set
$G(\alpha ,\beta)$:
$$h:g\mapsto g+d(h)+\beta h+h\alpha,$$
for $h\in \cI ^{-1}, g\in G(\alpha ,\beta)$. We call two maps {\it
homotopic}, if they lie in the same $\cI ^{-1}$-orbit.
\end{defi}

To make the category $\cM \cC (\cB ,\cI)$ well defined we need to
prove a lemma.

\begin{lemma} Let $\alpha _1, \alpha _2, \alpha _3, \alpha
_4\in MC(\cI)$ and $g_1\in G(\alpha _1,\alpha _2)$, $g_1,g_3\in
G(\alpha _2,\alpha _3)$, $g_4\in G(\alpha _3 ,\alpha _4)$. If
$g_2$ and $g_3$ are homotopic, then so are $g_2g_1$ and $g_3g_1$
(resp. $g_4g_2$ and $g_4g_3$).
\end{lemma}

\begin{proof} Omit.
\end{proof}

Let $\cC$ be another DG algebra with a nilpotent DG ideal $\cJ
\subset \cC$. A homomorphism of DG algebras $\psi :\cB \to \cC$ such
that $\psi (\cI)\subset \cJ$ induces the functor
$$\psi ^*:\cM \cC (\cB ,\cI)\to \cM \cC (\cC ,\cJ).$$

\begin{defi} Let $\cB$ be a DG algebra and $\cR$ be an artinian DG
algebra with the maximal ideal $m\subset \cR$. Denote by $\cM\cC
_{\cR}(\cB)$ the Maurer-Cartan groupoid $\cM\cC(\cB \otimes
\cR,\cB\otimes m)$. A homomorphism of artinian DG algebras $\phi
:\cR \to \cQ$ induces the functor $\phi ^*:\cM\cC_{\cR}(\cB)\to
\cM\cC_{\cQ}(\cB)$. Thus we obtain the Maurer-Cartan 2-functor
$$\cM\cC(\cB):\dgart \to {\bf Gpd}.$$ We denote by
$\cM\cC _+(\cB)$, $\cM\cC _-(\cB)$, $\cM\cC _0(\cB)$, $\cM\cC
_{\cl}(\cB)$ the restrictions of the functor $\cM \cC (\cB)$ to
subcategories $\dgart _+$, $\dgart _-$, $\art $, $\cart $.
\end{defi}

\begin{remark} A homomorphism of DG algebras $\psi:\cC \to \cB$
induces a morphism of functors
$$\psi ^*:\cM\cC (\cC)\to \cM \cC (\cB).$$
\end{remark}

\section{Description of functors $\Def ^{\h}(E)$ and $\coDef ^{\h}(E)$}

We are going to give a description of the functor $\Def ^{\h}$ and
hence also of the functor $\coDef ^{\h}$ via the Maurer-Cartan
functor $\cM\cC$.

\begin{prop} Let $\cA$ be a DG category and $E\in \cA^0\text{-mod}$.
Denote by $\cB$ the DG algebra $\End(E)$. Then there exists an
equivalence of functors $\theta =\theta ^E: \cM\cC(\cB)\to
\Def^{\h}(E)$. (Hence also $\cM\cC(\cB)$ and $\coDef^{\h}(E)$ are
equivalent.)
\end{prop}

\begin{proof} Fix an artinian DG algebra $\cR$ with the maximal ideal $m$.
Let us define an equivalence of groupoids
$$\theta _{\cR}:\cM\cC_{\cR}(\cB)\to \Def ^{\h}_{\cR}(E).$$

Denote by $S_0=p^*E\in \cA _{\cR}^0\text{-mod}$ the trivial
$\cR$-deformation of $E$ with the differential
$d_{E,\cR}=d_E\otimes 1+1\otimes d_{\cR}$. There is a natural
isomorphism of DG algebras $\End(S_0)=\cB \otimes \cR$.

Let $\alpha \in \cM\cC(\cB\otimes m)=\cM\cC_{\cR}(\cB)$. Then in
particular $\alpha \in \End ^1(S_0)$. Hence
$d_{\alpha}:=d_{E,\cR}+\alpha$ is an endomorphism of degree 1 of the
graded module $S_0^{\gr}$. The Maurer-Cartan condition on $\alpha$
is equivalent to $d_{\alpha}^2=0$. Thus we obtain an object
$S_{\alpha}\in \cA^0_{\cR}\text{-mod}$. Clearly $i^*S_{\alpha}=E$,
so that
$$\theta _{\cR}(\alpha):=(S_{\alpha},\id)\in \Def _{\cR}^{\h}(E).$$

One checks directly that this map on objects extends naturally to a
functor $\theta _{\cR}:\cM\cC_{\cR}(\cB)\to \Def ^{\h}_{\cR}(E)$.
Indeed, maps  between Maurer-Cartan objects induce isomorphisms of
the corresponding deformations; also homotopies between such maps
become allowable homotopies between the corresponding isomorphisms.

It is clear that the functors $\theta _{\cR}$ are compatible with
the functors $\phi ^*$ induced by morphisms of DG algebras $\phi
:\cR \to \cQ$. So we obtain a morphism of functors
$$\theta :\cM\cC(\cB)\to \Def^{\h}(E).$$

It suffices to prove that $\theta _{\cR}$ is an equivalence for
each $\cR$.

\medskip

\noindent{\bf Surjective.} Let $(T,\tau )\in \Def ^{h}_{\cR}(E)$. We
may and will assume that $T^{gr}=S_0^{gr}$ and $\tau =\id$. Then
$\alpha _T:=d_T-d_{\cR,E}\in \End ^1(S_0)=(\cB\otimes \cR)^1$  is an
element in $MC(\cB \otimes \cR)$. Since $i^*\alpha _T=0$ it follows
that $\alpha _T\in \cM\cC _{\cR}(\cB)$. Thus $(T,\tau )=\theta
_{\cR}(\alpha _T)$.

\medskip

\noindent{\bf Full.} Let $\alpha, \beta \in \cM \cC _{\cR}(\cB)$. An
isomorphism between the corresponding objects $\theta
_{\cR}(\alpha)$ and $\theta _{\cR}(\beta)$ is defined by an element
$f\in \End (S_0)=(\cB \otimes \cR)$ of degree zero. The condition
$i^*f=\id _Z$ means that $f\in 1+(\cB\otimes m)$. Thus $f\in
G(\alpha ,\beta)$.

\medskip

\noindent{\bf Faithful.} Let $\alpha, \beta \in \cM \cC _{\cR}(\cB)$
and $f,g\in G(\alpha ,\beta)$. One checks directly that $f$ and $g$
are homotopic (i.e. define the same morphism in $\cM\cC_{\cR}(\cB)$)
if and only if there exists an allowable homotopy between $\theta
_{\cR}(f)$ and $\theta _{\cR}(g)$. This proves the proposition.
\end{proof}

\begin{cor} For $E\in \cA^0\text{-mod}$ the functors $\Def
^{\h}(E)$ and $\coDef ^{\h}(E)$ depend (up to equivalence) only on the DG algebra
$\End(E)$.
\end{cor}

We will prove a stronger result in Corollary 8.2 below.

\begin{example} Let $E\in \cA ^0\text{-mod}$ and denote $\cB
=\End(E)$. Consider $\cB$ as a (free) right $\cB$-module, i.e.
$\cB\in \cB ^{0}\text{-mod}$. Then $\Def ^h(\cB)\simeq \Def ^h(E)$
($\simeq \coDef ^h(\cB)\simeq \coDef ^h(E)$) because $\End
(\cB)=\End (E)=\cB$. We will describe this equivalence directly in
Section 9 below.
\end{example}

\section{Obstruction Theory}

It is convenient to describe the obstruction theory for our
(equivalent) deformation functors $\Def ^h$ and $\coDef ^h$ using
the Maurer-Cartan functor $\cM\cC(\cB)$ for a fixed  DG algebra
$\cB$.

Let $\cR$ be an artinian DG algebra with a maximal ideal $m$, such
that $m^{n+1}=0$. Put $I=m^n$, $\overline{\cR}=\cR/I$ and $\pi :\cR
\to \overline{\cR}$ the projection morphism. We have $mI=Im=0$.

Note that the kernel of the homomorphism $1 \otimes \pi:\cB \otimes
\cR\to \cB \otimes \overline{\cR }$ is the (DG) ideal $\cB\otimes
I$.  The next proposition describes the obstruction theory for
lifting objects and morphisms along the functor
$$\pi^*:\cM\cC _{\cR}(\cB)\to \cM \cC _{\overline{\cR}}(\cB).$$
It is essentially copied from [GW]. Note however a small difference
in part 3) since we do not assume that out DG algebras live in
nonnegative dimensions (and of course we work with DG algebras and
not with DG Lie algebras).

\begin{prop}
 1). There exists a map $o_2:Ob\cM\cC _{\overline{\cR}}(\cB)\to
 H^2(\cB\otimes I)$ such that $\alpha \in Ob\cM\cC
 _{\overline{\cR}}(\cB)$ is in the image of $\pi ^*$ if and only
 if $o_2(\alpha)=0$. Furthermore if $\alpha ,\beta \in Ob\cM\cC
 _{\overline{\cR}}(\cB)$ are isomorphic, then $o_2(\alpha)=0$ if
 and only if $o_2(\beta)=0$.

 2). Let $\xi \in Ob\cM\cC_{\overline{\cR}}(\cB)$. Assume that the
 fiber $(\pi ^*)^{-1}(\xi)$ is not empty. Then there exists a simply
 transitive action of the group $Z^1(\cB\otimes I)$ on the set
 $Ob(\pi ^*)^{-1}(\xi)$. Moreover the composition of the difference
 map
 $$Ob(\pi ^*) ^{-1}(\xi)\times Ob(\pi ^*) ^{-1}(\xi)\to Z^1(\cB\otimes I)$$
 with the projection
 $$Z^1(\cB\otimes I)\to H^1(\cB\otimes I)$$
 which we denote by
 $$o_1:Ob(\pi ^*)^{-1}(\xi)\times Ob(\pi ^*)^{-1}(\xi)\to H^1(\cB\otimes
 I)$$ has the following property: for $\alpha ,\beta \in Ob(\pi
 ^*)^{-1}(\xi)$ there exists a morphism $\gamma :\alpha \to \beta$
 s.t. $\pi ^*(\gamma)=\id _{\xi}$
 if and only if $o_1(\alpha ,\beta)=0$.

 3). Let $\tilde{\alpha },\tilde{\beta}\in Ob\cM\cC _{\cR}(\cB)$ be
 isomorphic objects and let $f:\alpha \to \beta$ be a morphism
 from $\alpha =\pi ^*(\tilde{\alpha})$ to $\beta =\pi
 ^*(\tilde{\beta})$. Then there is a simply transitive action of
 the group $H^0(I\otimes \cB)$ on the set $(\pi ^*)^{-1}(f)$ of
 morphisms $\tilde{f}:\tilde{\alpha}\to \tilde{\beta}$ such that
 $\pi ^*(\tilde{f})=f$. In particular the difference map
 $$o_0:(\pi ^*)^{-1}(f)\times (\pi ^*)^{-1}(f)\to H^0(\cB\otimes I)$$
 has the property: if $\tilde{f},\tilde{f}^\prime\in (\pi
  ^*)^{-1}(f)$, then $\tilde{f}=\tilde{f}^\prime$ if and only if
 $o_0(\tilde{f},\tilde{f}^\prime)=0$.
 \end{prop}

 \begin{proof} 1) Let $\alpha \in Ob\cM\cC
 _{\overline{\cR}}(\cB)=MC(\cB\otimes (m/I))$. Choose
 $\tilde{\alpha}\in (\cB\otimes m)^1$ such that $\pi
 (\tilde{\alpha})=\alpha$. Consider the element
 $$Q(\tilde{\alpha})=d\tilde{\alpha}+\tilde{\alpha}^2\in
 (\cB\otimes m)^2.$$
 Since $Q(\alpha)=0$ we have
 $Q(\tilde{\alpha})\in (\cB\otimes I)^2$. We claim that
 $dQ(\tilde{\alpha})=0$. Indeed,
 $$dQ(\tilde{\alpha})=d(\tilde{\alpha}^2)=d(\tilde{\alpha})\tilde{\alpha}-\tilde{\alpha}d(\tilde{\alpha}).$$
 We have $d(\tilde{\alpha})\equiv \tilde{\alpha}^2(mod(\cB\otimes
 I)).$ Hence
 $dQ(\tilde{\alpha})=-\tilde{\alpha}^3+\tilde{\alpha}^3=0$ (since
 $I\cdot m=0$).

 Furthermore suppose that $\tilde{\alpha}^\prime\in
 (\cB \otimes m)^1$ is another lift of $\alpha$, i.e.
 $\tilde{\alpha}^\prime -\tilde{\alpha}\in (\cB \otimes I)^1$. Then
 $$Q(\tilde{\alpha}^\prime)-Q(\tilde{\alpha})=d(\tilde{\alpha}^\prime-\tilde{\alpha})+
 (\tilde{\alpha}^\prime -\tilde{\alpha})(\tilde{\alpha}^\prime
 +\tilde{\alpha})=d(\tilde{\alpha}^\prime -\tilde{\alpha}).$$
Thus the cohomology class of the cocycle $Q(\tilde{\alpha})$ is
independent of the lift $\tilde{\alpha}$. We denote this class by
$o_2(\alpha)\in H^2(\cB\otimes I)$.

If $\alpha =\pi ^*(\tilde{\alpha})$ for some $\tilde{\alpha}\in
Ob\cM\cC _{\cR}(\cB)$, then clearly $o_2(\alpha)=0$. Conversely,
suppose $o_2(\alpha)=0$ and let $\tilde{\alpha}$ be as above. Then
$dQ(\tilde{\alpha})=d\tau$ for some $\tau \in (\cB\otimes I)^1$.
Put $\tilde{\alpha}^\prime=\tilde{\alpha}-\tau$. Then
$$Q(\tilde{\alpha}^\prime)=d\tilde{\alpha}-d\tau
+\tilde{\alpha}^2-\tilde{\alpha}\tau -\tau \tilde{\alpha}+\tau
^2=Q(\tilde{\alpha})-d\tau=0.$$

Let us prove the last assertion in 1). Assume that $\pi
^*(\tilde{\alpha})=\alpha$ and $\beta =g(\alpha)$ for some $g\in
1+(\cB\otimes m/I)^0$. Choose a lift $\tilde{g}\in 1+(\cB\otimes
m)^0$ of $g$ and put $\tilde{\beta}:=\tilde{g}(\tilde{\alpha})$.
Then $\pi ^*(\tilde{\beta})=\beta$. This proves 1).

2). Let $\alpha \in Ob(\pi ^*)^{-1}(\xi)$ and $\eta \in
Z^1(\cB\otimes I)$. Then
$$Q(\alpha +\eta)=d\alpha +d\eta +\alpha
^2 +\alpha \eta +\eta \alpha +\eta ^2=Q(\alpha)+d\eta =0.$$ So
$\alpha +\eta \in Ob(\pi ^*)^{-1}(\xi)$. This defines the action
of the group $Z^1(\cB \otimes I)$ on the set $Ob(\pi
^*)^{-1}(\xi)$.

Let $\alpha ,\beta \in Ob(\pi ^*)^{-1}(\xi)$. Then $\alpha -\beta
\in (\cB \otimes I)^1$ and
$$d(\alpha -\beta)=d\alpha -d\beta +\beta (\alpha -\beta)+(\alpha
-\beta)\beta +(\alpha -\beta )^2=Q(\alpha )-Q(\beta )=0.$$ Thus
$Z^1(\cB\otimes I)$ acts simply transitively on $Ob(\pi
^*)^{-1}(\xi)$. Now let $o_1(\alpha ,\beta)\in H^1(\cB\otimes I)$
be the cohomology class of $\alpha -\beta$. We claim that there
exists a morphism $\gamma :\alpha \to \beta$ covering $\id_{\xi}$
if and only if $o_1(\alpha ,\beta)=0$.

Indeed, let $\gamma$ be such a morphism. Then by definition the
morphisms $\pi ^*(\gamma)$ and $\id_{\xi}$ are homotopic. That is
there exists $h\in (\cB \otimes (m/I))^{-1}$ such that
$$\id_{\xi}=\pi ^*(\gamma)+d(h)+\xi h+h\xi.$$
Choose a lifting $\tilde{h}\in (\cB \otimes m)^{-1}$ on $h$ and
replace the morphism $\gamma$ by the homotopical one
$$\delta =\gamma +d(\tilde{h})+\beta \tilde{h}+\tilde{h}\alpha.$$
Thus $\delta =1+u$, where $u\in (\cB \otimes I)^0$. But then
$$\beta =\delta \alpha \delta ^{-1}+\delta d(\delta ^{-1})=\alpha
-du,$$ so that $o_1(\alpha ,\beta)=0$.

Conversely, let $\alpha -\beta=du$ for some $u\in (\cB\otimes I)^0$.
Then $\delta =1+u$ is a morphism from $\alpha $ to $\beta$ and $\pi
^*(\delta)=\id _{\xi}$. This proves 2).

3). Let us define the action of the group $Z^0(\cB\otimes I)$ on
the set $(\pi ^*)^{-1}(f)$. Let $\tilde{f}:\tilde{\alpha}\to
\tilde{\beta}$ be a lift of $f$, and $v\in Z^0(\cB\otimes I)$.
Then $\tilde{f}+v$ also belongs to $(\pi ^*)^{-1}(f)$. If $v=du$
for $u\in (\cB\otimes I )^{-1}$, then
$$\tilde{f}+v=\tilde{f}+du+\tilde{\beta} u+u \tilde{\alpha}$$
and hence morphisms $\tilde{f}$ and $\tilde{f}+v$ are homotopic.
This induces the action of $H^0(\cB\otimes I)$ on the set $(\pi
^*)^{-1}(f)$.

To show that this action is simply transitive let
$\tilde{f}^\prime :\tilde{\alpha }\to \tilde{\beta}$ be another
morphism in $(\pi ^*)^{-1}(f)$. This means by definition that
there exists $h\in (\cB\otimes (m/I))^{-1}$ such that
$$f=\pi ^*(\tilde{f}^\prime)+dh+\beta h+h\alpha.$$
Choose a lifting $\tilde{h}\in (\cB \otimes m)^{-1}$ of $h$ and
replace $\tilde{f}^\prime$ by the homotopical morphism
$$\tilde{g}=\tilde{f}^\prime+d\tilde{h}+\tilde{\beta}\tilde{h}+
\tilde{h}\tilde{\alpha}.$$ Then $\tilde{g}=\tilde{f}+v$ for $v\in
(\cB \otimes I)^0$. Since $\tilde{f}, \tilde{g}:\tilde{\alpha}\to
\tilde{\beta}$ we must have that $v\in Z^0(\cB \otimes I)$. This
shows the transitivity and proves 3).
\end{proof}

\section{Invariance theorem and its implications}

\begin{theo} Let $\phi :\cB
\to \cC$ be a quasi-isomorphism of  DG algebras. Then the induced
morphism of functors $$\phi ^*:\cM\cC(\cB )\to \cM\cC(\cC)$$ is an
equivalence.
\end{theo}

\begin{proof}
The proof is almost the same as that of Theorem 2.4 in [GW]. We
present it for reader's convenience and also because of the slight
difference in language: in [GW] they work with DG Lie algebras as
opposed to DG algebras.

Fix an artinian DG algebra $\cR$ with the maximal ideal $m\subset
\cR$, such that $m^{n+1}=0$. We prove that
$$\phi ^*:\cM\cC _{\cR}(\cB)\to \cM\cC _{\cR}(\cC)$$
is an equivalence by induction on $n$. If $n=o$, then both
groupoids contain one object and one morphism, so are equivalent.
Let $n>0$. Put $I=m^n$ with the projection $\pi :\cR \to \cR
/I=\overline{\cR}$. We have the commutative functorial diagram
$$\begin{array}{ccc}
\cM\cC _{\cR}(\cB) & \stackrel{\phi ^*}{\rightarrow} & \cM\cC
_{\cR}(\cC)\\
\pi ^*\downarrow & & \downarrow \pi ^*\\
\cM \cC _{\overline{\cR}}(\cB) & \stackrel{\phi ^*}{\rightarrow} &
\cM \cC _{\overline{\cR}}(\cC).
\end{array}$$
By induction we may assume that the bottom functor is an
equivalence. To prove the same about the top one we need to
analyze the fibers of the functor $\pi _*$. This has been done by
the obstruction theory.

We will prove that the functor
$$\phi ^*:\cM\cC _{\cR}(\cB)\to \cM\cC _{\cR}(\cC)$$
is surjective on the isomorphism classes of objects, is full and
is faithful.

\medskip

\noindent{\bf Surjective on isomorphism classes.} Let $\beta \in
Ob\cM\cC _{\cR}(\cC)$. Then $\pi ^*\beta \in Ob\cM\cC
_{\overline{\cR} }(\cC)$. By the induction hypothesis there exists
$\alpha ^\prime \in Ob\cM\cC _{\overline{\cR}}(\cC)$ and an
isomorphism $g: \phi ^*\alpha ^\prime \to \pi ^* \beta$. Now
$$H^2(\phi)o_2(\alpha ^\prime)=o_2(\phi ^*\alpha ^\prime )=
o_2(\pi ^*\beta )=0.$$ Hence $o_2(\alpha ^\prime)=0$, so there
exists $\tilde{\alpha }\in Ob\cM\cC _{\cR}(\cB)$ such that $\pi
^*\tilde{\alpha}=\alpha ^\prime$, and hence
$$\phi ^*\pi ^*\tilde{\alpha}=\pi ^*\phi ^*\tilde{\alpha}=\phi
^*\alpha ^\prime.$$

Choose a lift $\tilde{g}\in 1+(\cC\otimes m)^0$ of $g$ and put
$\tilde{\beta}=\tilde{g}^{-1}(\beta)$. Then
$$\pi ^*(\tilde{\beta})=\pi ^*(\tilde{g}^{-1}(\beta))=g^{-1}\pi
^*\beta=\phi ^*\alpha ^\prime.$$

The obstruction to the existence of an isomorphism $\phi
^*\tilde{\alpha}\to \tilde{\beta}$ covering $\id_{\pi^*(\alpha
^\prime)}$ is an element $o_1(\phi ^*(\tilde{\alpha}),
\tilde{\beta})\in H^1(\cC\otimes I)$. Since $H^1(\phi)$ is
surjective there exists a cocycle $u\in Z^1(\cB\otimes I)$ such that
$H^1(\phi)[u]=o_1(\phi ^*(\tilde{\alpha}),\tilde{\beta})$. Put
$\alpha =\tilde{\alpha}-u\in Ob\cM\cC _{\cR}(\cB)$. Then
$$\begin{array}{rcl}o_1(\phi^*\alpha ,\tilde{\beta}) & = &
o_1(\phi ^*\alpha ,\phi ^*\tilde{\alpha})+o_1(\phi
^*\tilde{\alpha},\tilde{\beta})\\
& = & H^1(\phi)o_1(\alpha
,\tilde{\alpha})+o_1(\phi^*\tilde{\alpha},\tilde{\beta})\\
& = & -H^1(\phi)[u]+o_1(\phi ^*\tilde{\alpha},\beta)=0
\end{array}
$$
This proves the surjectivity of $\phi ^*$ on isomorphism classes.

\medskip

\noindent{\bf Full.} Let $f:\phi ^*\alpha _1\to \phi ^*\alpha _2$
be a morphism in $\cM\cC _{\cR}(\cC)$. Then $\pi ^*f$ is a
morphism in $\cM\cC _{\overline{\cR}}(\cC)$:
$$\pi ^*(f):\phi^*\pi^*\alpha _1\to \phi ^*\pi ^*\alpha _2.$$

By induction hypothesis there exists $g:\pi ^*\alpha _1\to \pi
^*\alpha _2$ such that $\phi ^*(g)=\pi^*(f)$. Let $\tilde{g}\in
1+(\cC \otimes m)^0$ be any lift of $g$. Then $\pi
^*(\tilde{g}\alpha _1)=\pi ^*\alpha _2$. The obstruction to the
existence of a morphism $\gamma :\tilde{g}\alpha _1\to \alpha _2$
covering $\id_{\pi ^*\alpha _2}$ is an element $o_1(\tilde{g}\alpha
_1,\alpha _2)\in H^1(\cB\otimes I)$. By assumption $H^1(\phi)$ is an
isomorphism and we know that
$$H^1(\phi)(o_1(\tilde{g}\alpha
_1,\alpha _2))=o_1(\phi^*\tilde{g}\alpha _1,\phi ^*\alpha _2)=0,$$
since the morphism $f\cdot (\phi ^*\tilde{g})^{-1}$ is covering the
identity morphism $\id_{\pi ^*\phi ^*\alpha _2}$. Thus
$o_1(\tilde{g}\alpha _1,\alpha _2)=0$ and $\gamma $ exists. Then
$\gamma \cdot \tilde{g}:\alpha _1\to \alpha _2$ is covering $g:\pi
^*\alpha _1\to \pi^*\alpha _2$. Hence both morphisms $\phi ^*(\gamma
\cdot \tilde{g})$ and $f$ are covering $\pi ^*(f)$. The obstruction
to their equality is an element $o_0(\phi ^*(\gamma \cdot
\tilde{g}),f)\in H^0(\cC \otimes I)$. Let $u\in Z^0(\cC\otimes I)$
be a representative of the inverse image under $H^0(\phi)$ of this
element. Then $\phi ^*(\gamma \cdot \tilde{g}+u)=f$.

\medskip

\noindent{\bf Faithful.} Let $\gamma _1,\gamma _2:\alpha _1\to
\alpha _2$ be morphisms in $\cM\cC _{\cR}(\cB)$ with $\phi^*\gamma
_1=\phi^*\gamma _2$. Then $\phi ^*\pi ^*\gamma _1=\phi ^*\pi
^*\gamma _2$. By the induction hypothesis $\pi ^*\gamma _1=\pi
^*\gamma _2$, so the obstruction $o_0(\gamma _1,\gamma _2)\in
H^0(\cB\otimes I)$ is defined. Now $H^0(\phi )o_0(\phi ^*\gamma
_1,\phi ^*\gamma _2)=0$. Since $H^0(\phi)$ is injective it follows
that $\gamma _1=\gamma _2$. This proves the theorem.
\end{proof}

\begin{cor} The homotopy (co-) deformation functor of $E\in \cA
^{0}\text{-mod}$ depends only on the quasi-isomorphism class of
the DG algebra $\End (E)$.
\end{cor}

\begin{proof} This follows from Theorem 8.1 and Proposition 6.1.
\end{proof}

The next proposition provides two examples of this situation. It was
communicated to us by Bernhard Keller.

\begin{prop} (Keller)  a) Assume that $E^\prime \in \cA
^{0}\text{-mod}$ is homotopy equivalent to $E$. Then the DG algebras
$\End (E)$ and $\End(E^\prime )$ are canonically quasi-isomorphic.

b) Let $P\in \cP(\cA)$ and $I\in \cI(\cA)$  be  quasi-isomorphic.
Then the DG algebras $\End(P)$ and $\End(I)$ are canonically
quasi-isomorphic.
\end{prop}

\begin{proof} a) Let $g:E\to E^\prime$ be a homotopy equivalence.
Consider its cone $C(g)\in \cA ^{0}\text{-mod}$. Let $\cC \subset
\End (C(g))$ be the DG subalgebra consisting of endomorphisms
which leave $E^\prime $ stable. There are natural projections
$p:\cC\to \End(E^\prime)$ and $q:\cC\to \End (E)$. We claim that
$p$ and $q$ are quasi-isomorphisms. Indeed, $\Ker(p)$ (resp. $\Ker
(q)$) is the complex $\Hom (E,C(g))$ (resp. $\Hom
(C(g),E^\prime)$). These complexes are acyclic, since $g$ is a
homotopy equivalence.

b) The proof is similar. Let $f:P\to I$ be a quasi-isomorphism. Then
the cone $C(f)$ is acyclic. We consider the DG subalgebra $\D\subset
\End (C(f))$ which leaves $I$ stable. Then $\D$ is quasi-isomorphic
to $\End(I)$ and $\End(P)$ because the complexes $\Hom (P,C(f))$ and
$\Hom (C(f),I)$ are acyclic.
\end{proof}

\begin{cor} a) If DG $\cA ^0$-modules $E$ and $E^\prime$ are
homotopy equivalent then the functors $\Def ^{\h}(E)$, $\coDef
^{\h}(E)$, $\Def ^{\h}(E^\prime)$, $\coDef ^{\h}(E^\prime)$ are
canonically equivalent.

b) Let $P\to  I$ be a quasi-isomorphism between $P\in \cP(\cA)$ and
$I\in \cI(\cA)$. Then the functors $\Def ^{\h}(P)$, $\coDef
^{\h}(P)$, $\Def ^{\h}(I)$, $\coDef ^{\h}(I)$ are canonically
equivalent.
\end{cor}

\begin{proof} Indeed, this follows from Proposition 8.3 and
Corollary 8.2.
\end{proof}

Actually, one can prove a more precise statement.

\begin{prop} Fix an artinian DG algebra $\cR$.

a) Let $g:E\to E^\prime $ be a homotopy equivalence of DG $\cA
^0$-modules. Assume that $(V,\id)\in \Def ^{\h}_{\cR}(E)$ and
$(V^\prime,\id)\in \Def ^{\h}_{\cR}(E^\prime)$ are objects that
correspond to each other via the equivalence $\Def
^{\h}_{\cR}(E)\simeq \Def ^{\h}(E^\prime)$ of Corollary 8.4. Then
there exists a homotopy equivalence $\tilde{g}:V\to V^\prime$ which
extends $g$, i.e. $i^*\tilde{g}=g$. Similarly for the objects of
$\coDef ^{\h}_{\cR}$ with $i^!$ instead of $i^*$.

b) Let $f:P\to I$ be a quasi-isomorphism with $P\in \cP(\cA)$, $I\in
\cI(\cA)$.  Assume that $(S,\id)\in \Def ^{\h}_{\cR}(E)$ and
$(T,\id)\in \Def ^{\h}_{\cR}(E^\prime)$ are objects that correspond
to each other via the equivalence $\Def ^{\h}_{\cR}(P)\simeq \Def
^{\h}_{\cR}(I)$ of Corollary 8.4. Then there exists a
quasi-isomorphism $\tilde{f}:S\to T$ which extends $f$, i.e.
$i^*\tilde{f}=f$. Similarly for the objects of $\coDef^{\h}$ with
$i^!$ instead of $i^*$.
\end{prop}

\begin{proof} a)
 Consider the DG algebra
 $$\cC   \subset \End(C(g))$$ as in the proof of Proposition 8.3.
 We proved there that the natural projections $\End
 (E)\leftarrow \cC \rightarrow \End(E^\prime)$ are
 quasi-isomorphisms. Hence the induced functors between groupoids
 $\cM \cC_{\cR}(\End (E))\leftarrow \cM \cC _{\cR}(\cC)\rightarrow \cM \cC_{\cR}(\End
 (E^\prime))$ are equivalences by Theorem 8.1. Using Proposition 6.1
 we may and will assume that deformations $(V, \id)$, $(V^\prime
 ,\id )$ correspond to elements $\alpha _E\in \cM \cC_{\cR}(\End
 (E))$, $\alpha _{E^\prime}\in \cM \cC_{\cR}(\End (E^\prime))$ which
 come from the same element $\alpha \in \cM \cC_{\cR}(\cC)$.

 Consider the DG modules $E\otimes \cR$, $E^\prime \otimes \cR$
 with the differentials $d_E\otimes 1+1\otimes d_{\cR}$ and
 $d_{E^\prime}\otimes 1+1\otimes d_{\cR}$ respectively and the morphism $g\otimes 1:E\otimes \cR \to
 E^\prime \otimes \cR$. Then
$$\cC\otimes \cR =\left(  \begin{array}{cc}
\End(E^\prime\otimes \cR) & \Hom (E\otimes \cR,E^\prime \otimes \cR) \\
0 & \End(E\otimes \cR)
\end{array}  \right)
  \subset \End(C(g\otimes 1)),$$
  and
 $$\alpha =\left( \begin{array}{cc}
 \alpha _{E^\prime} & t\\
 0                  & \alpha _{E}
 \end{array}\right).$$

Recall that the differential in the DG module $C(g\otimes 1)$
 is of the form $(d_{E^\prime}\otimes 1, d_E[1]\otimes 1+g[1]\otimes
 1)$. The element $\alpha$ defines a new differential $d_{\alpha}$ on
$C(g\otimes 1)$ which is $(d_{E^\prime}\otimes 1+\alpha _{E^\prime},
(d_E[1]\otimes 1+\alpha _E)+(g[1]\otimes
 1+t))$. The fact that $d_{\alpha }^2=0$ implies that $\tilde{g}:=g\otimes
 1+t[-1]:V\to V^\prime$ is a closed morphism of degree zero and hence
 the DG module $C(g\otimes 1)$ with the differential $d_\alpha$ is
 the cone $C(\tilde{g})$ of this morphism.

 Clearly,
 $i^*\tilde{g}=g$ and it remains to prove that $\tilde{g}$ is a
 homotopy equivalence. This in turn is equivalent to the
 acyclicity of the DG algebra $\End(C(\tilde{g}))$. But recall
 that the differential in $\End(C(\tilde{g}))$ is an "$\cR$-deformation" of the
 differential in the  DG algebra
 $\End(C(g))$ which is acyclic, since $g$ is a homotopy
 equivalence. Therefore $\End(C(\tilde{g}))$ is also acyclic. This
 proves the first statement in a). The last statement follows by the
 equivalence of groupoids $\Def ^{\h}_{\cR}\simeq \coDef
 ^{\h}_{\cR}$ (Proposition 4.7).

The proof of b) is similar: exactly in the same way we construct a
closed morphism of degree zero $\tilde{f}:S\to T$ which extends $f$.
Then $\tilde{f}$ is a quasi-isomorphism, because $f$ is such.
\end{proof}

\begin{prop} a) Let $F:\cA \to \cC$ be a DG functor which induces an
equivalence of derived categories $\bL F^*:D(\cA)\to D(\cC)$. (For
example, this is the case if $F$ induces a quasi-equivalence
$F^{\pre-tr}:\cA ^{\pre-tr}\to \cC ^{\pre-tr}$ (Corollary 3.15)).

a) Let $P\in \cP (\cA)$. Then the map of DG algebras $F^*:\End
(P)\to \End (F^*(P))$ is a quasi-isomorphism. Hence the deformation
functors $\Def ^{\h}$ and $\coDef ^{\h}$ of $P$ and $F^*(P)$ are
equivalent.

b) Let $I\in \cI(\cA)$. Then the map of DG algebras $F^!:\End(I)\to
\End (F^!(I))$ is a quasi-isomorphism. Hence the deformation
functors $\Def ^{\h}$ and $\coDef ^{\h}$ of $I$ and $F^!(I)$ are
equivalent.
\end{prop}

\begin{proof} a) By Lemma 3.6 we have $F^*(P)\in \cP (\cC)$. Hence the assertion follows from
Theorems 3.1 and 8.1.

b) The functor $\bR F^!:D(\cA)\to D(\cC)$ is also an equivalence
because of adjunctions $(F_*,\bR F^!),(\bL F^* ,F_*)$. Also
$F^!(I)\in \cI(\cC)$ (Lemma 3.6). Hence the assertion follows from
Theorems 3.1 and 8.1.
\end{proof}

\begin{theo} Let $\phi :\cB \to \cC$ be a morphism of  DG
algebras, such that the induced map $H^i(\phi):H^i(\cB)\to H^i(\cC)$
is an isomorphism for $i\geq 0$. Then the induced morphism of
functors
$$\phi ^*:\cM\cC _-(\cB)\to \cM\cC _-(\cC)$$
is an equivalence.
\end{theo}

\begin{proof} Let $V^{\cdot}$ be a complex such that $V^i=0$ for
$i>0$. Then $H^i(\phi \otimes 1):H^i(\cB \otimes V^{\cdot})\to
H^i(\cC \otimes V^{\cdot})$ is an isomorphism for $i\geq 0$. Now the
proof is the same as that of Theorem 8.1. Indeed, in that proof we
only considered cohomology groups $H^0$, $H^1$ and $H^2$. Hence the
same reasoning applies, since we restrict ourselves to artinian DG
algebras $\cR \in \dgart _-$.
\end{proof}

\section{Direct relation between functors $\Def ^{\h}(F)$ and $\Def
^{\h}(\cB)$ ($\coDef ^{\h}(F)$ and $\coDef ^{\h}(\cB)$)}

\subsection{DG functor $\Sigma$} Let $F\in \cA ^0\text{-mod}$ and
put $\cB =\End(F)$. Recall the DG functor from Example 3.14
$$\Sigma =\Sigma ^F:\cB ^0\text{-mod}\to \cA ^0\text{-mod},\quad
\Sigma (M)=M\otimes _{\cB}F.$$ For each artinian DG algebra $\cR$ we
obtain the corresponding DG functor
$$\Sigma _{\cR}:(\cB\otimes \cR) ^0\text{-mod}\to \cA _{\cR}^0\text{-mod},\quad
\Sigma _{\cR}(M)=M\otimes _{\cB}F.$$

\begin{lemma} The DG functors $\Sigma _{\cR}$ have the following
properties.

a) If a DG $(\cB\otimes \cR) ^0$-module $M$ is graded $\cR$-free
(resp. graded $\cR$-cofree), then so is the DG $\cA_{\cR} ^0$-module
$\Sigma _{\cR}(M)$.

b) Let $\phi :\cR \to \cQ$ be a homomorphism of artinian DG
algebras. Then there are natural isomorphisms of DG functors
$$\Sigma _{\cQ}\cdot \phi ^*=\phi ^*\cdot \Sigma _{\cR}, \quad
\Sigma _{\cR}\cdot \phi _*=\phi _*\cdot \Sigma _{\cQ}.$$ In
particular,
$$\Sigma \cdot i ^*=i ^*\cdot \Sigma _{\cR}.$$

c) There is a natural isomorphism of DG functors
$$\Sigma _{\cQ}\cdot \phi ^!=\phi ^!\cdot \Sigma _{\cR}$$
on the full DG subcategory of DG $(\cB\otimes \cR)^0$-modules $M$
such that $M^{\gr}\simeq M_1^{\gr}\otimes M_2^{\gr}$ for a
$\cB^0$-module $M_1$ and an $\cR ^0$-module $M_2$. (This subcategory
includes in particular graded $\cR$-cofree modules.) Therefore
$$\Sigma \cdot i ^!=i ^!\cdot \Sigma _{\cR}$$
on this subcategory.

d) For a graded $\cR$-free DG $(\cB\otimes \cR) ^0$-module $M$ there
is a functorial isomorphism
$$\Sigma _{\cR}(M\otimes _{\cR} \cR ^*)=\Sigma _{\cR}(M)\otimes _{\cR}\cR ^*$$
\end{lemma}

\begin{proof} The only nontrivial assertion is c). For any DG $(\cB
\otimes \cR)^0$-module $M$ there is a natural closed morphism of
degree zero of DG $\cA ^0_{\cR}$-modules
$$\gamma _M:\Hom _{\cR ^0}(Q, M)\otimes _{\cB}F\to \Hom _{\cR ^0}(Q,
M\otimes _{\cB}F), \quad \gamma(g\otimes
f)(q)=(-1)^{\bar{f}\bar{q}}g(q)\otimes f.$$ Since $\cQ$ is a finite
$\cR ^0$-module $\gamma _M$ is an isomorphism if $M^{\gr}\simeq
M_1^{\gr}\otimes M_2^{\gr}$ for a $\cB^0$-module $M_1$ and an $\cR
^0$-module $M_2$.
\end{proof}

\begin{prop} a) For each artinian DG algebra $\cR$ the DG functor
$\Sigma _{\cR}$ induces functors between groupoids
$$\Def ^{\h}(\Sigma _{\cR}):\Def ^{\h}_{\cR}(\cB)\to \Def ^{\h}_{\cR
}(F),$$
$$\coDef ^{\h}(\Sigma _{\cR}):\coDef ^{\h}_{\cR}(\cB)\to \coDef ^{\h}_{\cR
}(F),$$

b) The collection of DG functors $\{\Sigma _{\cR}\}_{\cR}$ defines
natural transformations
$$\Def ^{\h}(\Sigma ):\Def ^{\h}(\cB)\to \Def ^{\h}(F),$$
$$\coDef ^{\h}(\Sigma ):\Def ^{\h}(\cB)\to \Def ^{\h}(F).$$

c) The morphism $\Def ^{\h}(\Sigma )$ is compatible with the
equivalence $\theta $ of Proposition 6.1. That is the functorial
diagram
$$\begin{array}{ccc}
\cM \cC (\cB) & = & \cM \cC (\cB)\\
\theta ^{\cB} \downarrow & & \downarrow \theta ^{F}\\
\Def ^{\h}(\cB) & \stackrel{\Def^{\h}(\Sigma)}{\rightarrow} & \Def
^{\h}(F)
\end{array}
$$
is commutative.

d) The morphisms $\Def ^{\h}(\Sigma )$ and $\coDef ^{\h}(\Sigma )$
are compatible with the equivalence $\delta$ of Proposition 4.7.
That is the functorial diagram
$$\begin{array}{ccc}
\Def ^{\h}(\cB) & \stackrel{\Def^{\h}(\Sigma)}{\rightarrow} & \Def
^{\h}(F)\\
\delta ^{\cB} \downarrow & & \downarrow \delta ^{F}\\
\coDef ^{\h}(\cB) & \stackrel{\coDef^{\h}(\Sigma)}{\rightarrow} &
\coDef ^{\h}(F)
\end{array}
$$
is commutative.

e) The natural transformations $\Def ^{\h}(\Sigma )$ and $\coDef
^{\h}(\Sigma )$ are equivalences, i.e. for each $\cR$ the functors
$\Def ^{\h}(\Sigma _{\cR})$ and $\coDef ^{\h}(\Sigma _{\cR})$ are
equivalences.
\end{prop}

\begin{proof} a) and b) follow from parts a),b),c) of Lemma 9.1; c) is obvious;
d) follows from part d) of Lemma 9.1; e) follows from c) and d).
\end{proof}

\subsection{DG functor $\psi^*$} Let $\psi :\cC \to \cB$ be a
homomorphism of DG algebras. Recall the corresponding DG functor
$$\psi ^*:\cC ^0\text{-mod}\to \cB ^0\text{-mod},\quad \psi
^*(M)=M\otimes _{\cC}\cB.$$ For each artinian DG algebra $\cR$ we
obtain a similar DG functor
$$\psi ^*_{\cR}:(\cC \otimes \cR) ^0\text{-mod}\to (\cB\otimes \cR) ^0\text{-mod},\quad \psi
^*(M)=M\otimes _{\cC}\cB.$$

The next lemma and proposition are complete analogues of Lemma 9.1
and Proposition 9.2.

\begin{lemma} The DG functors $\psi ^* _{\cR}$ have the following
properties.

a) If a DG $(\cC\otimes \cR) ^0$-module $M$ is graded $\cR$-free
(resp. graded $\cR$-cofree), then so is the DG $(\cB\otimes \cR)
^0$-module $\psi ^*_{\cR}(M)$.

b) Let $\phi :\cR \to \cQ$ be a homomorphism of artinian DG
algebras. Then there are natural isomorphisms of DG functors
$$\psi^* _{\cQ}\cdot \phi ^*=\phi ^*\cdot \psi^* _{\cR}, \quad
\psi^* _{\cR}\cdot \phi _*=\phi _*\cdot \psi ^* _{\cQ}.$$ In
particular,
$$\psi^* \cdot i ^*=i ^*\cdot \psi^* _{\cR}.$$

c) There is a natural isomorphism of DG functors
$$\psi^* _{\cQ}\cdot \phi ^!=\phi ^!\cdot \psi^* _{\cR}$$
on the full DG subcategory of  DG $(\cC\otimes \cR)^0$-modules $M$
such that $M^{\gr}\simeq M_1^{\gr}\otimes M_2^{\gr}$ for a
$\cC^0$-module $M_1$ and an $\cR ^0$-module $M_2$. (This subcategory
includes in particular graded $\cR$-cofree modules.) Therefore
$$\psi^* \cdot i ^!=i ^!\cdot \psi^* _{\cR}$$
on this subcategory.

d) For a graded $\cR$-free DG $(\cC\otimes \cR) ^0$-module $M$ there
is a functorial isomorphism
$$\psi^* _{\cR}(M\otimes _{\cR} \cR ^*)=\psi^* _{\cR}(M)\otimes _{\cR}\cR ^*$$
\end{lemma}

\begin{proof} This is a special case of Lemma 9.1.
\end{proof}

\begin{prop} a) For each artinian DG algebra $\cR$ the DG functor
$\psi^* _{\cR}$ induces functors between groupoids
$$\Def ^{\h}(\psi^* _{\cR}):\Def ^{\h}_{\cR}(\cC)\to \Def ^{\h}_{\cR
}(\cB),$$
$$\coDef ^{\h}(\psi^* _{\cR}):\coDef ^{\h}_{\cR}(\cC)\to \coDef ^{\h}_{\cR
}(\cB),$$

b) The collection of DG functors $\{\psi^* _{\cR}\}_{\cR}$ defines
natural transformations
$$\Def ^{\h}(\psi^* ):\Def ^{\h}(\cC)\to \Def ^{\h}(\cB),$$
$$\coDef ^{\h}(\psi^* ):\Def ^{\h}(\cC)\to \Def ^{\h}(\cB).$$

c) The morphism $\Def ^{\h}(\psi^* )$ is compatible with the
equivalence $\theta $ of Proposition 6.1. That is the functorial
diagram
$$\begin{array}{ccc}
\cM \cC (\cC) & \stackrel{\psi^*}{\rightarrow} & \cM \cC (\cB)\\
\theta ^{\cC} \downarrow & & \downarrow \theta ^{\cB}\\
\Def ^{\h}(\cC) & \stackrel{\Def^{\h}(\psi^*)}{\rightarrow} & \Def
^{\h}(\cB)
\end{array}
$$
is commutative.

d) The morphisms $\Def ^{\h}(\psi^* )$ and $\coDef ^{\h}(\psi^* )$
are compatible with the equivalence $\delta$ of Proposition 4.7.
That is the functorial diagram
$$\begin{array}{ccc}
\Def ^{\h}(\cC) & \stackrel{\Def^{\h}(\psi^*)}{\rightarrow} & \Def
^{\h}(\cB)\\
\delta ^{\cC} \downarrow & & \downarrow \delta ^{\cB}\\
\coDef ^{\h}(\cC) & \stackrel{\coDef^{\h}(\psi^*)}{\rightarrow} &
\coDef ^{\h}(\cB)
\end{array}
$$
is commutative.

e) Assume that $\psi$ is a quasi-isomorphism. Then the natural
transformations $\Def ^{\h}(\psi^* )$ and $\coDef ^{\h}(\psi^* )$
are equivalences, i.e. for each $\cR$ the functors $\Def
^{\h}(\psi^* _{\cR})$ and $\coDef ^{\h}(\psi^* _{\cR})$ are
equivalences.

f) Assume that the induced map $H^i(\psi):H^i(\cB)\to H^i(\cC)$ is
an isomorphism for $i\geq 0$.  Then the natural transformations
$\Def ^{\h}_-(\psi^* )$ and $\coDef ^{\h}_-(\psi^* )$ are
equivalences, i.e. for each $\cR \in \dgart _-$ the functors $\Def
^{\h}(\psi^* _{\cR})$ and $\coDef ^{\h}(\psi^* _{\cR})$ are
equivalences.
\end{prop}

\begin{proof} a) and b) follow from parts a),b),c) of Lemma 9.3; c)
is obvious; d) follows from part d) of Lemma 9.3; e) follows from
c),d) and Theorem 8.1; f) follows from c),d) and Theorem 8.7.
\end{proof}

Later we will be especially interested in the following example.

\begin{lemma}(Keller). a) Assume that the DG algebra $\cB$ satisfies
the following conditions: $H^i(\cB)=0$ for $i<0$, $H^0(\cB)=k$
(resp. $H^0(\cB)=k$). Then there exists a DG subalgebra $\cC\subset
\cB$ with the properties: $\cC^i=0$ for $i<0$, $\cC^0=k$, and the
embedding $\psi:\cC\hookrightarrow \cB$ is a quasi-isomorphism
(resp. the induced map $H^i(\psi):H^i(\cC)\to H^i(\cB)$ is an
isomorphism for $i\geq 0$).
\end{lemma}

\begin{proof} Indeed, put $\cC^0=k$, $\cC^1=K\oplus L$, where
$d(K)=0$ and $K$ projects isomorphically to $H^1(\cB)$, and
$d:L\stackrel{\sim}{\to}d(\cB ^1)\subset \cB ^2$. Then take $\cC
^i=\cB^i$ for $i\geq 2$ and $\cC ^i=0$ for $i<0$.
\end{proof}

\section{The derived deformation and co-deformation functors}

\subsection{The functor $\Def (E)$} Fix a DG category $\cA$ and an object $E\in \cA^0\text{-mod}$. We
are going to define a 2-functor $\Def(E)$ from the category $\dgart$
to the category ${\bf Gpd}$ of groupoids. This functor assigns to a
DG algebra $\cR$ the groupoid $\Def _{\cR}(E)$ of $\cR$-deformations
of $E$ in the {\it derived} category $D(\cA)$.

\begin{defi} Fix an artinian DG algebra $\cR$. An object of the
groupoid $\Def _{\cR}(E)$ is a pair $(S,\sigma)$, where $S\in
D(\cA _{\cR})$ and $\sigma$ is an isomorphism (in $D(\cA)$)
$$\sigma :\bL i^*S\to E.$$
A morphism $f:(S,\sigma)\to (T,\tau)$ between two
$\cR$-deformations of $E$ is an isomorphism (in $D(\cA _{\cR})$)
$f:S\to T$, such that
$$\tau \cdot \bL i^*(f)=\sigma.$$
This defines the groupoid $\Def _{\cR}(E)$. A homomorphism of
artinian DG algebras $\phi:\cR \to \cQ$ induces the functor
$$\bL\phi ^*:\Def _{\cR}(E)\to \Def _{\cQ}(E).$$
Thus we obtain a 2-functor
$$\Def (E):\dgart \to {\bf Gpd}.$$
\end{defi}

We call $\Def (E)$ the functor of derived deformations of $E$.

\begin{remark} A quasi-isomorphism $\phi :\cR\to \cQ$ of artinian
DG algebras induces an equivalence of groupoids
$$\bL\phi ^*:\Def _{\cR}(E)\to \Def _{\cQ}(E).$$ Indeed,
$\bL\phi ^*:D(\cA _{\cR})\to D(\cA _{\cQ})$ is an equivalence of
categories (Proposition 3.7) which commutes with the functor $\bL
i^*$.
\end{remark}

\begin{remark} A quasi-isomorphism $\delta:E_1\to E_2$ of
DG $\cA^0$-modules induces an equivalence of functors
$$\delta _*:\Def (E_1)\to \Def (E_2)$$
by the formula $\delta _*(S,\sigma)=(S,\delta \cdot \sigma)$.
\end{remark}

\begin{prop} Let $F:\cA \to \cA ^\prime$ be a DG functor which
induces a quasi-equivalence  $F^{\pre-tr}:\cA ^{\pre-tr}\to
\cA^{\prime \pre-tr}$. Then for any $E\in D(\cA )$ the deformation
functors $\Def (E)$ and $\Def (\bL F^*(E))$ are canonically
equivalent. (Hence also $\Def (F_*(E^\prime))$ and $\Def (E^\prime)$
are equivalent for any $E^\prime \in D(\cA ^\prime)$).
\end{prop}

\begin{proof} For any artinian DG algebra $\cR$ the functor $F$
induces a commutative functorial diagram
$$\begin{array}{ccc}
D(\cA _{\cR}) & \stackrel{\bL (F \otimes \id)^*}{\longrightarrow} &
D(\cA _{\cR}^\prime)\\
\downarrow \bL i^* & & \downarrow \bL i^*\\
D(\cA ) & \stackrel{\bL F^*}{\longrightarrow} & D(\cA ^\prime)
\end{array}
$$
where $\bL F^*$ and $\bL (F\otimes \id)^*$ are equivalences by
Corollary 3.15.  The horizontal arrows define a functor
$F^*_{\cR}:\Def _{\cR}(E)\to \Def _{\cR}(\bL F^*(E))$. Moreover
these functors are compatible with the functors $\bL \phi ^*:\Def
_{\cR}\to \Def _{\cQ}$ induced by morphisms $\phi :\cR \to \cQ$ of
artinian DG algebras. So we get the morphism $F^*:\Def (E)\to \Def
(\bL F^*(E))$ of 2-functors. It is clear that for each $\cR$ the
functor $F^*_{\cR}$ is an equivalence. Thus $F^*$ is also such.
\end{proof}

\begin{example} Suppose that $\cA ^\prime$ is a pre-triangulated DG category (so that
the homotopy category $\Ho (\cA ^\prime)$ is triangulated). Let
$F:\cA \hookrightarrow \cA ^\prime$ be an embedding of a full DG
subcategory so that the triangulated category $\Ho (\cA ^\prime)$ is
generated by the collection of objects $Ob\cA$. Then the assumption
of the previous proposition holds.
\end{example}

\begin{remark} In the definition of the functor $\Def (E)$ we
could work with the homotopy category of h-projective
DG modules instead of the derived category. Indeed, the functors
$i^*$ and $\phi ^*$ preserve h-projective DG modules.
\end{remark}

\begin{defi} Denote by $\Def _+(E)$, $\Def _-(E)$, $\Def
_0(E)$, $\Def _{\cl}(E)$ the restrictions of the functor $\Def (E)$
to subcategories $\dgart _+$, $\dgart _-$, $\art$, $\cart$
respectively.
\end{defi}

\subsection{The functor $\coDef (E)$}
Now we define the 2-functor $\coDef (E)$ of {\it derived
co-deformations} in a similar way replacing everywhere the functors
$(\cdot )^*$ by $(\cdot )^!$.

\begin{defi} Fix an artinian DG algebra $\cR$. An object of the
groupoid $\coDef _{\cR}(E)$ is a pair $(S,\sigma)$, where $S\in
D(\cA _{\cR})$ and $\sigma$ is an isomorphism (in $D(\cA)$)
$$\sigma :E\to \bR i^!S.$$
A morphism $f:(S,\sigma)\to (T,\tau)$ between two
$\cR$-deformations of $E$ is an isomorphism (in $D(\cA _{\cR})$)
$f:S\to T$, such that
$$ \bR i^!(f)\cdot \sigma=\tau.$$
This defines the groupoid $\coDef _{\cR}(E)$. A homomorphism of
artinian DG algebras $\phi:\cR \to \cQ$ induces the functor
$$\bR\phi ^!:\coDef _{\cR}(E)\to \coDef _{\cQ}(E).$$
Thus we obtain a 2-functor
$$\coDef (E):\dgart \to {\bf Gpd}.$$
\end{defi}

We call $\coDef (E)$ the functor of derived co-deformations of $E$.

\begin{remark} A quasi-isomorphism $\phi :\cR\to \cQ$ of artinian
DG algebras induces an equivalence of groupoids
$$\bR\phi ^!:\coDef _{\cR}(E)\to \coDef _{\cQ}(E).$$ Indeed,
$\bR\phi ^!:D(\cA _{\cR})\to D(\cA _{\cQ})$ is an equivalence of
categories (Proposition 3.7) which commutes with the functor $\bR
i^!$.
\end{remark}

\begin{remark} A quasi-isomorphism $\delta:E_1\to E_2$ of
$\cA$-DG-modules induces an equivalence  of functors
$$\delta ^*:\coDef (E_2)\to \coDef (E_1)$$
by the formula $\delta ^*(S,\sigma)=(S,\sigma \cdot \delta)$.
\end{remark}

\begin{prop} Let $F:\cA \to \cA ^\prime$ be a DG functor as in Proposition 10.4 above.
 Consider the induced equivalence of derived
categories $\bR F^!:D(\cA )\to D(\cA ^\prime)$ (Corollary 3.15).
Then for any $E\in D(\cA )$ the deformation functors $\coDef (E)$
and $\coDef (\bR F^!(E))$ are canonically equivalent. (Hence also
$\coDef (F_*(E^\prime))$ and $\coDef (E^\prime)$ are equivalent for
any $E^\prime \in D(\cA ^\prime)$).
\end{prop}

\begin{proof} For any artinian DG algebra $\cR$ the functor $F$
induces a commutative functorial diagram
$$\begin{array}{ccc}
D(\cA _{\cR}) & \stackrel{\bR ((F \otimes \id)^!)}{\longrightarrow}
&
D(\cA _{\cR}^\prime)\\
\downarrow \bR i^! & & \downarrow \bR i^!\\
D(\cA ) & \stackrel{\bR F^!}{\longrightarrow} & D(\cA ^\prime),
\end{array}
$$
where $\bR (\cR \otimes F)^!$ is an equivalence by Corollary 3.15.
The horizontal arrows define a functor $ F^!_{\cR}:\coDef
_{\cR}(E)\to \coDef _{\cR}(\bR F^!(E))$. Moreover these functors are
compatible with the functors $\bR \phi ^!:\coDef _{\cR}\to \coDef
_{\cQ}$ induced by morphisms $\phi :\cR \to \cQ$ of artinian DG
algebras. So we get the functor $F^!:\coDef (E)\to \coDef (\bL
F(E))$. It is clear that for each $\cR$ the functor $F^!_{\cR}$ is
an equivalence. Thus $F^!$ is also such.
\end{proof}

\begin{example} Let $F:\cA ^\prime \to \cA $ be as in Example 10.5 above.
Then the assumption of the previous proposition holds.
\end{example}

\begin{remark} In the definition of the functor $\coDef (E)$ we
could work with the homotopy category of h-injective
DG modules instead of the derived category. Indeed, the functors
$i^!$ and $\phi ^!$ preserve h-injective DG modules.
\end{remark}

\begin{defi} Denote by $\coDef _+(E)$, $\coDef _-(E)$, $\coDef
_0(E)$, $\coDef _{\cl}(E)$ the restrictions of the functor $\coDef
(E)$ to subcategories $\dgart _+$, $\dgart _-$, $\art$, $\cart$
respectively.
\end{defi}

\begin{remark} The functors $\Def (E)$ and $\coDef (E)$ are not always equivalent
(unlike their homotopy counterparts $\Def ^{\h}(E)$ and $\coDef
^{\h}(E)$). In fact we expect that the functors $\Def$ and $\coDef $
are the "right ones" only in case they can be expressed in terms of
the functors $\Def ^{\h}$ and $\coDef ^{\h}$ respectively. (See the
next section).
\end{remark}

\section{Relation between functors $\Def$ and $\Def ^h$ (resp.
$\coDef $ and $\coDef ^h$)}

The ideal scheme that should relate these deformation functors is
the following. Let $\cA$ be a DG category, $E\in \cA ^0\text{-mod}$.
Choose quasi-isomorphisms $P\to E$ and $E\to I$, where $P\in
\P(\cA)$ and $I\in \cI(\cA)$. Then there should exist natural
equivalences
$$\Def(E)\simeq \Def ^h (P),\quad\quad \coDef(E)\simeq \coDef
^h(I).$$ Unfortunately, this does not always work.

\begin{example} Let $\cA$ be just a graded algebra $A=k[t]$, i.e.
$\cA$ contains a single object with the endomorphism algebra $k[t]$,
$\deg (t)=1$ (the differential is zero). Take the artinian DG
algebra $\cR$ to be $\cR=k[\epsilon]/(\epsilon ^2)$,
$\deg(\epsilon)=1$ (again the differential is zero). Let $E=A$ and
consider a DG $\cA _{\cR}^0$-module $M=E\otimes \cR$ with the
differential $d_M$ which is the multiplication by $\epsilon$.
Clearly, $M$ defines an object in $\Def ^h_{\cR}(E)$ which is not
isomorphic to the trivial deformation. However, one can check that
$\bL i^*M$ is not quasi-isomorphic to $E$ (although $i^*M=E$), thus
$M$ does not define an object in $\Def (E)$. This fact and the next
proposition show that the groupoid $\Def _{\cR}(E)$ is connected
(contains only the trivial deformation), so it is not the "right"
one.
\end{example}

\begin{prop} Assume that $\Ext ^{-1}(E,E)=0$.

1) Fix a quasi-isomorphism $P\to E$, $P\in \cP(\cA)$. Let $\cR $  be
an artinian DG algebra and $(S, \id)\in \Def ^{\h}(P)$. The
following conditions are equivalent:

a) $S\in \P(\cA _{\cR})$,

b) $i^*S=\bL i^*S$,

c) $(S,\id)$ defines an object in the groupoid $\Def _{\cR}(E)$.

The functor $\Def (E)$ is equivalent to the full subfunctor of
$\coDef ^{\h}(P)$ consisting of objects $(S,\id) \in  \Def
^{\h}(P)$, where $S$ satisfies a) (or b)) above.

2) Fix a quasi-isomorphism  $E\to I$ with $I\in \cI(\cA)$. Let $\cR
$  be an artinian DG algebra and $(T, \id)\in \coDef ^{\h}(I)$. The
following conditions are equivalent:

a') $T\in \cI(\cA _{\cR})$,

b') $i^!T=\bR i^!T$,

c') $(T,\id)$ defines an object in the groupoid $\coDef _{\cR}(E)$.

The functor $\Def (E)$ is equivalent to the full subfunctor of
$\coDef ^{\h}(I)$ consisting of objects $(T,\id) \in  \coDef
^{\h}(I)$, where $T$ satisfies a') (or b')) above.
  \end{prop}

\begin{proof} 1) It is clear that a) implies b) and b) implies c).
We will prove that c) implies a). We may and will replace the
functor $\Def (E)$ by an equivalent functor $\Def (P)$ (Remark
10.3).

Since $(S, \id)$ defines an object in $\Def _{\cR}(P)$ there exists
a quasi-isomorphism $g:\tilde{S}\to S$ where $\tilde{S}$ has
property (P) (hence $\tilde{S}\in \P(\cA _{\cR})$), such that $i^*g:
i^*\tilde{S}\to i^*S=P$ is also a quasi-isomorphism. Denote
$Z=i^*\tilde{S}$. Then $Z\in \P(\cA)$ and hence $i^*g$ is a homotopy
equivalence. Since both $\tilde{S}$ and $S$ are graded $\cR$-free,
the map $g$ is also a homotopy equivalence (Proposition 3.12d)).
Thus $S\in \P(\cA _{\cR})$.

Let us prove the last assertion in 1).

 Fix an object $(\overline{S} ,\tau) \in \Def _{\cR}(P)$. Replacing $(\overline{S},
\tau)$ by an isomorphic object we may and will assume that
$\overline{S}$ satisfies property (P). In particular,
$\overline{S}\in \cP(\cA _{\cR})$ and $\overline{S}$ is graded
$\cR$-free. This implies that $(\overline{S},\id)\in \Def^{\h}
_{\cR}(W)$ where $W=i^*\overline{S}$. We have $W\in \P(\cA)$. The
quasi-isomorphism $\tau :W\to S$ is therefore a homotopy
equivalence. By Corollary 8.4a) and Proposition 8.5a) there exists
an object $(S^\prime,\id)\in \Def ^{\h}_{\cR}(P)$ and a homotopy
equivalence $\tau ^\prime :\overline{S}\to S^\prime$ such that
$i^*(\tau ^\prime)=\tau$. This shows that $(\overline{S}, \tau)$ is
isomorphic (in $\Def _{\cR}(P)$) to an object $(S^\prime ,\id )\in
\Def ^{\h}_{\cR}(P)$, where $S^\prime \in \P(\cA _{\cR})$.

Let $(S ,\id ),(S^\prime, \id) \in \Def _{\cR}^{\h}(P)$ be two
objects such that $S,S^\prime \in \cP(\cA _{\cR})$.
 Consider the
obvious map
$$\delta :\Hom _{\Def ^{\h}_{\cR}(P)}((S,\id ),(S^\prime,\id ))\to
\Hom _{\Def _{\cR}(P)}((S,\id),(S^\prime,\id)).$$ It suffices to
show that $\delta$ is bijective.

Let $f:(S,\id ) \to (S^\prime, \id )$ be an isomorphism in $\Def
_{\cR}(P)$. Since $S,S^\prime\in \cP(\cA _{\cR})$ and $P \in
\cP(\cA)$ this isomorphism $f$ is a homotopy equivalence $f:S\to
S^\prime$ such that $i^*f$ is homotopic to $\id _P$. Let $h:i^*f\to
\id$ be a homotopy. Since $S$, $S^\prime$ are graded $\cR$-free the
map $i^*:\Hom (S, S^\prime )\to \Hom (P,P)$ is surjective
(Proposition 3.12a)). Choose a lift $\tilde{h}:S\to S^\prime[1]$ of
$h$ and replace $f$ by $\tilde{f}=f-d\tilde{h}$. Then
$i^*\tilde{f}=id$. Since $S$ and $S^\prime$ are graded $\cR$-free
$\tilde{f}$ is an isomorphism (Proposition 3.12d)). This shows that
$\delta$ is surjective.

Let $g_1,g_2:S\to S^\prime$ be two isomorphisms (in $\cA
_{\cR}^0\text{-mod}$) such that $i^*g_1=i^*g_2=\id _P$. That is
$g_1,g_2$ represent morphisms in $\Def _{\cR}^{\h}(P)$. Assume that
$\delta (g_1)=\delta (g_2)$, i.e. there exists a homotopy $s:g_1\to
g_2$. Then $d(i^*s)=i^*(ds)=0$. Since by our assumption $H^{-1}\Hom
(P,P)=0$ there exists $t\in \Hom ^{-2}(P,P)$ with $dt=i^*s$. Choose
a lift $\tilde{t}\in \Hom ^{-2}(S,S^\prime)$ of $t$. Then
$\tilde{s}:=s-d\tilde{t}$ is an allowable homotopy between $g_1$ and
$g_2$. This proves that $\delta $ is injective and finishes the
proof of 1).

The proof of 2) is very similar, but we present it for completeness.
Again it is clear that a') implies b') and b') implies c'). We will
prove that c') implies a') We may and will replace the functor
$\coDef (E)$ by an equivalent functor $\coDef (I)$ (Remark 10.10).

Since $(T,\id)$ defines an object in $\coDef _{\cR}(I)$, there
exists a quasi-isomorphism $g:T\to \tilde{T}$ where $\tilde{T}$ has
property (I) (hence $\tilde{T} \in \cI(\cA _{\cR})$), such that
$i^!g:I=i^!T\to i^!\tilde{T}$ is also a quasi-isomorphism. Denote
$K=i^!\tilde{T}$. Then $K\in \cI(\cA)$ and hence $i^!g$ is a
homotopy equivalence. Since both $T$ and $\tilde{T}$ are graded
$\cR$-cofree, the map $g$ is also a homotopy equivalence
(Proposition 3.12d)). Thus $T\in \cI(\cA _{\cR})$.

Let us prove the last assertion in 2).

 Fix an object $(\overline{T} ,\tau) \in \coDef _{\cR}(I)$.
 Replacing $(\overline{T}, \tau)$ by an isomorphic object we may and will
assume that $\overline{T}$ satisfies property (I). In particular,
$\overline{T}\in \cI(\cA _{\cR})$ and $\overline{T}$ is graded
$\cR$-cofree. This implies that $(\overline{T},\id)\in \coDef^{\h}
_{\cR}(L)$ where $L=i^!\overline{T}$. We have $L\in \cI (\cA)$ and
hence the quasi-isomorphism $\tau :I\to L$ is a homotopy
equivalence. By Corollary 8.4a) and Proposition 8.5a) there exist an
object $(T^\prime,\id)\in \coDef ^{\h}_{\cR}(I)$ and a homotopy
equivalence $\tau ^\prime : T^\prime\to \overline{T}$ such that
$i^!\tau ^\prime =\tau$. In particular, $T^\prime\in \cI(\cA
_{\cR})$. This shows that $(\overline{T},\tau)$ is isomorphic (in
$\coDef _{\cR}(I)$) to an object $(T^\prime ,\id) \in \coDef
_{\cR}^{\h}(I)$ where $T^\prime \in \cI(\cA _{\cR})$.

Let $(T ,\id ),(T^\prime, \id) \in \coDef _{\cR}^{\h}(I)$ be two
objects such that $T,T^\prime \in \cI(\cA _{\cR})$.
 Consider the
obvious map
$$\delta :\Hom _{\coDef ^{\h}_{\cR}(I)}((T,\id ),(T^\prime,\id ))\to
\Hom _{\coDef _{\cR}(I)}((T,\id),(T^\prime,\id)).$$ It suffices to
show that $\delta$ is bijective.

Let $f:(T,\id ) \to (T^\prime, \id )$ be an isomorphism in $\coDef
_{\cR}(I)$. Since $T,T^\prime\in \cI(\cA _{\cR})$ and $I \in
\cI(\cA)$ this isomorphism $f$ is a homotopy equivalence $f:T\to
T^\prime$ such that $i^!f$ is homotopic to $\id _I$. Let $h:i^!f\to
\id$ be a homotopy.  Since $T$, $T^\prime$ are graded $\cR$-cofree
the map $i^!:\Hom (T, T^\prime )\to \Hom (I,I)$ is surjective
(Proposition 3.12a)). Choose a lift $\tilde{h}:T\to T^\prime[1]$ of
$h$ and replace $f$ by $\tilde{f}=f-d\tilde{h}$. Then
$i^!\tilde{f}=id$. Since $T$ and $T^\prime$ are graded $\cR$-free
$\tilde{f}$ is an isomorphism (Proposition 3.12d)). This shows that
$\delta$ is surjective.

Let $g_1,g_2:T\to T^\prime$ be two isomorphisms (in $\cA
_{\cR}^0\text{-mod}$) such that $i^!g_1=i^!g_2=\id _I$. That is
$g_1,g_2$ represent morphisms in $\coDef _{\cR}^{\h}(I)$. Assume
that $\delta (g_1)=\delta (g_2)$, i.e. there exists a homotopy
$s:g_1\to g_2$. Then $d(i^!s)=i^!(ds)=0$. Since by our assumption
$H^{-1}\Hom (I,I)=0$ there exists $t\in \Hom ^{-2}(I,I)$ with
$dt=i^!s$. Choose a lift $\tilde{t}\in \Hom ^{-2}(T,T^\prime)$ of
$t$. Then $\tilde{s}:=s-d\tilde{t}$ is an allowable homotopy between
$g_1$ and $g_2$. This proves that $\delta $ is injective.
\end{proof}

\begin{remark} In the situation of Proposition 11.2
using Corollary 8.4b)  also obtain full and faithful embeddings of
functors $ \Def (E)$, $\coDef (E)$ in each of the equivalent
functors $\Def ^{\h}(P)$, $\coDef ^{\h}(P)$, $\Def ^{\h}(I)$,
$\coDef ^{\h}(I)$.
\end{remark}

The next two theorems provide important examples when the functors
$\Def _-$ and $\Def ^{\h}_-$ (resp. $\coDef _- $ and $\coDef
^{\h}_-$) are equivalent.

\begin{defi} An object $M\in \cA ^0\text{-mod}$ is called bounded above (resp. below) if
there exists $i$ such that $M(A)^j=0$ for all $A\in \cA$ and all $j\geq i$ (resp. $j\leq i$).
\end{defi}

\begin{theo} Assume  that $\Ext ^{-1}(E,E)=0$.

a) Suppose that there exists a bounded above $P\in \P(\cA)$  and a
quasi-isomorphism $P\to E$. Then the functors $\Def _-(E)$ and $\Def
_-^{\h}(P)$  are equivalent.

 b) Suppose that
there exists a bounded below $I\in \cI(\cA)$  and a
quasi-isomorphism $E\to I$. Then the functors $\coDef _-(E)$ and
$\coDef _-^{\h}(I)$  are equivalent.
\end{theo}

\begin{proof} Fix $\cR \in \dgart _-$ with the maximal ideal $m$.

In both cases it suffices to show that the embedding of groupoids
$\Def _{\cR} (E)\simeq \Def _{\cR}(P) \subset \Def _{\cR}^{\h}(P)$
(resp. $\coDef _{\cR} (E)\simeq \coDef _{\cR}(I) \subset \coDef
_{\cR}^{\h}(I)$) in Proposition 11.2 is essentially surjective.

a)  Choose an object $(S,\id)\in \Def _{\cR}^{\h}(P)$. It suffices
to prove the following lemma.

\begin{lemma}  The DG $\cA _{\cR}^0$-module $S$ is acyclic for the functor $i^*$, i.e.
$\bL i^*S=i^*S$.
\end{lemma}

Indeed the lemma implies that $S$ defines an object in $\Def
_{\cR}(P)$.

\begin{proof} Choose a quasi-isomorphism $f:Q\to S$ where $Q\in
\cP(\cA _{\cR})$. We need to prove that $i^*f$ is a
quasi-isomorphism. It suffices to prove that $\pi _!i^*f$ is a
quasi-isomorphism (Example 3.13). Recall that $\pi _!i^*=i^*\pi _!$.
Thus it suffices to prove that $\pi _!f$ is a homotopy equivalence.
Clearly $\pi _!f$ is a quasi-isomorphism. The DG $\cR^0$-module $\pi
_!Q$ is h-projective (Example 3.13). We claim that the DG $\cR
^0$-module $\pi _!S $ is also h-projective. Indeed, $\pi _!S$ is
bounded above and since $\cR \in \dgart _-$ this DG $\cR ^0$-module
has an increasing filtration with subquotients being free DG $\cR
^0$-modules. Thus $\pi _!S$ satisfies property (P) and hence is
h-projective.  It follows that the quasi-isomorphism $\pi _!f:\pi_!
Q\to \pi _!S$ is a homotopy equivalence. Hence $i^*\pi _!f=\pi
_!i^*f$ is also such.
\end{proof}

b) The following lemma implies that an object in $\coDef
_{\cR}^{\h}(I)$ is also an object in $\coDef _{\cR}(I)$, which
proves the theorem.
\end{proof}

\begin{lemma} Let $T\in \cA _{\cR}^0\text{-mod}$ be graded cofree
and bounded below. Then  $T$ is acyclic for the functor $i^!$, i.e.
$\bR i ^!T=i ^!T$.
\end{lemma}

\begin{proof} Denote $N=i^!T\in \cA ^0\text{-mod}$. Choose a quasi-isomorphism $g:T\to J$ where $J\in \cI(\cA
_{\cR})$. We need to prove that $i^!g$ is a quasi-isomorphism. It
suffices to show that $\pi _* i^!g$ is a quasi-isomorphism. Recall
that $\pi _*i^!=i^!\pi _*$. Thus it suffices to prove that $\pi _*g$
is a homotopy equivalence. Clearly it is a quasi-isomorphism.

Recall that the DG $\cR^0$-module $\pi _*J$ is h-injective (Example
3.13) We claim that $\pi _*T$ is also such. Indeed, since $\cR \in
\dgart _-$ the DG $\cR ^0$-module $\pi _*T$ has a decreasing
filtration
$$G^0\supset G^1\supset G^2\supset
...,$$ with
$$\gr \pi _*T=\oplus _{j}(\pi _*N)^j \otimes \cR ^*.$$
A direct sum of shifted copies of the DG $\cR ^0$-module $\cR ^*$ is
h-injective (Lemma 3.18). Thus each $(\pi _*N)^j \otimes \cR ^*$ is
h-injective and hence each quotient $\pi _*T/G^j$ is h-injective.
Also
$$\pi _*T=\lim_{\leftarrow}\pi _*T/G^j.$$
Therefore $\pi _*T$ is h-injective by Remark 3.5.

It follows that $\pi _*g$ is a homotopy equivalence, hence also
$i^!\pi _*g$ is such.
\end{proof}

The last theorem allows us to compare the functors $\Def _-$ and
$\coDef _-$ in some important special cases. Namely we have the
following corollary.

\begin{cor} Assume that

a) $\Ext^{-1}(E,E)=0$;

b) there exists a bounded above $P\in \cP(\cA)$ and a
quasi-isomorphism $P\to E$;

c) there exists a bounded below $I\in \cI(\cA)$ and a
quasi-isomorphism $E\to I$.

Then the functors $\Def _-(E)$  and $\coDef _-(E)$ are equivalent.
\end{cor}

\begin{proof} We have a quasi-isomorphism $P\to I$. Hence by Proposition 8.3 the DG algebras
$\End (P)$ and $\End (I)$ are quasi-isomoprhic. Therefore, in
particular, the functors $\Def _-^{\h}(P)$ and $\coDef _-^{\h}(I)$
are equivalent (Corollary 8.4b)). It remains to apply the last
theorem.
\end{proof}

In practice in order to find the required bounded resolutions one
might need to pass to a "smaller" DG category. So it is useful to
have the following stronger corollary.

\begin{cor} Let $F:\cA \to \cA ^\prime$ be a DG functor which induces a
quasi-equivalence $F^{\pre-tr}:\cA ^{\pre-tr}\to \cA ^{\prime
\pre-tr}$. Consider the corresponding equivalence $F_*:D(\cA
^\prime)\to D(\cA)$ (Corollary 3.15). Let $E\in \cA ^{\prime
0}\text{-mod}$ be such that

a) $\Ext^{-1}(E,E)=0$;

b) there exists a bounded above $P\in \cP(\cA)$ and a
quasi-isomorphism $P\to F_*(E)$;

c) there exists a bounded below $I\in \cI(\cA)$ and a
quasi-isomorphism $F_*(E)\to I$.

Then the functors $\Def _-(E)$  and $\coDef _-(E)$ are equivalent.
\end{cor}

\begin{proof} By the above corollary the functors
$\Def _-(F_*(E))$  and $\coDef _-(F_*(E))$ are equivalent. By
Proposition 10.4 the functors $\Def _-(E)$ and $\Def _-(F_*(E))$ are
equivalent. Since the functor $\bR F^!:D(\cA)\to D(\cA ^\prime)$ is
also an equivalence, we conclude that the functors $\coDef _-(E)$
and $\coDef _-(F_*(E))$ are equivalent by Proposition 10.11.
\end{proof}

\begin{example} If in the above corollary the DG category $\cA
^\prime$ is pre-triangulated, then one can take for $\cA$ a full DG
subcategory of $\cA ^\prime$ such that $\Ho (\cA ^\prime)$ is
generated as a triangulated category by the subcategory $\Ho (\cA)$.
One can often choose $\cA$ to have one object.
\end{example}

\begin{theo} Assume  that $\Ext ^{-1}(E,E)=0$. Let $E\to J$ be a
quasi-isomorphism with $J\in \cI(\cA )$. Suppose that there exists a
bounded below $P\in \P(\cA)$ and a quasi-isomorphism $P\to E$. Then

a) the functors $\coDef _-(E)$ and $\coDef _-^{\h}(J)$  are
equivalent;

b) the functors $\coDef _-(E)$ and $\coDef _-^{\h}(P)$  are
equivalent.
\end{theo}

\begin{proof} We may replace the functor $\coDef _-(E)$ by
an equivalent functor $\coDef _-(J)$ or $\coDef _-(P)$.

Consider the quasi-isomorphism $f:P\to J$. Given $(T,id)\in \coDef
^{\h}_{\cR}(J)$ by Corollary 8.4b) and Proposition 8.5b) there
exists $(S,\id)\in \coDef ^{\h}_{\cR}(P)$ and a quasi-isomorphism
$\tilde{f}:S\to T$ which lifts $f$, i.e. $i^!\tilde{f}=f$. By Lemma
11.7  we have $\bR i^!S=i^!S$. Hence also $\bR i^!T=i^!T$. Therefore
the functors $\coDef _-(J)$ and $\coDef _-^{\h}(J)$ can be
identified by part 2) of Proposition 11.2. This proves a). Now b)
follows from a) and the canonical equivalence of functors $\coDef
_-^{\h}(P)\simeq \coDef _-^{\h}(J)$ (Corollary 7.5b)).
\end{proof}

\subsection{Relation between functors $\Def _-(E)$, $\coDef
_-(E)$ and $\Def _-(\cC)$, $\coDef _-(\cC)$}

The next proposition  follows immediately from our previous results.

\begin{prop} Let $\cA$ be a DG category and $E\in \cA ^0\text{-mod}$. Assume that

a) $\Ext ^{-1}(E,E)=0$;

b) there exists a quasi-isomorphism $F\to E$ with $F\in \cP (\cA)$
or a quasi-isomorphism $E\to F$ with $F\in \cI(\cA)$ such that $F$
is bounded below (resp. there exists a quasi-isomorphism $F\to E$
with $F\in \cP(\cA)$ bounded above);

c) there exists a bounded below (resp. bounded above) DG algebra
$\cC$ which is quasi-isomorphic to $\End (F)$.

Then the functors $\coDef _-(E)$ and $\coDef _-(\cC)$ (resp.
functors $\Def _-(E)$ and $\Def _- (\cC)$) are equivalent.
\end{prop}

\begin{proof} By Theorem 11.5, 11.11b) above we have the following equivalences
of functors $\coDef _-(E)\simeq \coDef _-^{\h}(F)$, $\coDef
_-(\cC)\simeq \coDef _-^{\h}(\cC)$ (resp. $\Def _-(E)\simeq \Def
_-^{\h}(F)$, $\Def _-(\cC)\simeq \Def _-^{\h}(\cC)$). By  Example
6.3 and Theorem 8.1 there are the equivalences $\coDef
^{\h}_-(\cC)\simeq \coDef ^{\h}_-(F)$, $\Def ^{\h}_-(\cC)\simeq \Def
^{\h}_-(F)$. This proves the proposition.
\end{proof}

\begin{remark}  The equivalences of functors
$\coDef ^{\h}_-(\cC)\simeq \coDef ^{\h}_-(F)$, $\Def
^{\h}_-(\cC)\simeq \Def ^{\h}_-(F)$ in the end of the proof of last
proposition can be made explicit. Put $\cB =\End(F)$. Assume, for
example, that $\psi :\cC \to \cB$ is a homomorphism of DG algebras
which is a quasi-isomorphism. Then the composition of DG functors
(Propositions 9.2,9.4)
$$\Sigma ^F\cdot \psi ^*:\cC ^0\text{-mod}\to \cA ^0\text{-mod}$$
induces equivalences of functors
$$\Def ^{\h} (\Sigma ^F \cdot \psi ^*): \Def ^{\h}(\cC)\simeq \Def
^{\h}(F)$$
$$\coDef ^{\h} (\Sigma ^F \cdot \psi ^*):\coDef ^{\h}(\cC)\simeq \coDef
^{\h}(F)$$  by Propositions 9.2e) and 9.4f).
\end{remark}

\part{Pro-representability}

\section{Coalgebras and the bar construction}

\subsection{Coalgebras and comodules} We will consider DG
coalgebras. For a DG coalgebra $\cG$ we denote by $\cG ^{\gr}$ the
corresponding graded coalgebra obtained from $\cG$ by forgetting the
differential. Recall that if $\cG$ is a DG coalgebra, then its
graded dual $\cG ^*$ is naturally a DG algebra. Also given a finite
dimensional DG algebra $\cB$ its dual $\cB ^*$ is a DG coalgebra.

A morphism of DG coalgebras $k\to \cG$ (resp. $\cG \to k$) is called
a co-augmentation (resp. a co-unit) of $\cG$ if it satisfies some
obvious compatibility condition. We denote by $\overline{\cG}$ the
cokernel of the co-augmentation map.

Denote by $\overline{\cG} _{[n]}$ the kernel of the $n$-th iterate
of the co-multiplication map $\Delta ^n:\overline{\cG} \to
\overline{\cG} ^{\otimes n}$. The DG coalgebra $\cG$ is called {\it
co-complete} if
$$\overline{\cG}=\bigcup_{n\geq 2}\overline{\cG} _{[n]}.$$

A $\cG$-{\it comodule} means a left DG comodule over $\cG$.

 A  $\cG^{\gr}$-comodule is {\it free} if it is isomorphic to
$\cG\otimes V$ with the obvious comodule structure for some graded
vector space $V$.

Denote by $\cG ^0$ the DG coalgebra with the opposite
co-multiplication.

Let $g:\cH \to \cG$ be a homomophism of DG coalgebras. Then $\cH$ is
a DG $\cG$-comodule with the co-action $g\otimes 1\cdot \Delta
_{\cH}:\cH\to \cG \otimes \cH$ and a DG $\cG ^0$-comodule with the
co-action  $1\otimes g\cdot \Delta _{\cH}:\cH\to \cH \otimes \cG$.

Let $M$ and $N$ be a right and left DG $\cG$-comodules respectively.
Their cotensor product $M\square _{\cG}N$  is defined as the kernel
of the map
$$\Delta _M\otimes 1-1\otimes \Delta _N:M\otimes N\to M\otimes \cG
\otimes N,$$ where $\Delta _M:M\to M\otimes \cG$ and $\Delta _N:N\to
\cG \otimes N$ are the co-action maps.

A DG coalgebra $\cG$ is a left and right DG comodule over itself.
Given a DG $\cG$-comodule $M$ the co-action morphism $M\to \cG
\otimes M$ induces an isomorphism $M=M\square _{\cG}\cG$. Similarly
for DG $\cG ^0$-modules.

\begin{defi} The dual $\cR ^*$ of an artinian DG algebra $\cR$ is
called an {\it artinian} DG coalgebra.
\end{defi}

Given an artinian DG algebra $\cR$, its augmentation $\cR \to k$
induces the co-augmentation $k\to \cR ^*$ and  its unit $k\to \cR$
induces the co-unit $\cR ^*\to k$.

\subsection{From comodules to modules}

If   $P$ is a DG comodule over a DG coalgebra $\cG$, then $P$ is
naturally a DG module over the DG algebra $(\cG ^*)^0$. Namely, the
$(\cG ^*)^0$-module structure is defined as the composition
$$P\otimes \cG^*\stackrel{\Delta_P\otimes 1}{\longrightarrow}\cG \otimes P\otimes \cG ^*
\stackrel{T\otimes 1}{\longrightarrow}P\otimes \cG \otimes \cG
^*\stackrel{1 \otimes \ev}{\longrightarrow}P,$$ where $T:\cG\otimes
P\to P\otimes \cG$ is the transposition map.

 Similarly, if $Q$ is
a DG  $\cG ^0$-comodule, then $P$ is a DG module over $\cG ^*$.

Let $P$ and $Q$ be a left and right DG $\cG$-comodules respectively.
Then $P\otimes Q$ is a  DG $\cG ^*$-bimodule, i.e. a DG $\cG^*
\otimes \cG ^{*0}$-module by the above construction. Note that its
center
$$Z(P\otimes Q):=\{ x\in P\otimes Q\ \vert \ ax=(-1)^{\bar{a}\bar{x}}xa
\ \  \text{for all} \ \  a\in \cG ^*\}$$ is isomorphic to the
cotensor product $Q\square _{\cG}P$.

\subsection{Twisting cochains and the bar construction}

Let $\cG$ be a DG coalgebra and $\cC$ be a DG algebra. Recall that
the complex $\Hom _k(\cG,\cC)$ is a DG algebra under the
convolution. That is for $f,g\in \Hom _k(\cG ,\cC)$ their product is
the composition
$$\cG \stackrel{\Delta }{\to} \cG \otimes \cG
\stackrel{f\otimes g}{\to} \cC \otimes \cC \stackrel{\mu }{\to}
\cC,$$ where $\Delta$ and $\mu$ denote co-multiplication and
multiplication in $\cG $ and $\cC$ respectively.

\begin{defi} The elements of the Maurer-Cartan cone $MC(\Hom
_k(\cG ,\cC))$ are called twisting cochains. Suppose that $\cC$ is
augmented and $\cG$ is co-augmented. Then a twisting cochain $\tau
\in MC(\Hom _k(\cG ,\cC))$ is called admissible if it comes from a
twisting cochain in $\Hom _k(\overline{\cG},\overline{\cC})$.
\end{defi}

\begin{example} Suppose in the above definition that
$\cG$ is concentrated in nonnegative degrees. Then a twisting
cochain in $MC(\Hom _k(\cG ,\cC))$ which comes from a twisting
cochain in $\Hom _k(\overline{\cG},\cC)$ is automatically admissible
(for degree reasons).
\end{example}

Let $\tau \in MC(\Hom _k(\cG ,\cC))$. If $f:\cH \to \cG$ is a
homomorphism of DG coalgebras (resp. $g:\cC \to \cD$ is a
homomorphism of DG algebras), then $\tau  f \in MC(\Hom _k(\cH
,\cC))$ (resp. $g \tau \in MC(\Hom _k(\cG ,\cD))$).

Note that the identity map $\id :\Hom _k(\cG ,\cC)\to \Hom _k(\cG
^0,\cC ^0)$ is an isomorphism of the DG algebra $\Hom _k(\cG ,\cC)$
with the DG algebra $\Hom _k(\cG ^0,\cC ^0)^0$. Hence $\tau \in
MC(\Hom _k(\cG ,\cC))$ if and only if $-\tau \in MC(\Hom _k(\cG
^0,\cC ^0))$.

\begin{defi} Let $\cC$ be an augmented DG algebra with the augmentation ideal
$\overline{\cC}$. We denote by $B\cC$ its bar construction, which is
a DG coalgebra. Recall that as a coalgebra $B\cC$ is the tensor
coalgebra on the graded vector space $\overline{\cC}[1]$.   Denote
by $\hat{S}=(B\cC)^*$ the dual DG algebra.
\end{defi}

The DG coalgebra $B\cC$ is co-augmented and co-complete. Hence the
DG algebra $\hat{S}$ is an augmented complete local DG algebra.
$B\cC$ is the union of its finite dimensional DG subcoalgebras.
Hence $\hat{S}$ is the inverse limit of its finite dimensional
quotients. Note that any finite dimensional quotient of $\hat{S}$ is
an artinian DG algebra, hence $\hat{S}$ is a pro-object in the
category $\dgart$. It follows that any finite dimensional DG
subcoalgebra of $B\cC$ is artinian (i.e. its dual is an artinian DG
algebra).

The following lemma is from [Le].

\begin{lemma} Let $\cC$ be an augmented DG algebra. Then the composition $\tau _{\cC}$ of the
projection with $B\cC\to \overline{\cC}[1]$ with the (shifted)
embedding $\overline{\cC}\hookrightarrow \cC$ is the {\it }
universal admissible twisting cochain for $\cC$. That is given a
co-augmented DG coalgebra $\cG$ and an admissible twisting cochain
$\tau :\cG \to \cC$ there exists a unique morphism $g_{\tau}:\cG \to
B\cC$ of DG coalgebras such that $\tau _{\cC}\cdot g_{\tau}=\tau.$
\end{lemma}

\subsection{From modules to comodules}

Let $\tau \in MC(\Hom _k(\cG ,\cC))$ be a twisting cochain and $M$
be a DG $\cC$-module. Consider the free $\cG^{\gr}$-comodule
$\cG\otimes M$ with the differential $d_{\tau}=d_{\cG}\otimes
1+1\otimes d_M+t_{\tau}$, where $t_{\tau}$ is the composition
$$\cG\otimes M\stackrel{\Delta\otimes 1 }{\longrightarrow}\cG\otimes
\cG\otimes M\stackrel{1\otimes \tau \otimes
1}{\longrightarrow}\cG\otimes \cC\otimes M\stackrel{1\otimes \mu
_M}{\longrightarrow} \cG\otimes M.$$ Here $\Delta$ is the
comultiplication in $\cG$ and $\mu _M$ is the $\cC$-module structure
of $M$. Then $d_{\tau }^2=0$ and $\cG\otimes M$ with the
differential $d_{\tau}$ is a DG $\cG$-comodule. Hence by the
construction in (12.2) it is a DG $(\cG ^*)^0$-module. To stress the
role of $\tau$ we will sometimes denote this DG module by $\cG
\otimes _{\tau}M$.

Similarly, if $N$ is a DG $\cC ^0$-module, then using $-\tau\in
MC(\Hom _k(\cG ,\cC))$ we obtain a structure of a DG $\cG
^0$-comodule (hence of a DG $\cG ^*$-module by 12.1) on $N\otimes
\cG$. We denote the resulting DG module by $N\otimes _{\tau}\cG$.

\begin{example} Put $M=C$. Then the $\cG$-comodule (or the $(\cG ^*)^0$-module) structure on
$\cG \otimes _{\tau}\cC$ commutes with the right multiplication on
$\cC$. Thus $\cG \otimes _{\tau}\cC$ is a DG $(\cC \otimes \cG
^*)^0$-module. Similarly, $\cC \otimes _{\tau}\cG$ is a DG $\cC
\otimes \cG ^*$-module.
\end{example}

Let $f:\cH\to \cG$ be a homomorphism of DG coalgebras and $\tau \in
MC(\Hom _k(\cG ,\cC))$. Then there are natural isomorphisms of DG
$\cC \otimes \cH ^*$- and $(\cC \otimes \cH ^*)^0$-modules
respectively
$$\cC \otimes _{\tau f}\cH \simeq (\cC \otimes _{\tau}\cG) \square
_{\cG}\cH,\quad \cH \otimes _{\tau f}\cC \simeq \cH \square
_{\cG}(\cG\otimes _{\tau}\cC).$$

\subsection{The bar complex}

\begin{defi} Let $\cC$ be an augmented DG algebra, $\cG=B\cC$ and $\tau
_{\cC}:B\cC\to \cC$ be the canonical (universal) twisting cochain.
The resulting DG  $(\cC \otimes \hat{S})^0$-module $B\cC\otimes
\cC=B\cC \otimes _{\tau _{\cC}}\cC$ is called the bar complex of
$\cC$. It is quasi-isomorphic to $k$. Similarly, if we consider
$\cC$ as a DG $\cC^0$-module then we obtain a DG $\cC \otimes
\hat{S}$-module $\cC \otimes _{\tau _{\cC}}B\cC$, which is also
quasi-isomorphic to $k$.
\end{defi}

\subsection{Deformations and the bar complex}

\begin{lemma} Let $\cC$ be a DG algebra and $\cR$ a finite dimensional DG
algebra. Then the DG algebras $\cC \otimes \cR$ and $\Hom _k(\cR ^*
,\cC)$ are naturally isomorphic.
\end{lemma}

\begin{proof} The isomorphism $\theta :\cC \otimes \cR\to
\Hom_k(\cR ^* ,\cC)$ is given by the formula
$$\theta (c\otimes r)(f)=
(-1)^{\bar{f}\bar{r}}cf(r)$$ for $f\in \cR ^*$, $c\in \cC$, $r\in
\cR $.
\end{proof}

\begin{remark} Assume that in the above lemma the DG algebra $\cR$ is
artinian with the maximal ideal $m\subset \cR$. Note that under the
isomorphism of this lemma the ideal $\cC \otimes m$ is mapped to the
ideal $\Hom _k(m^*, \cC)$. Thus in particular there is a natural
bijection between the set $MC(\cC \otimes m)$ and the collection of
twisting cochains in $\Hom _k(m^*,\cC)$.
\end{remark}

\begin{remark} Assume in the above remark that the DG algebra $\cC$ is augmented and that  $\cR \in
\dgart _-$ (hence $\cR ^*$ is concentrated in nonnegative degrees).
Then elements of  $MC(\cC \otimes m)$ are in bijection with the set
of admissible twisting cochains in  $\Hom _k(\cR ^*,\cC)$.
\end{remark}

\subsubsection{} Let $\cC$ be an augmented DG algebra, $\cR$ an
artinian DG algebra with the maximal ideal $m$. Let $\tau \in \Hom
_k(\cR ^* ,\cC)$ be an admissible twisting cochain and $\cR
^*\otimes _{\tau}\cC$ be the corresponding DG $(\cC \otimes
\cR)^0$-module.

Let $\alpha \in MC(\cC\otimes m)$ be the element corresponding to
$\tau$ by Remark 12.9 and denote by $\cC\otimes _{\alpha}\cR ^*$ the
corresponding object in the groupoid $\coDef ^{\h}_{\cR}(\cC)$. Thus
in particular $\cC\otimes _{\alpha}\cR ^*$ is a DG $(\cC \otimes
\cR)^0$-module.

\begin{lemma} The DG $(\cR \otimes
\cC)^0$-modules $\cR ^*\otimes _{\tau}\cC$ and $\cC\otimes
_{\alpha}\cR ^*$ are isomorphic.
\end{lemma}

\begin{proof} The isomorphism is simply the transposition of the two
factors.
\end{proof}

In the above notation consider the homomorphism of DG coalgebras
$g_{\tau}:\cR ^* \to B\cC$ induced by $\tau$ as in Lemma 12.3. This
induces a homomorphism of the dual DG algebras
$g_{\tau}^*:\hat{S}\to \cR$. Thus in particular $\cR$ becomes a DG
$\hat{S}^0$-module. Notice that this $\hat{S}^0$-module structure on
$\cR$ coincides with the one coming from the homomorphism
$g_{\tau}:\cR ^* \to B\cC$ by the construction in 12.2 and the
identification $\cR =\cR ^{**}$.

Consider the bar complex $B\cC \otimes \cC$ (which is a DG
$(\cC\otimes \hat{S})^0$-module) and the DG $(\cC \otimes
\cR)^0$-module $\Hom _{\hat{S}^0}(\cR ,B\cC\otimes\cC).$

\begin{lemma} The DG $(\cC
\otimes \cR)^0$-modules $\Hom _{\hat{S}^0}(\cR ,B\cC\otimes\cC)$ and
$\cR ^*\otimes _{\tau}\cC$ are isomorphic.
\end{lemma}

\begin{proof} Indeed, by 12.4 and 12.2 we have the isomorphisms of
DG $(\cC \otimes \cR)^0$-modules
$$\cR ^*\otimes _{\tau}\cC =\cR ^* \square _{B\cC}(B\cC \otimes
\cC)=Z(\cR^*\otimes (B\cC \otimes \cC)),$$ where $Z(\cR^*\otimes
(B\cC \otimes \cC))$ is the center of the $\hat{S}\otimes
\hat{S}^0$-module $\cR^*\otimes (B\cC \otimes \cC)$. Now, given a DG
$(\cC \otimes \hat{S})^0$-module $M$, there is an isomorphism of DG
$(\cC \otimes \cR)^0$-modules
$$\theta :Z(\cR ^* \otimes M)\to \Hom _{\hat{S}^0}(\cR ,M)$$ defined
by $\theta (f\otimes m)(r)=(-1)^{\bar{r}\bar{m}}f(r)m$.
\end{proof}

Later on we will use the following corollary.

\begin{cor} Let $\cC$ be an augmented DG algebra. Then for every object $(T ,\id)\in \coDef _{\cR}^{\h}(\cC)$
there exists a homomorphism of DG algebra $\hat{S}\to \cR$ and an
isomorphism of DG $(\cC\otimes \cR )^0$-modules $T=\Hom
_{\hat{S}^0}(\cR ,B\cC \otimes \cC)$.
\end{cor}

\begin{proof} This follows from the last two lemmas.
\end{proof}

\section{Some functors defined by the bar complex}

\begin{defi} An augmented DG algebra $\cC$ is called admissible if

a) it is
 {\it nonnegative}, i.e $\cC ^i=0$
for $i<0$;

b) it is {\it connected}, i.e $\cC ^0=k$;

c) and it is locally finite, i.e. $\dim _k\cC ^i  <\infty$ for all
$i$.
\end{defi}

\subsection{The functor $\Delta$}
Fix an augmented DG algebra $\cC$. Consider the bar construction
$B\cC$, the corresponding DG algebra $\hat{S}=(B\cC)^*$ and the DG
$(\hat{S}\otimes \cC)^0$-module $B\cC \otimes \cC$=$B\cC \otimes
_{\tau _{\cC}} \cC$ (the bar complex). If $\cC$ is connected and
nonnegative, then $B\cC$ is concentrated in nonnegative degrees and
consequently $\hat{S}$ is concentrated in nonpositive degrees.

Denote by $D_f(\hat{S})\subset D(\hat{S})$ the full triangulated
subcategory consisting of DG modules with finite dimensional
cohomology.

\begin{lemma} Assume that the DG algebra $\cC$ is connected and nonnegative.
Then the category $D_f(\hat{S})$ is the triangulated envelope of the
DG $\hat{S^0}$-module $k$.
\end{lemma}

\begin{proof} Denote by $\langle k\rangle\subset D(\hat{S})$ the
triangulated envelope of $k$. By our assumption the DG algebra
$\hat{S}$ is concentrated in nonpositive degrees.

 Let $M$ be a DG $\hat{S}^0$-module with finite
dimensional cohomology. First assume that $M$ is concentrated in one
degree. Then $\dim M<\infty$. Since $\hat{S}^{\gr}$ is a local
algebra the module $M$ has a filtration with subquotients isomorphic
to $k$. Thus $M\in \langle k\rangle$.

In the general case by Lemma 3.19 we may and will assume that
$M^i=0$ for $\vert i\vert
>>0$. Let $s$ be the least integer such that $M^s\neq 0$. The kernel
$K$ of the differential $d:M^s\to M^{s+1}$ is a DG
$\hat{S}^0$-submodule. By the above argument $K\in \langle
k\rangle$. If $K\neq 0$ then by induction on the dimension of the
cohomology we obtain that $M/K \in \langle k\rangle$. Hence also $M
\in \langle k\rangle$. If $K=0$, then the DG $\hat{S}^0$-submodule
$\tau _{< s+1}M$ (Lemma 3.19) is acyclic, and hence $M$ is
quasi-isomorphic to $\tau _{\geq s+1}M$. But we may assume that
$\tau _{\geq s+1}M \in \langle k\rangle$ by descending induction on
$s$.
\end{proof}

Choose a quasi-isomorphism of DG $\hat{S}^0$-modules $B\cC \otimes
\cC \to J$, where $J$ satisfies the property (I) (hence is
h-injective).

Consider the contravariant DG functor $\Delta
:\hat{S}^0\text{-mod}\to \cC ^0\text{-mod}$ defined by
$$\Delta (M):=\Hom _{\hat{S}^0}(M,J)$$
This functor extends trivially to derived categories $\Delta
:D(\hat{S})\to D(\cC )$.

\begin{theo} Assume that the DG algebra $\cC$ is admissible. Then

 a) The contravariant
functor $\Delta $ is full and faithful on the category
$D_f(\hat{S})$.

b) $\Delta(k)$ is isomorphic to $\cC$.
\end{theo}

\begin{proof} By Lemma 13.2 the category $D_f(\hat{S})$ is the triangulated
envelope of the DG $\hat{S}^0$-module $k$. So for the first
statement of the theorem it suffices to prove that the map $\Delta
:\Ext _{\hat{S}^0}(k,k)\to \Ext _{\cC ^0}(\Delta (k),\Delta (k))$ is
an isomorphism. The following proposition implies the theorem.

\begin{prop} Under the assumptions of the above theorem the following holds.

a) The complex $\bR \Hom _{\hat{S}^0}(k,k)$ is quasi-isomorphic to
$\cC$.

b) The natural morphism of complexes $\Hom _{\hat{S}^0}(k,B\cC
\otimes \cC)\to \Hom _{\hat{S}^0}(k,J)$ is a quasi-isomorphism.

c) $\Delta (k)$ is quasi-isomorphic to $\cC$.

d) $\Delta :\Ext _{\hat{S}^0}(k,k)\to \Ext _{\cC ^0}(\Delta
(k),\Delta (k))$ is an anti-isomorphism.
\end{prop}

\begin{proof} a) Recall the DG $\hat{S}\otimes \cC$-module
$\cC \otimes B\cC=\cC \otimes _{\tau _{\cC}}B\cC$ (Definition 12.7).
Consider the corresponding DG $(\hat{S}\otimes \cC)^0$-module
$P:=\Hom _k(\cC \otimes B\cC, k)$. Since $\cC$ is locally finite and
bounded below and $B\cC$ is bounded below the graded
$\hat{S}^0$-module $P^{\gr}$ is isomorphic to $(\hat{S}\otimes \Hom
_k(\cC ,k))^{\gr}$. Since the complex $\Hom _k(\cC ,k)$ is bounded
above and the DG algebra $\hat{S}$ is concentrated in nonnegative
degrees the DG $\hat{S}^0$-module $P$ has the property (P) (and
hence is h-projective). Thus $\bR \Hom _{\hat{S}^0}(k,k)=\Hom
_{\hat{S}^0}(P,k)=\Hom _k(\Hom _k(\cC ,k),k)=\cC ^0$. This proves
a).

b): Since $\Hom _{\hat{S}^0}(k,B\cC \otimes \cC)=\cC$ the assertion
follows from a).

c) follows from b).

d) follows from a) and c).
\end{proof}

This proves the theorem.
\end{proof}

\begin{remark}
Notice that for any augmented DG algebra $\cC$ we have $\Hom
_{\hat{S}^0}(k, B\cC \otimes \cC)=\cC$. Thus the DG $(\cC\otimes
\hat{S})^0$-module is a "homotopy $\hat{S}$-co-deformation" of
$\cC$. The Proposition 13.4 implies that for an admissible  $\cC$
this DG $(\cC\otimes \hat{S})^0$-module is a "derived
$\hat{S}$-co-deformation" of $\cC$.
 (Of course we have only defined co-deformations along artinian DG
algebras.)
\end{remark}

\subsection{The functor $\nabla$}

Now we define another functor $\nabla :D(\hat{S})\to D(\cC)$, which
is closely related to $\Delta$.

Denote by $m$ the augmentation ideal of $\hat{S}$. For a DG
$\hat{S}^0$-module $M$ denote $M_n:=M/m^nM$ and
$$\hat{M}=\lim_{\stackrel{\longleftarrow}{n}}M_n.$$
 Fix a DG
$\hat{S}^0$-module $N$. Choose a quasi-isomorphism $P\to N$ with an
h-projective $P$. Define
$$\nabla
(N):=\lim_{\rightarrow}\Delta(P_n)=\lim_{\rightarrow}\Hom_{\hat{S}^0}(P_n,J).$$

Denote by $\Perf(\hat{S})\subset D(\hat{S})$ the full triangulated
subcategory which is generated by the DG $\hat{S}^0$-module
$\hat{S}$.

\begin{theo}
Assume that the DG algebra $\cC$ is admissible and finite
dimensional. Then

a) The contravariant functor $\nabla :D(\hat{S})\to D(\cC)$ is full
and faithful on the subcategory $\Perf(\hat{S})$.

b) $\nabla (\hat{S})$ is isomorphic to $k$.
\end{theo}

\begin{proof}
Denote by $m\subset \hat{S}^0$ the maximal ideal and put
$S_n:=\hat{S}^0/m^n\hat{S}^0$. Since the DG algebra $\cC$ is finite
dimensional $S_n$ is also finite dimensional for all $n$. We need a
few lemmas.

\begin{lemma} Let $K$ be a DG $\hat{S}^0$-module such that $\dim
_kK<\infty$. Then the natural morphism of complexes
$$\Hom _{\hat{S}^0}(K,B\cC \otimes \cC)\to \Hom _{\hat{S}^0}(K,J)$$
is a quasi-isomorphism.
\end{lemma}

\begin{proof} Notice that since the algebra $\hat{S}$ is local,
every element $x\in m$ acts on $K$ as a nilpotent operator. Hence in
particular $m^nK=0$ for $n>>0$. For the same reason the DG
$\hat{S}^0$-module $K$ has a filtration with subquotients isomorphic
to $k$. Thus we may prove the assertion by induction on $\dim K$. If
$K=k$, then this is part b) of Proposition 13.4. Otherwise we can
find a short exact sequence of DG $\hat{S}^0$-modules
$$0\to M\to K\to N\to 0,$$
such that $\dim M,\dim N <\dim K$.

\medskip

\noindent{\bf Sublemma.} {\it The sequence of complexes
$$0\to \Hom _{\hat{S}^0}(N, B\cC \otimes \cC)\to \Hom _{\hat{S}^0}(K, B\cC \otimes
\cC) \to \Hom _{\hat{S}^0}(M, B\cC \otimes \cC) \to 0$$ is  exact.}

\medskip

\begin{proof} We only need to prove the surjectivity of the map
$$\Hom _{\hat{S}^0}(K, B\cC \otimes
\cC) \to \Hom _{\hat{S}^0}(M, B\cC \otimes \cC).$$

  Let
$n>>0$ be such that $m^nK=m^nM=0$. Let ${}_n(B\cC \otimes
\cC)\subset (B\cC \otimes \cC)$ denote the DG $\hat{S}^0$-submodule
consisting of elements $x$ such that $m^nx=0$. Then ${}_n(B\cC
\otimes \cC)$
 is a DG $S_n$-module and $\Hom _{\hat{S}^0}(K, B\cC \otimes
 \cC)=\Hom_{S_n}(K,{}_n(B\cC \otimes \cC))$ and similarly for $M$.

Note that ${}_n(B\cC \otimes \cC)$
 as a graded $S_n$-module is isomorphic to $S^*_n\otimes \cC$,
 hence is a finite direct sum of shifted copies of the injective graded module
 $S_n ^*$.  Hence the above map of complexes is surjective.
 \end{proof}

 Now we can prove the lemma.

 Consider the commutative diagram of complexes
 $$\begin{array}{ccccccccc} 0 & \to & \Hom _{\hat{S}^0}(N, B\cC \otimes \cC) & \to &  \Hom _{\hat{S}^0}(K, B\cC \otimes
\cC) & \to & \Hom _{\hat{S}^0}(M, B\cC \otimes \cC) &  \to &  0\\
 & & \downarrow \alpha & & \downarrow \beta & & \downarrow \gamma &
 & \\
 0 & \to &  \Hom _{\hat{S}^0}(N, J) & \to & \Hom _{\hat{S}^0}(K, J) & \to & \Hom _{\hat{S}^0}(M, J) & \to &
 0,
 \end{array}$$
where the bottom row is exact since $J^{gr}$ is an injective graded
$\hat{S}^0$-module (because $J$ satisfies property (I)). By the
induction assumption $\alpha $ and $\gamma $ are quasi-isomorphisms.
Hence also $\beta $ is such.
\end{proof}

We are ready to prove the theorem.

It follows from Lemma 13.6 that $\nabla(\hat{S})$ is
quasi-isomorphic to
$$\lim_{\rightarrow}\Hom
_{\hat{S}^0}(S_n,B\cC\otimes \cC)=\lim_{\rightarrow}\Hom
_{S_n}(S_n,{}_n(B\cC\otimes
\cC))=\lim_{\rightarrow}({}_n(B\cC\otimes \cC))=B\cC\otimes \cC.$$
This proves the second assertion. The first one follows from the
next lemma.

\begin{lemma} For any augmented DG algebra $\cC$
the complex $\bR \Hom _{\cC ^0}(k,k)$ is quasi-isomorphic to $\hat{S}$.
\end{lemma}

\begin{proof}  Note that the
DG $\cC ^0$-module $B\cC \otimes \cC$ has the property (P). Hence
$$\bR\Hom _{\cC ^0}(k,k)=\Hom _{\cC ^0}(B\cC \otimes \cC ,k)=\Hom _k
(B\cC ,k)=\hat{S}.$$
\end{proof}

This proves the theorem.
\end{proof}

\subsection{The functor $\Psi$}

Finally consider the covariant functor $\Psi :D(\hat{S}^0)\to
D(\cC)$ defined by
$$\Psi (M):=(B\cC \otimes \cC )\stackrel{\bL}{\otimes
}_{\hat{S}}M.$$

\begin{theo} For any augmented DG algebra $\cC$ the following holds.

a) The functor $\Psi$ is full and faithful on the subcategory
$\Perf(\hat{S})$.

b) $\Psi(\hat{S})=k$
\end{theo}

\begin{proof} b) is obvious and a) follows from Lemma 13.8 above.
\end{proof}

\section{Pro-representability of the co-deformation functor}

\subsection{The 2-category $2\text{-}\dgalg$ and deformation functor $\coDEF$}

\begin{defi} We define the 2-category $2\text{-}\dgalg$ of DG algebras as follows.
The objects are augmented DG algebras. For DG algebras $\cB, \cC$
the collection of 1-morphisms $1\text{-}\Hom(\cB,\cC)$ consists of
pairs $(M,\theta)$, where
\begin{itemize}
\item $M\in D(\cB \otimes \cC^0)$ is such that there exists an
isomorphism (in $D(\cC^0)$)  $\cC\to \nu _*M$ (where $\nu _*:D(\cB
\otimes \cC^0)\to D(\cC^0)$ is the functor of restriction of scalars
corresponding to the natural homomorphism $\nu :\cC  \to \cB ^0
\otimes \cC$);
\item and $\theta :k{\stackrel{\bL}{\otimes }}_{\cC}M\to k$ is an isomorphism in
$D(\cB)$.
\end{itemize}
The  composition of 1-morphisms
$$1\text{-}\Hom(\cB,\cC)\times 1\text{-}\Hom(\cC,\cD)\to
1\text{-}\Hom(\cB,\cD)$$
 is defined by the tensor product $\cdot
 \stackrel{\bL}{\otimes}_{\cC}\cdot $.
 Given
1-morphisms $(M_1,\theta _1), (M_2,\theta _2)\in
1\text{-}\Hom(\cB,\cC)$ a 2-morphism $f: (M_1,\theta _1)\to
(M_2,\theta _2)$ is an isomorphism (in $D(\cB \otimes \cC^0)$)
$f:M_2\to M_1$ (not from $M_1$ to $M_2$!) such that $\theta _1\cdot
k{\stackrel{\bL}{\otimes }}_{\cC}(f)=\theta _2$. So in particular
the category $1\text{-}\Hom(\cB ,\cC)$ is a groupoid. Denote by
$2\text{-}\dgart$ the full 2-subcategory of $2\text{-}\dgalg$
consisting of artinian DG algebras. Similarly we define the full
2-subcategories $2\text{-}\dgart_+$, $2\text{-}\dgart_-$,
$2\text{-}\art$, $2\text{-}\cart$ (Definition 2.3).
\end{defi}

\begin{remark} Assume that augmented DG algebras $\cB$ and $\cC$ are
such that $\cB ^i=\cC ^i=0$ for $i>0$ and $\dim \cB ^i, \dim \cC
^i<\infty$ for all $i$. Denote by $\langle k\rangle \subset D(\cB
\otimes \cC ^0)$ the triangulated envelope of the DG $\cB ^0\otimes
\cC$-module $k$. Let $(M,\theta )\in 1\text{-}\Hom (\cB ,\cC)$. Then
by Corollary 3.22 $M\in \langle k\rangle$.
\end{remark}

For any DG algebra $\cB$ we obtain a 2-functor $h_{\cB}$ between the
2-categories $2\text{-}\dgalg$ and {\bf Gpd} defined by
$h_{\cB}(\cC)=1\text{-}\Hom (\cB ,\cC)$.

Note that a usual homomorphism of augmented
 DG algebras $\gamma :\cB
\to \cC$ defines the structure of a DG $\cB^0$-module on $\cC$ with
the canonical isomorphism of DG $\cB^0$-modules $\id :
k{\stackrel{\bL}{\otimes }}_{\cC}\cC \to k$. Thus it defines a
1-morphism $(\cC ,\id)\in 1\text{-}\Hom (\cB ,\cC)$. This way we get
a 2-functor $\cF:\adgalg \to 2\text{-}\dgalg$, which is the identity
on objects.

\begin{lemma} Assume that augmented DG algebras $\cB$ and $\cC$ are concentrated in degree zero
(hence have zero differential). Then

a) the map $\cF:\Hom (\cB ,\cC)\to \pi _0(1\text{-}\Hom (\cB,\cC))$
is surjective, i.e. every 1-morphism from $\cB$ to $\cC$ is
isomorphic to $\cF(\gamma)$ for a homomorphism of algebras $\gamma$;

b) the 1-morphisms $\cF(\gamma _1)$ and $\cF(\gamma _2)$ are
isomorphic if and only if $\gamma _2$ is the composition of $\gamma
_1$ with the conjugation by an invertible element in $\cC$;

c) in particular, if $\cC$ is commutative then the map of sets
$\cF:\Hom (\cB ,\cC)\to \pi _0(1\text{-}\Hom (\cB,\cC))$ is a
bijection.
\end{lemma}

\begin{proof} a) For any  $(M, \theta)\in 1\text{-}\Hom (\cB,\cC)$ the DG
$\cB ^0 \otimes \cC$-module $M$ is isomorphic (in $D(\cB  \otimes
\cC^0)$ to  $H^0(M)$. Thus we may assume that $M$ is concentrated in
degree 0. By assumption there exists an isomorphism of $\cC$-modules
 $\cC \to M$. Multiplying this isomorphism by a scalar we may assume that
 it is compatible with the isomorphisms $\id :k{\stackrel{\bL}{\otimes }}_{\cC}\cC\to k$
 and $\theta :k{\stackrel{\bL}{\otimes }}_{\cC}M\to k$.
  A choice of such an isomorphism defines a homomorphism of algebras
 $\cB ^0\to \End _{\cC}(\cC)=\cC ^0$.  Thus
$(M,\theta)$ is isomorphism to $\cF (\gamma)$.

b) Let $\gamma _1,\gamma _2:\cB \to \cC$ be homomorphisms of
algebras. A 2-morphism $f:\cF (\gamma _1)\to \cF (\gamma _2)$ is
simply an isomorphism of the corresponding $\cB ^0\otimes \cC
$-modules $f:\cC \to \cC$, which commutes with the augmentation.
Being an isomorphism of $\cC$-modules it is the right multiplication
by an invertible element $c\in \cC$. Hence for every $b \in \cB$ we
have $c^{-1}\gamma _1(b)c=\gamma _2(b)$.

c) This follows from a) and b).
\end{proof}

\begin{remark} If in the definition of 1-morphisms $1\text{-}\Hom(\cB,\cC)$ we do not fix an
isomorphism $\theta$, then we obtain a special case of a
"quasi-functor" between the DG categories $\cB \text{-mod}$
 and $\cC \text{-mod}$. This notion was first introduced by Keller
 in [Ke] for DG modules over general DG categories.
 \end{remark}

The next proposition asserts that the deformation functor $\coDef$
has a natural "lift" to the 2-category $2\text{-}\dgart$.

\begin{prop} There exist a 2-functor  $\coDEF(E)$ from
$2\text{-}\dgart$ to {\bf Grp} and an equivalence of 2-functors
$\coDef(E)\simeq \coDEF(E)\cdot \cF$. For an artinian DG algebra
$\cR$ the groupoids  $\coDef _{\cR}(E)$ and $\coDEF _{\cR}(E)$) are
isomorphic.
\end{prop}

\begin{proof} Given artinian DG
algebras $\cR$, $\cQ$  and $M=(M,\theta )\in 1\text{-}\Hom (\cR
,\cQ)$ we need to define the corresponding functor
$$M^!:\coDef _{\cR}(E)\to \coDef _{\cQ}(E).$$
Let $S=(S,\sigma)\in \coDef _{\cR}(E)$. Put
$$M^!(S):=\bR \Hom _{\cR^0}(M,S)\in D(\cA _{\cQ}).$$
We claim that $M^!(S)$ defines an object in $\coDef _{\cQ}(E)$, i.e.
$\cR \Hom _{\cQ^0}(k,M^!(S))$ is naturally isomorphic to $E$ (by the
isomorphisms $\theta $ and $\sigma$).

 Indeed, choose quasi-isomorphisms $P\to k$ and $S\to I$
for $P\in \cP(\cA _{\cQ})$ and $I\in \cI(\cA _{\cR})$. Then
$$\bR\Hom _{\cQ^0}(k,M^!(S))=\Hom _{\cQ^0}(P,\Hom _{\cR^0}(M,I)).$$
By Lemma 3.17 the last term is equal to $\Hom _{\cR}(P\otimes
_{\cQ}M,I)$. Now the isomorphism $\theta$ defines an isomorphism
between $P\otimes _{\cQ}M=k\stackrel{\bL}{\otimes} _{\cQ}M$ and $k$,
and we compose it with the isomorphism $\sigma :E \to \bR\Hom
_{\cR}(k,I)=i^!S$.

So $M^!$ is a functor from $\coDef_{\cR}(E)$ to $\coDef _{\cQ}(E)$.

Given another artinian DG algebra $\cQ ^\prime$ and $M^\prime \in
1\text{-}\Hom (\cQ ,\cQ ^\prime)$ there is a natural isomorphism of
functors
$$(M^\prime
\stackrel{\bL}{\otimes}_{\cQ}M)^!(-) \simeq M^{\prime !}\cdot M
^!(-).$$ (This follows again from Lemma 3.17).

Also a 2-morphism $f\in 2\text{-}\Hom(M,M_1)$ between objects
$M,M_1\in 1\text{-}\Hom (\cR ,\cQ)$ induces an isomorphism of the
corresponding functors $M^!\stackrel{\sim}{\to}M_1^!$.

Thus we obtain a 2-functor $\coDEF (E):2\text{-}\dgart \to {\bf
Gpd}$, such that $\coDEF (E)\cdot \cF=\coDef(E)$.
\end{proof}

We denote by $\coDEF_+(E)$, $\coDEF_-(E)$, $\coDEF_{\art}(E)$,
$\coDEF_{\cl}(E)$ the restriction of the functor $\coDEF (E)$ to
subcategories $2\text{-}\dgart _+$, $2\text{-}\dgart _-$,
$2\text{-}\art$ and $2\text{-}\cart$ respectively.

The next result is the analogue of Proposition 11.12 for the functor
$\coDEF _-$.

\begin{prop} Let $E\in \cA ^0\text{-mod}$. Assume that

a) $\Ext ^{i}(E,E)=0$ for $i<0$ and $\Ext ^0(E,E)=k$;

b) there exists a quasi-isomorphism $F\to E$ with $F\in \cP (\cA)$
or a quasi-isomorphism $E\to F$ with $F\in \cI(\cA)$ such that $F$
is bounded below.

Then for any DG algebra $\cC$ which is bounded below and
quasi-isomorphic to $\End(F)$ the functors $\coDEF _-(\cC)$ and
$\coDEF _-(E)$ are equivalent.
\end{prop}

\begin{proof} If DG algebras $\cC$ and $\End(F)$ are
quasi-isomorphic, then there exists a DG algebra $\cD$ and
homomorphisms of DG algebras $\phi :\cD \to \cC$ and $\psi :\cD \to
\End(F)$ which are quasi-isomorphisms. By Lemma 9.5 we may assume
that $\cD$ is also bounded below. We will prove that the
homomorphism $\psi$ induces an equivalence of functors $\coDEF
_-(\cD)$ and $\coDEF _-(F)$. The same argument proves that the
homomorphism $\phi$ induces an equivalence of functors $\coDEF
_-(\cD)$ and $\coDEF _-(\cC)$. Thus we may and will assume that
$\cC=\cD$.

 Consider the DG functor
$$\cL :=\Sigma ^F\cdot \psi ^*:\cC ^0\text{-mod}\to \cA
^0\text{-mod}, \quad \cL(N)=N\otimes _{\cC}F$$ as in Remark 11.13.
It induces the equivalence of functors
$$\coDef ^{\h}(\cL):\coDef ^{\h}_-(\cC)\stackrel{\sim}{\to}\coDef
^{\h}_-(F).$$ For every artinian DG algebra $\cR \in \dgart _-$ the
corresponding DG functor
$$\cL _{\cR}:(\cC\otimes \cR )^0\text{-mod}\to \cA
_{\cR}^0\text{-mod}$$ induces the equivalence of groupoids
$\coDef^{\h} _{\cR}(\cC)\stackrel{\sim}{\to}\coDef^{\h}_{\cR}(F)$
(Propositions 9.2,9.4). By Theorems 11.5b) and 11.11b) there are
natural equivalences of functors $$\coDef ^{\h}_-(F)\simeq \coDef
_-(E), \quad \coDef ^{\h}_-(\cC)\simeq \coDef _-(\cC).$$ Hence the
functor $\bL\cL$ induces the equivalence
$$\bL\cL:\coDef _-(\cC)\stackrel{\sim}{\to}\coDef _-(E).$$

Fix $\cR,\cQ \in 2\text{-}\dgart _-$ and $M\in
1\text{-}\Hom(\cR,\cQ)$. We need to show that there exists an
isomorphism between functors from $\coDef _{\cR}(\cC)$ to $\coDef
_{\cQ}(E)$
$$\bL\cL _{\cQ}\cdot M^!\simeq M^!\cdot \bL\cL _{\cR}.$$

Since the cohomology of $M$ is finite dimensional, and the DG
algebra $\cR \otimes \cQ$ has no components in positive degrees, by
Corollary 3.21 we may assume that $M$ is finite dimensional.

\begin{lemma} Let $(S, \id)$ be an object in $\coDef
^{\h}_{\cR}(\cC)$ or in $\coDef ^{\h}_{\cR}(F)$. Then $S$ is acyclic
for the functor $M^!$, i.e. $M^!(S)=\Hom _{\cR ^0}(M,S)$.
\end{lemma}

\begin{proof} In the proof of Lemma 11.7 we showed that $S$ is
h-injective when considered as a DG $\cR^0$-module.
\end{proof}

Choose $(S,\id)\in \coDef ^{\h}(\cC)$. By the above lemma
$M^!(S)=\Hom _{\cR^0}(M,S)$.

We claim that the DG $\cC^0$-module $\Hom _{\cR^0}(M,S)$ is
h-projective. Indeed, first notice that the graded $\cR ^0$-module
$S$ is injective being isomorphic to a direct sum of copies of
shifted graded $\cR^0$-module $\cR ^*$ (the abelian category of
graded $\cR ^0$-modules is locally notherian, hence a direct sum of
injectives is injective). Second, the DG $\cR ^0$-module $M$ has a
(finite) filtration with subquotients isomorphic to $k$. Thus the DG
$\cC ^0$-module $\Hom _{\cR ^0}(M,S)$ has a filtration with
subquotients isomorphic to $\Hom _{\cR ^0}(k,S)=i^!S\simeq \cC$. So
it has property (P).

Hence $\bL\cL\cdot M^!(S)=\Hom _{\cR^0}(M,S)\otimes _{\cC}F$. For
the same reasons $M^!\cdot \bL\cL _{\cR}(S)=\Hom _{\cR^0}(P,
S\otimes _{\cC}F)$. The isomorphism
$$\Hom_{\cR^0}(M,S)\otimes _{\cC}F=\Hom _{\cR^0}(M, S\otimes _{\cC}F)$$
follows from the fact that $S$ as a graded module is a tensor
product of graded $\cC^0$ and $\cR^0$ modules and also because $\dim
_kM<\infty$.
\end{proof}

\subsection{The 2-category $2^\prime\text{-}\dgalg$ and deformation
functor $\DEF$}

It turns out that the deformation functor $\Def$ lifts naturally to
a different version of a 2-category of DG algebras. We denote this
2-category $2^\prime\text{-}\dgalg$. It differs from
$2\text{-}\dgalg$ in two respects: the 1-morphisms are objects in
$D(\cB ^0\otimes \cC)$ (instead of $D(\cB \otimes \cC^0)$) and
2-morphisms go in the opposite direction.

\begin{defi} We define the 2-category $2^\prime \text{-}\dgalg$ of DG algebras as follows.
The objects are DG algebras. For DG algebras $\cB, \cC$ the
collection of 1-morphisms $1\text{-}\Hom(\cB,\cC)$ consists of pairs
$(M,\theta)$, where
\begin{itemize}
\item $M\in D(\cB ^0 \otimes \cC)$ and there exists an isomorphism (in
$D(\cC)$)  $\theta :\cC\to \nu _*M$ (where $\nu _*:D(\cB^0 \otimes
\cC)\to D(\cC)$ is the functor of restriction of scalars
corresponding to the natural homomorphism $\nu :\cC ^0  \to \cB
\otimes \cC ^0$);
\item and $\theta :M{\stackrel{\bL}{\otimes }}_{\cC}k\to k$ is an isomorphism
in $D(\cB^0)$.
\end{itemize}
The composition of 1-morphisms
$$1\text{-}\Hom(\cB,\cC)\times 1\text{-}\Hom(\cC,\cD)\to
1\text{-}\Hom(\cB,\cD)$$
 is defined by the tensor product $\cdot
 \stackrel{\bL}{\otimes}_{\cC}\cdot $.
 Given
1-morphisms $(M_1,\theta _1), (M_2,\theta _2)\in
1\text{-}\Hom(\cB,\cC)$ a 2-morphism $f: (M_1,\theta _1)\to
(M_2,\theta _2)$ is an isomorphism (in $D(\cB ^0 \otimes \cC)$)
$f:M_1\to M_2$ such that $\theta _1=\theta _2\cdot
(f){\stackrel{\bL}{\otimes}}_{\cC}k$. So in particular the category
$1\text{-}\Hom(\cB ,\cC)$ is a groupoid. Denote by
$2^\prime\text{-}\dgart$ the full 2-subcategory of
$2^\prime\text{-}\dgalg$ consisting of artinian DG algebras.
Similarly we define the full 2-subcategories
$2^\prime\text{-}\dgart_+$, $2^\prime\text{-}\dgart_-$,
$2^\prime\text{-}\art$, $2^\prime\text{-}\cart$ (Definition 2.3).
\end{defi}

\begin{remark} The exact analogue of Remark 14.2 holds for the
category $2^\prime\text{-}\dgalg$.
\end{remark}

For any DG algebra $\cB$ we obtain a 2-functor $h_{\cB}$ between the
2-categories $2^\prime\text{-}\dgalg$ and {\bf Gpd} defined by
$h_{\cB}(\cC)=1\text{-}\Hom (\cB ,\cC)$.

Note that a usual homomorphism of DG algebras $\gamma :\cB \to \cC$
defines the structure of a $\cB$-module on $\cC$ with the canonical
isomorphism of DG $\cB$-modules
$\cC{\stackrel{\bL}{\otimes}}_{\cC}k$. Thus it defines a 1-morphism
$(\cC ,\id)\in 1\text{-}\Hom (\cB ,\cC)$. This way we get a
2-functor $\cF^\prime:\adgalg \to 2^\prime\text{-}\dgalg$, which is
the identity on objects.

\begin{remark}
The precise analogue of Lemma 14.3 holds for the 2-category
$2^\prime \text{-}\dgalg$ and the functor $\cF ^\prime$.
\end{remark}

\begin{prop} There exist a 2-functor  $\DEF(E)$ from
$2^\prime\text{-}\dgart$ to {\bf Grp} and an equivalence of
2-functors $\Def(E)\simeq \DEF(E)\cdot \cF ^\prime$. For an artinian
DG algebra $\cR$ the groupoids  $\Def _{\cR}(E)$ and $\DEF
_{\cR}(E)$ are isomorphic.
\end{prop}

\begin{proof} Let $\cR$,
$\cQ$ be artinian DG algebras. Given $(M,\theta )\in 1\text{-}\Hom
(\cR ,\cQ)$ we define the corresponding functor
$$M^*: \Def _{\cR}(E)\to \Def _{\cQ}(E)$$ as follows
$$M^*(S):=S\stackrel{\bL}{\otimes}_{\cR}M$$
for $(S,\sigma)\in \Def _{\cR}(E)$. Then we have the canonical
isomorphism
$$M^*(S)\stackrel{\bL}{\otimes }_{\cQ}k=
S\stackrel{\bL}{\otimes }_{\cR}(M\stackrel{\bL}{\otimes }_{\cQ}k)
\stackrel{\theta}{\to} S\stackrel{\bL}{\otimes
}_{\cR}k\stackrel{\sigma}{\to}E.$$
 So that $M^*(S)\in \Def
_{\cQ}(E)$ indeed.

Given another artinian DG algebra $\cQ ^\prime$ and $M^\prime \in
1\text{-}\Hom (Q,Q^\prime)$ there is a natural isomorphism of
functors
$$M^{\prime *}\cdot M^*=(M\stackrel{\bL}{\otimes
}_{\cQ}M^\prime)^*.$$ Also a 2-morphism $f\in 2\text{-}\Hom (M,M_1)$
between $M,M_1\in 1\text{-}\Hom (\cR ,\cQ)$ induces an isomorphism
of corresponding functors $M^*\stackrel{\sim}{\to}M^*_1$.

Thus we obtain a 2-functor $\DEF(E):2^\prime \text{-}\dgart \to {\bf
Gpd}$, such that $\DEF(E)\cdot \cF ^\prime =\Def(E)$.
\end{proof}

We denote by $\DEF_+(E)$, $\DEF_-(E)$, $\DEF_{\art}(E)$,
$\DEF_{\cl}(E)$ the restriction of the functor $\DEF (E)$ to
subcategories $2^\prime\text{-}\dgart _+$, $2^\prime\text{-}\dgart
_-$, $2^\prime\text{-}\art$ and $2^\prime\text{-}\cart$
respectively.

The next proposition is the analogue of Proposition 14.6 for the
functor $\DEF _-$.

\begin{prop} Let $E\in \cA ^0\text{-mod}$. Assume that

a) $\Ext ^{-1}(E,E)=0$;

b) there exists a quasi-isomorphism  $F\to E$ with $F\in \cP(\cA)$
bounded above;

Then a homomorphism of  DG algebras $\psi :\cC \to \End(F)$, such
that $\cC$ is bounded above and $\psi$ is a quasi-isomorphism,
induces an equivalence of functors $\DEF _-(E)$ and $\DEF _-
(\cC)$).
\end{prop}

\begin{proof} The proof is exactly parallel to the proof of
Proposition 13.6. We omit it.
\end{proof}

\subsection{Pro-representability theorem}

\begin{theo}  Let $\cC$ be an admissible DG algebra (Definition 13.1).
Then the functor
$\coDEF _-(\cC)$ is pro-representable by the DG algebra
$\hat{S}=(B\cC)^*$. That is there exists an equivalence of functors
$\coDEF _-(\cC)\simeq h_{\hat{S}}$ from $2\text{-}\dgart_-$ to {\bf
Gpd}.
\end{theo}

As a corollary we obtain the following theorem.

\begin{theo} Fix $E\in \cA ^0\text{-mod}$.
Assume that $E$ is quasi-isomorphic to a bounded below $F$ which is
h-projective or h-injective. Also assume that an admissible DG
algebra $\cC$  is quasi-isomorphic to $\End(F)$. Then the functor
$\coDEF _-(E)$ is pro-representable by the DG algebra
$\hat{S}=(B\cC)^*$.
\end{theo}

\begin{proof} Indeed, by Proposition 14.6 the functors $\coDEF
_-(E)$ and $\coDEF _-(\cC)$ are equivalent, so it remains to apply
Theorem 14.13.
\end{proof}

\begin{proof} Let us prove Theorem 14.13.

Consider the DG $(\cC\otimes \hat{S})^0$-module $B\cC \otimes \cC$.
Choose a quasi-isomorphism $B\cC\otimes \cC \to J$, where $J$ is an
h-injective DG $(\cC\otimes \hat{S})^0$-module. By Lemma 3.23 $J$ is
also h-injective when considered as a DG $\hat{S}^0$- or $\cC
^0$-module via the restriction of scalars.

We first define a morphism of functors $\Theta :\h_{\hat{S}}\to
\coDEF_-(\cC)$. Namely, given an artinian DG algebra $\cR \in
2\text{-}\dgart _-$ and a 1-morphism $M=(M,\theta) \in 1\text{-}\Hom
(\hat{S},\cR)$ we define
$$\Theta (M):=\Hom _{\hat{S}^0}(M,J).$$
By Lemma 3.17 we have $\bR\Hom _{\cR^0}(k,\Hom
_{\hat{S}^0}(M,J))=\bR\Hom
_{\hat{S}^0}(k{\stackrel{\bL}{\otimes}}_{\cR}M,J)$. Hence the
quasi-isomorphism $\theta :k{\stackrel{\bL}{\otimes}} _{\cR}M\to k$
induces a quasi-isomorphism
$$\bR\Hom _{\cR^0}(k,\Theta (M))\simeq
\bR\Hom _{\hat{S}^0}(k, J)=\Hom _{\hat{S}^0}(k, J),$$ and by
Proposition 13.4 the last term is canonically quasi-isomorphic to
$\cC$ as a DG module over $\cC^0$. Thus we obtain a canonical
quasi-isomorphism $\sigma (M) :\bR\Hom _{\cR^0}(k,\Theta (M))\to
\cC$, which means that $(\Theta (M), \sigma (M))\in \coDEF
_{\cR}(\cC)$.

Given another artinian DG algebra $Q\in 2\text{-}\dgart _-$ and a
1-morphism $N=(N,\delta)\in 1\text{-}\Hom (\cR ,\cQ)$ we claim that
the object $\Theta (N{\stackrel{\bL}{\otimes }}_{\cR}M)\in \coDEF
_{\cQ}(\cC)$ is canonically isomorphic to the object $\bR\Hom
_{\cR^0}(N,\Theta (M))$. This follows again from Lemma 3.17. Thus
$\Theta$ is indeed a morphism of functors.

It remains to prove that for each $\cR \in 2\text{-}\dgart _-$ the
induced functor $\Theta _{\cR}: 1\text{-}\Hom (\hat{S},\cR)\to
\coDEF _{\cR}(\cC)$ is an equivalence of groupoids. So fix a DG
algebra $\cR \in 2\text{-}\dgart _-$.

\medskip

\noindent{\bf Surjective on isomorphism classes.} By Theorem 11.11b)
we know that groupoids $\coDEF _{\cR}(\cC)=\coDef _{\cR}(\cC)$ and
$\coDef ^{\h} _{\cR}(\cC)$ are equivalent. Given an object $(T ,\id
)\in \coDef ^{\h}_{\cR}(\cC)$ it is known by Corollary 12.13 that
there exists a homomorphism of DG algebras $\phi: \hat{S}\to \cR$
such that the DG $(\cR\otimes \cC)^0$-module $T$ is isomorphic to
$\Hom _{\hat{S}^0}(\cR ,B\cC \otimes \cC)$. It follows from Lemma
13.7  that $\Hom _{\hat{S}^0}(\cR ,B\cC \otimes \cC)=\bR\Hom
_{\hat{S}^0}(\cR ,B\cC \otimes \cC)$. Therefore $(T,\id)$,
considered as an object in $\coDEF _{\cR}(\cC)$ is isomorphic to
$\Theta (M)$, where $M=\cR$ is a DG module over $\cR\otimes
\hat{S}^0$ via the homomorphism $\phi$.

\medskip

\noindent{\bf Full and faithful.}

Consider the above functor $\Theta$ as a contravariant DG functor
from $\cR \otimes \hat{S}^0\text{-mod}$ to $(\cC \otimes
\cR)^0\text{-mod}$ given by
$$\Theta (M)=\Hom _{\hat{S}^0}(M,J).$$

Define the contravariant DG functor $\Phi :(\cC \otimes
\cR)^0\text{-mod}\to \cR \otimes \hat{S}^0\text{-mod}$ by the
similar formula
$$\Phi (N)=\Hom _{\cC ^0}(N,J).$$
These DG functors induce the corresponding functors between the
derived categories
$$\Theta :D(\cR ^0\otimes \hat{S})\to D(\cC \otimes \cR),\quad \quad
\Phi :D(\cC \otimes \cR ) \to D(\cR ^0\otimes \hat{S}).$$ Denote by
$\langle k\rangle\subset D(\cR ^0\otimes \hat{S})$ and $\langle \cC
\rangle \subset D(\cC \otimes \cR)$ the triangulated envelopes of
the DG $\cR\otimes \hat{S}^0$-module $k$ and the DG $(\cC \otimes
\cR)^0$-module $\cC$ respectively.

\begin{lemma} The functors $\Theta $ and $\Phi$ induce mutually
inverse anti-equivalences of the triangulated categories $\langle
k\rangle$ and $\langle \cC \rangle$.
\end{lemma}

\begin{proof} For $M\in \cR \otimes \hat{S}^0\text{-mod}$ and $N\in
(\cC \otimes \cR)^0\text{-mod}$ we have the functorial morphisms of
DG modules
$$\beta _M: M\to \Phi (\Theta (M)), \quad \beta _M(x)(f)=(-1)^{\bar{f}\bar{x}}f(x),$$
$$\gamma _N: N\to \Theta (\Phi (N)), \quad \gamma _N(y)(g)=(-1)^{\bar{g}\bar{y}}g(y).$$

By Proposition 13.4a) the DG $(\cC \otimes \cR)^0$-module $\Theta
(k)$ is quasi-isomorphic to $\cC$ and hence $\Phi (\Theta (k))$ is
quasi-isomorphic to $k$, so that $\beta _k$ is a quasi-isomorphism.
Also we have $\Phi (\cC)=J$ and again by Proposition 13.4a) $\Theta
(J)$ is quasi-isomorphic to $\cC$. So $\gamma _{\cC}$ is a
quasi-isomorphism. This proves the lemma.
\end{proof}

Notice that for $(M,\theta )\in 1\text{-}\Hom (\hat{S},\cR)$ (resp.
for $(S,\sigma)\in \coDEF _{\cR}(\cC)=\coDef ^{\h}_{\cR}(\cC)$)
$M\in \langle k\rangle$ by Remark 14.2 (resp. $S\in \langle
\cC\rangle$). Hence the fact that the functor $\Theta _{\cR}:
1\text{-}\Hom (\hat{S},\cR)\to \coDEF _{\cR}(\cC)$ is full and
faithful follows from the last lemma and the following remark.

\begin{remark} The following diagram of DG functors is commutative
$$
\begin{array}{rcl}
\cR \otimes \hat{S}^0\text{-mod} &
\stackrel{\Theta}{\longrightarrow} & (\cC \otimes
\cR)^0\text{-mod}\\
k\otimes _{\cR}(\cdot)\downarrow & & \downarrow \Hom _{\cR
^0}(k,\cdot)\\
\hat{S}^0\text{-mod} & \stackrel{\Theta}{\longrightarrow} &
\cC^0\text{-mod}.
\end{array}
$$
Indeed, this follows from Lemma 3.17.
\end{remark}

This proves the theorem.
\end{proof}

\end{document}